\documentclass[12pt,spanish]{amsart}

\usepackage{amsaddr}
\usepackage{amsmath,amsthm, amsfonts,amssymb}
\usepackage[spanish]{babel}
\usepackage{appendix} 
\usepackage{cancel}
\usepackage[all]{xy}
\usepackage[colorinlistoftodos,shadow]{todonotes}
\usepackage{setspace}
\usepackage[utf8]{inputenc}
\usepackage{tikz}
\usetikzlibrary{mindmap,trees}
\usepackage{verbatim}
\usepackage{nag}
\usepackage{hyperref}
\usepackage{tikz-cd}
\usepackage{stmaryrd}

\newcommand{\bb}[1]{\mathbb{#1}}
\newcommand{\fk}[1]{\mathfrak{#1}}

\newcommand{\N}{\bb N}

\newcommand{\set}[1]{\left\{{#1}\right\}}

\newcommand{\sqrd}[8]{
 \xymatrix{%
 {#1} \ar[r]^-{#5}\ar[d]_-{#8} & {#2}\ar[d]^-{#6} \\ %
 {#4}\ar[r]_-{#7} & {#3} %
 }%
 } 

\def\Spec{\mathop{\mbox{\normalfont Spec}}\nolimits}

\newtheorem{Theorem}{Teorema}[section]
\newtheorem{prop}[Theorem]{Proposición}
\newtheorem{lem}[Theorem]{Lema}
\newtheorem{cor}[Theorem]{Corolario}
\theoremstyle{definition}
    \newtheorem{Def}[Theorem]{Definición}
    \newtheorem{exa}[Theorem]{Ejemplo}
    
\theoremstyle{remark}
    \newtheorem{obs}[Theorem]{Observación}

\begin{document}


\author{J. Rogelio Pérez-Buendía}
\thanks{Quiero agradecer al Profesor Arthur Ogus, de la Universidad de California en Berkeley, por haberme facilitado una copia de su libro ``Lectures on Logarithmic Algebraic Geometry'', que fue de gran ayuda para este trabajo.}
\address{CONACYT-CIMAT Mérida}
\email{rogelio.perez@cimat.mx}

\author{Ernesto Antonio Reyes Ramírez}
\address{DEMAT-Universidad de Guanajuato}
\email{ernesto.reyes@cimat.mx}
\title{La Geometría de los Monoides}


\date{Julio 2020}
\maketitle
\begin{abstract}
    En este artículo  presentamos las bases de la teoría de monoides desde el punto de vista categórico, haciendo énfasis en las analogías y diferencias entre la teoría de módulos sobre anillos conmutativos. Se presenta la generalización de esquema afín a esquema afín monoidal; estudiamos la relación que éstos tienen con las variedades tóricas y sentamos las bases para el estudio de la Geometría Logarítmica de Fontaine-Kato-Illusie~\cite{kato1989logarithmic}, ampliamente usada en Geometría Aritmética. 
\end{abstract}
\tableofcontents

\section{Introducción}
Un monoide consiste en un conjunto M, una operación binaria asociativa $m:M\times M \to M$ y un elemento identidad $e_{M}$. Si además la operación binaria es conmutativa llamaremos a M un monoide conmutativo. Claramente todo grupo es un monoide por lo que el desarrollo de esta teoría tendrá consigo muchas analogías con la teoría de grupos o en su caso con la teoría de anillos. 

De igual forma que en la teoría de anillos conmutativos, en la teoría de monoides conmutativos se definen sus ideales e ideales primos los cuales dan paso a sus espectros. A dichos espectros, como en el caso de anillos, los podemos dotar de una topología, la topología de Zariski, y así desarrollar equivalencias entre propiedades geométricas y algebraicas.

Aún con todas estas analogías entre distintas estructuras algebraicas existen diferencias específicas las cuales hacen importante el estudio de los monoides conmutativos. Una de las diferencias más grandes es la carencia de un análogo del lema de Nakayama ya que en monoides conmutativos no necesariamente se tienen leyes de cancelación. Mientras que en la teoría de anillos tenemos el primer teorema de isomorfismos, no existe un análogo de este para monoides. Como veremos el perder ciertas propiedades nos hacen ganar otras muy importantes y en cierto sentido más generales.


El objetivo principal de este trabajo es el de sentar las bases para el estudio de la Geometría Logarítmica de Fontaine-Kato-Illusie como un acercamiento a la Geometría Aritmética. La Geometría Logarítmica descansa fuertemente en el uso de gavillas de monoides de ahí que es fundamental entender, en principio, a la categoría de monoides y explorar sus propiedades algebraicas y geométricas. Aunque este trabajo se basa fuertemente en el libro~\cite{Ogus:2018aa}, nosotros nos hemos dado a la tarea de exponer el material desde un punto de vista que esperamos sea más accesible para el lector que inicia en el estudio de la geometría algebraica y la aritmética, sin embargo, cualquier error tanto conceptual como tipográfico es responsabilidad de nosotros, los autores. 
\section{Bases de la teoría de monoides}

\subsection{Límites en la categoría de monoides}
 \begin{Def}\label{res:2.1}
Un monoide es una tripleta $(M,\bigstar,e_{M})$ que consiste de un conjunto M, una operación binaria asociativa $\bigstar$ y un elemento identidad $e_{M}$ de M.
 \end{Def}
 
\begin{Def}\label{res:2.2}
Un homomorfismo de monoides es una función $\theta: M \to N$ tal que $\theta(e_{M})=e_{N}$ y $\theta(m\bigstar n)= \theta(m)\bigstar \theta(n)$.
\end{Def}

Notemos que aunque el elemento identidad $e_{M}$ es único, la compatibilidad de $\theta$ con $e_{M}$ no es inmediata de la compatibilidad con $\bigstar$. Por ejemplo consideremos los monoides $(\mathbb{Z},\cdot,1)$, $(\mathbb{Z} \times \mathbb{Z},\cdot,(1,1))$ y la función $\theta:\mathbb{Z}\to \mathbb{Z} \times \mathbb{Z}:a \mapsto (a,0)$. Entonces, 
$$\theta(ab) = (ab,0) = (a,0)(b,0) =  \theta(a)\theta(b)$$
Pero notemos que $\theta(1)=(1,0) \neq (1,1)$.\\
Todos los monoides que vamos a considerar aquí son conmutativos a menos que se diga lo contrario, y escribiremos \textbf{Mon} para denotar la categoría de monoides conmutativos y homomorfismos de monoides. \\
Escribiremos solamente M o $(M,\bigstar)$ en lugar de $(M,\bigstar,e_{M})$ cuando no haya confusión. De manera similar, si a y b son dos elementos de un monoide $(M,\bigstar,e_{M})$, escribiremos ab o $a+b$ en lugar de $a\bigstar b$, y 1 o 0 en lugar de $e_{M}$. \\
El ejemplo más básico de un monoide es el conjunto $\bb{N}$ de números naturales con la suma. Sea M un monoide y $m\in M$, veamos que existe un único homomorfismo $\bb{N} \to M:1 \mapsto m$. Sea $\theta:\bb{N}\to M$ con $\theta(n)=\sum_{i=1}^{n}m$ y $\sum_{i=1}^{0}m := 0_{M}$. Claramente es un homomorfismo de monoides y tal que $\theta(1)=m$. Sea $\phi:\bb{N} \to M$ homomorfismo de monoides tal que $\phi(1)=m$. Notemos que,
$$\phi(n) = \phi(\sum_{i=1}^{n}1) = \sum_{i=1}^{n}\phi(1) = \sum_{i=1}^{n}m = \theta(n)$$
De modo que $\theta$ es único. Por tanto, $\bb{N}$ es un monoide libre con conjunto generador $\{1\}$. De manera más general, si S es un conjunto, el conjunto $\bb{N}^{(S)}$ de funciones $I:S \to \bb{N}$ tales que $I(s)=0$ para todos los $s\in S$ menos un número finito, equipado con la adición puntual de funciones como operación, es un monoide libre con base $S \subset \N^{(S)}$. Entonces
$$Hom_{Mon}(\bb{N}^{(S)},M) \cong Hom_{set}(S,M),$$
es decir, $S \mapsto \bb{N}^{(S)}$ es adjunto izquierdo del funtor olvidadizo de la categoría de monoides a la categoria de conjuntos. \\

Límites arbitrarios existen en la categoría de monoides, y su formación conmuta con el funtor olvidadizo a la categoría de conjuntos. En particular, la intersección de un conjunto de submonoides de M es un submonoide.

\begin{Def}\label{res:2.3}
Sea M un monoide y S un subconjunto de M. Definiremos el submonoide generado por S como la intersección de todos los monoides que contienen a S. 
\end{Def}

Dicho submonoide es el submonide más pequeño de M conteniendo a S. Si existe un subconjunto finito S de M que genere a M, diremos que M es finitamente generado como monoide. \\
Colimites arbitrarios de monoides también existen. Las sumas directas son construidas fácilmente: la suma directa $\oplus_{i\in I} M_{i}$ de una familia $\{M_{i}\}_{i\in I}$ de monoides es el submonoide del producto $\prod_{i}M_{i}$ que consiste de aquellos elementos $m=(m_{i})$ tales que $m_{i}=0$ para todo $i\in I$ excepto por un número finito. La existencia de límites y colímites se debe a resultados más generales de la Teoría de categorías los cuales se pueden encontrar en ~\cite{vakil2017rising}.

Vamos a investigar ahora los cocientes y las relaciones de equivalencia. Sea $\theta:P \to Q$ un homomorfismo de monoides. Notemos que el kernel de $\theta$ no nos ayuda mucho. Por ejemplo, el kernel del homomorfismo $\theta:\bb{N}\times \bb{N} \to \bb{N}:(a,b) \mapsto a+b$ es solo $\{0\}$ pero $\theta$ no es inyectiva ya que $\theta(1,2)=1+2=2+1=\theta(2,1)$. En su lugar, vamos a considerar el conjunto $E(\theta)=\{(p,q):\theta(p)=\theta(q)\}$, el cual es una relación de equivalencia en P, ya que la igualdad lo es. El hecho de que $\theta$ es un homomorfismo de monoides implica que $E(\theta)$ es un submonoide de $P\times P$.

\begin{Def}\label{res:2.4}
Sea P un monoide. Una relación de equivalencia en P que también es un submonoide de $P\times P$ es llamada una relación de congruencia en P. 
\end{Def}

Uno puede ver fácilmente que si E es una relación de congruencia en P, entonces el conjunto $P/E$ de clases de equivalencia tiene una única estructura de monoide que hace a la proyección $P \to P/E$ un homomorfismo de monoides. Entonces existe una biyección entre las relaciones de congruencia en P y los homomorfismo de monoides con dominio P. Todo esto se resume en la siguiente proposición. \\ 

\begin{lem}\label{res:2.5}\label{lem:ecualizador}
Sean $u_{1},u_{2}:P \to Q$ dos homomorfismos de monoides. Entonces el ecualizador de dichos morfismos es el conjunto $E=\{a\in P:u_{1}(a)=u_{2}(a)\}$, con el morfismo inclusión e.
\end{lem}

\begin{proof}
Claramente dado $a\in E$, tenemos que 
$$u_{1}\circ e(a)= u_{1}(a) = u_{2}(a) = u_{2}\circ e(a).$$
Sea D monoide con $v:D \to P$ tal que $u_{2}\circ v = u_{1}\circ v$. Sea $a\in D$, entonces $u_{2}\circ v(a) = u_{1}\circ v(a)$, por lo que $(v(a),v(a))\in E$. Por lo que definimos $g:D \to E:a \mapsto (v(a),v(a))$. Dicha función es claramente un homomorfismo de monoides único. Por tanto, E junto con el morfismo inclusión conforman el ecualizador de $u_{1}$y $u_{2}$.
\end{proof}

\begin{prop}\label{res:2.6}
1. Sea $\pi: P \to Q$ un homomorfismo de monoides sobre, y sea E el ecualizador de los dos mapeos 
$$\begin{tikzcd}
    P\times P \arrow{r}{p_{1},p_{2}} & P\arrow{r}{\pi} & Q
\end{tikzcd}$$
(a) E es una relación de congruencia en P. \\
(b) Q es el ecualizador de los dos mapeos
$$\begin{tikzcd}
    E \arrow{r}{e} & P\times P \arrow{r}{p_{1},p_{2}} & P
\end{tikzcd}$$
Entonces, el diagrama

$$\sqrd {E}{P}{Q}{P}{p_{1} \circ e}{\pi}{\pi}{p_{2} \circ e}$$

es cartesiano y cocartesiano. \\
2. Sea $E \subset P\times P$ una relación de congruencia en P, $Q:= P/E$ el conjunto de clases de equivalencia, y sea $\pi:P \to Q$ la función que toma un elemento en P y lo envía a su clase de equivalencia. \\
(a) Existe una única estructura de Q como monoide tal que $\pi:P \to Q$ es un homomorfismo de monoides. \\
(b) La inclusión $e:E \to P\times P$ es el ecualizador de los dos homomorfismos 
$$\begin{tikzcd}
    P\times P \arrow{r}{p_{1},p_{2}} & P\arrow{r}{\pi} & Q
\end{tikzcd}$$
y Q es el coecualizador de los dos homomorfismos 
$$\begin{tikzcd}
    E \arrow{r}{e} & P\times P \arrow{r}{p_{1},p_{2}} & P
\end{tikzcd}$$
Entonces, el diagrama de (1b) es cartesiano y cocartesiano. \\
\end{prop}

\begin{proof}
(1a) Por el Lema~\ref{lem:ecualizador}, el ecualizador es
\begin{equation*}
\begin{split}
E &= \{(a,b)\in P\times P:\pi\circ p_{1}(a,b)=\pi\circ p_{2}(a,b)\} \\
&= \{(a,b)\in P\times P: \pi(a)=\pi(b)\}
\end{split}
\end{equation*}
Veamos que E es una relación de congruencia en P. Para todo $p\in P$ se tiene que $\pi(p)=\pi(p)$, por lo que $(p,p)\in E$. Sea $(a,b)\in E$. Entonces $\pi(a)= \pi(b)$, es decir, $\pi(b)= \pi(a)$, así que $(b,a)\in E$. Si $(a,b),(b,c)\in E$, entonces $\pi(a)= \pi(b)$ y $\pi(b)= \pi(c)$, por lo que $\pi(a)= \pi(c)$. De este modo, $(a,c)\in E$. Con $0\in P$, tenemos que $(0,0)\in P$. Sean $(a,b),(c,d)\in E$. Entonces, $\pi(a)=\pi(b)$ y $\pi(c)=\pi(d)$. Como $\pi$ es un homomorfismo de monoides, tenemos que $\pi(a+c)=\pi(b+d)$, por lo que $(a+c,b+d)\in E$. Por tanto, E es una relación de congruencia en P. \\
(1b) Notemos que para todo $(a,b)\in E$,
$$\pi\circ p_{1}\circ e(a,b)=\pi(a) = \pi(b) = \pi\circ p_{2}\circ e(a,b)$$
Sea $f:P \to M$ tal que $f\circ p_{1}\circ e = f\circ p_{2}\circ e$. Entonces, para $(a,b)\in E$, tenemos que $f(a)=f(b)$. Sea $q\in Q$, como $\pi$ es sobre, existe $r\in P$ tal que $\pi(r)=q$. Definamos entonces $g:Q \to M:q \mapsto f(r)$. Notemos que está bien definida ya que si $\pi(r)=\pi(t)=q$, entonces $(r,t)\in E$, y por lo observación anterior $f(r)=f(t)$. Claramente es un homomorfismos de monoides único. Por tanto Q es el coecualizador de dichos mapeos. \\
(2a) La estructura de Q será la siguiente. Dado $a,b\in P$ definimos,
$$[a]+[b] = [a+b]$$ 
Y $[0]$ el elemento identidad en Q. Es fácil ver que esto define una estructura monoidal en Q. Además,
$$\pi(0) = [0]$$
Y,
$$\pi(a+b) = [a+b] = [a]+[b] = \pi(a)+\pi(b)$$
Por lo que $\pi$ es un homomorfismo de monoides. Es fácil ver que dicha estructura es única. \\ 
(2b) Notemos que para $(a,b)\in E$
$$\pi\circ p_{1} \circ e(a,b) = [a] = [b] = \pi\circ p_{2}\circ e(a,b).$$
Sea $f:M \to P\times P$ tal que $\pi\circ p_{1}\circ f = \pi\circ p_{2}\circ f$. Esto implica que para $r\in M$, si $f(r)=(a,b)$, entonces $[a]=[b]$, es decir, $(a,b)\in E$. Por lo que definimos $g:M \to E:r \mapsto (a,b)$. Claramente g es un homomorfismo de monoides único. Por tanto, E es el ecualizador de los dos morfismos. \\
Sea $(a,b)\in E$. Entonces, 
$$\pi\circ p_{1}\circ e(a,b)= [a] = [b] = \pi\circ p_{2}\circ e(a,b).$$
Sea $f:P \to M$ tal que $f\circ p_{1}\circ e = f\circ p_{2}\circ e$. Entonces, para $(a,b)\in E$, tenemos que $f(a)=f(b)$. Sea $[q]\in Q$. Así, definamos $g:Q \to M:[q] \mapsto f(q)$. Notemos que si $[q] = [p]$, entonces por lo anterior $f(p)=f(q)$. De modo que está bien definida y claramente es un homomorfismo único. Por tanto, Q es el coecualizador de los dos morfismos.
\end{proof}

Esta proposición induce una biyección entre las relaciones de congruencia en P y los homomorfismos sobre con dominio P de la siguiente forma,
\begin{equation*}
\begin{split}
\{Congruencias~en~P\} &\leftrightarrows \{Morfismos~sobre~con~dominio~P\} \\
E &\mapsto \pi:P \to P/E \\
(P\times P)\times_{Q}(P\times P) &\mapsfrom \theta: P \to Q
\end{split}
\end{equation*}

\begin{obs}\label{res:2.7}
Si $u_{1}:P \to Q$ es sobre y $u_{2}:Q^{\prime} \to Q$ es un homomorfismo de monoides, entonces el pullback $P\times_{Q} Q^{\prime}$ con morfismos $v_{1},v_{2}$ es nuevamente sobre. \\
Para ver esto sea $b\in Q^{\prime}$. Entonces, $u_{2}(b)\in Q$. Como $u_{1}$ es sobre existe $a\in P$ tal que $u_{1}(a)=u_{2}(b)$. Sabemos que el pullback en la categoría de monoides es  $P\times_{Q} Q^{\prime}= \{(x,y):u_{1}(x)=u_{2}(y)\}$. Entonces $(a,b)\in P\times_{Q} Q^{\prime}$, por lo que $v_{2}(a,b)=b$. 
Esto implica que $P \to Q$ y E son universalmente efectivos. Por otro lado, no todo epimorfismo en la categoria de monoides es sobre. De hecho, un homomorfismo de monoides es universalmente un epimorfismo si y solo si es sobre. \\
\end{obs}

\begin{prop}\label{res:2.8}
Sea M un monoide y $(E_{i})_{i\in I}$ una familia de relaciones de congruencia en M. Entonces la intersección $P=\cap_{i\in I}E_{i}$ es una relación de congruencia en M.
\end{prop}

\begin{proof}
Sabemos que la intersección de submonoides de M es un submonoide, por lo que veamos solamente que es una relación de equivalencia. \\
Tenemos que $(a,a)\in M_{i}$, para todo i y todo $a\in P$. Entonces $(a,a)\in P$. Sea $(a,b)\in P$, por lo que $(a,b)\in M_{i}$ para todo i. Como cada $M_{i}$ es una relación de congruencia $(b,a)\in M_{i}$, así $(b,a)\in P$. Finalmente supongamos que $(a,b),(b,c)\in P$. Entonces $(a,b),(b,c)\in M_{i}$ para todo i. Al ser relaciones de cogruencia tenemos que $(a,c)\in M_{i}$. Por tanto $(a,c)\in P$.
\end{proof}

\begin{Def}\label{res:2.9}
Sea M un monoide, E una relación de congruencia y S un subconjunto de E. Entonces diremos que $(S)$, la intersección de las relaciones de congruencia que contienen a S, es la relación de congruencia generada por S. 
\end{Def}

A continuación veremos una caracterización útil de las relaciones de congruencia generadas por un subconjunto de $P\times P$. \\ 

\begin{prop}\label{res:2.10}\label{prop:cong-gen}
Sea P un monoide conmutativo. \\
1. Una relación de equivalencia $E \subset P\times P$ en P es una relación de congruencia si y solo si $(a+p,b+p)\in E$ para todo $(a,b)\in E$ y $p\in P$. \\
2. Si S es un subconjunto de $P\times P$, sea 
$$S_{P}:= \{(a+p,b+p):(a,b)\in S,p\in P\}$$
Entonces la relación de congruencia E generada por S es la relación de equivalencia generada por $S_{P}$. Explicitamente, E es la unión de la diagonal con el conjunto de pares $(x,y)$ para los cuales existe una sucesión finita $(p_{0},...,p_{n})$ con $p_{0}=x$, $p_{n}=y$ y tal que para todo $i>0$, $(p_{i-1},p_{i})$ o $(p_{i},p_{i+1})$ pertenecen a $S_{P}$.
\end{prop}

\begin{proof}
1) $\Rightarrow$ Sea E una relación de equivalencia en P que es estable bajo la adición de elementos en P. Supón que (a,b) y (c,d) pertenecen a E. Entonces $(a+c,b+c)\in E$ y $(c+b,d+b)\in E$, y dado que P es conmutativo y E es transitiva se tiene que $(a,b)+(c,d)=(a+c,b+d)\in E$. Entonces es cerrada bajo suma. Dado que E es simétrica y $0\in P$, entonces $(0,0)\in E$. Por tanto, E es un submonoide de $P\times P$, es decir, una relación de congruencia en P.  \\
$\Leftarrow$ Sea E una relación de congruencia en P. Entonces, para todo $p\in P$, se tiene que $(p,p)\in E$ por ser E simétrica. Por ser un submonoide es cerrada bajo suma, entonces para todo $(a,b)\in E$, tendremos que $(a+p,b+p)=(a,b)+(p,p)\in E$. \\

2) Sea E la relación de congruencia generada por S y $E^{\prime}$ la relación de equivalencia generada por $S_{P}$. Dado que $S_{P} \subset E$ y E es una relación de equivalencia, se sigue que $E^{\prime} \subset E$. Es claro que $S_{P}$ es estable bajo la suma de elementos de la diagonal de $P\times P$. Entonces si $(p_{0},...,p_{n})$ es una sucesión tal que $(p_{i-1},p_{i})$ o $(p_{i},p_{i+1})\in S_{P}$ para todo $i>0$, entonces $(p_{0}+p,...,p_{n}+p)$ tiene la misma propiedad. Entonces si $(x,y)\in E^{\prime}$ y $p\in P$ tenemos que $(x+p,y+p)\in E^{\prime}$. Se sigue de (1) que $E^{\prime}$ es una relación de congruencia, por tanto $E^{\prime}=E$.
\end{proof}

\begin{obs}\label{res:2.11}
 Si P es un grupo abeliano y $E \subset P\times P$ es una relación de congruencia en P,entonces la imagen de E bajo el homomorfismo $h:P\oplus P \to P:(p_{1},p_{2}) \to p_{2}-p_{1}$ es un subgrupo K de P, y $E=h^{-1}(K)$. Por otro lado, la imagen inversa bajo h de cualquier subgrupo  de P es una relación de congruencia en P. Esto define una biyección entre los subgrupos de P y las relaciones de congruencia en P dada por,
\begin{equation*}
\begin{split}
\{subgrupos~de~P\} &\leftrightarrows \{Relaciones~de~congruencia~en~P\} \\
K &\mapsto h^{-1}(K) \\
h(E) &\mapsfrom E
\end{split}
\end{equation*}
\end{obs}

\begin{lem}\label{res:2.12}
Sean $\theta_{1},\theta_{2}: P \to Q$ homomorfismos de monoides. Entonces el coecualizador de $\theta_{1},\theta_{2}$ es el cociente de Q por la relación de congruencia en Q generada por el conjunto $\{(\theta_{1}(p),\theta_{2}(p)):p\in P\}$.
\end{lem}

\begin{proof}
Sea E la relación de congruencia generada por el conjunto $\{(\theta_{1}(p),\theta_{2}(p)):p\in P\}$ y $j:Q \to Q/E$ el morfismo proyección. Entonces para todo $p\in P$ tenemos que,
$$j\circ \theta_{1}(p)=[\theta_{1}(p)]=[\theta_{2}(p)]= j\circ \theta_{2}(p).$$
Sea $f:Q \to M$ tal que $f\circ \theta_{1}= f \circ \theta_{2}$. Sea $g:Q/E \to M: [a] \mapsto f(a)$. Claramente g es un homomorfismo que cumple que $g\circ j(a)=g([a])=f(a)$, y además es único. De modo que $(Q/E,j)$ cumplen la propiedad universal de los coecualizadores. 
\end{proof}

La existencia de colímites arbitrarios se sigue de la existencia de sumas directas y coecualizadores de pares de homomorfismos utilizando la siguiente construcción. Sea $\{P_{i},\theta_{a}\}$ un funtor de una categoría pequeña a la categoría de monoides, donde i corre sobre todos los objetos de I y a sobre las flechas $a:i(a) \to j(a)$ de I. Dado que la categoría de monoides tiene suma directa, podemos definir $Q:= \oplus \{P_{i}:i\in Ob(I)\}$ y $R:=\oplus \{P_{i(a)}:a\in Flechas(I)\}$, con los homomorfismos canónicos, 
$$\{u_{i}:P_{i}\to Q:i \in Ob(I)\} ~~~~~ \{u_{a}:P_{i(a)} \to R: a\in Flechas(I)\}.$$ 
Entonces existen homomorfismos únicos $\theta_{1},\theta_{2}:R \to Q$ tal que para todo a, $\theta_{1} \circ u_{a}= u_{i(a)}$ y $\theta_{2} \circ u_{a}= u_{j(a)}\circ \theta_{a}:$ \\

$$\sqrd{R}{P_{i(a)}}{P_{j(a)}}{Q}{u_{a}}{\theta_{a}}{u_{j(a)}}{\theta_{2}}$$

El colímite del funtor $\{P_{i},\theta_{a}\}$ es el coecualizador de $\theta_{1},\theta_{2}$. \\

\begin{Def}\label{res:2.13}
Una presentación de un monoide M es un diagrama coecualizado
$$L_{1} \rightrightarrows L_{0} \rightarrow M$$
donde $L_{0},L_{1}$ son libres.
\end{Def}

El monoide M se dice que es de presentación finita si admite una presentación con $L_{0},L_{1}$ libres y finitamente generados como monoides. \\

\begin{Def}\label{res:2.14}
La suma amalgamada de un par de homomorfismos de monoides $u_{i}:P \to Q_{i}$, denotada simplemente por $Q_{1}\oplus_{P}Q_{2}$ junto con dos homomorfismos $v_{i}:Q_{i} \to Q_{1}\oplus_{P}Q_{2}$, es el colímite del diagrama $Q_{1}\leftarrow P \rightarrow Q_{2}$. El par $(v_{1},v_{2})$ hacen que universalmente el diagrama 

$$\sqrd{P}{Q_{1}}{Q}{Q_{2}}{u_{1}}{v_{1}}{v_{1}}{v_{2}}$$

conmute.
\end{Def}

Esta suma amalgamada puede ser vista como el pushout de $u_{1}$ y $u_{2}$. También se puede ver como el coecualizador de los dos mapeos $(u_{1},0)$ y $(0,u_{2})$ que van de P a $Q_{1}\oplus Q_{2}$. \\
La siguiente proposición describe la suma amalgamada explícitamente. Su descripción se vuelve más simple si alguno de los monoides en cuestión es un grupo.  \\

\begin{prop}\label{res:2.15}\label{prop:suma-amalgamada}
 Sea $u_{i}:P \to Q_{i}$ un par de homomorfisos de monoides, sea Q su suma amalgamada, y sea E la relación de congruencia dada por el morfismo natural sobre $Q_{1}\oplus Q_{2}\to Q .$  \\
1. Sea S el conjunto de pares 
$$((q_{1},q_{2}),(q_{1}^{\prime},q_{2}^{\prime}))\in (Q_{1} \times Q_{2})\times (Q_{1}\times Q_{2})$$
tal que existe $p\in P$ con $q_{1}^{\prime}=u_{1}(p)+q_{1}$ y $q_{2}=u_{2}(p)+q_{2}^{\prime}$. Entonces E es el conjunto de pares  
$$(a,b)\in (Q_{1}\times Q_{2}) \times (Q_{1}\times Q_{2})$$
tal que existe una sucesión $(r_{0},...,r_{n})$ en $(Q_{1}\oplus Q_{2})$ tal que $(a,b)=(r_{0},r_{n})$ y $(r_{i},r_{i+1})$ está en S si i es par y en $S^{t}:=\{(a,b):(b,a) \in S\}$ si i es impar. \\
2. Sea $E^{\prime}$ el conjunto de pares $((q_{1},q_{2}),(q_{1}^{\prime},q_{2}^{\prime}))$ tales que existe p y $p^{\prime}$ en P con $q_{1}+u_{1}(p^{\prime})=q_{1}^{\prime}+u_{1}(p)$ y $q_{2}+u_{2}(p)=q_{2}^{\prime}+u_{2}(p^{\prime})$. Entonces $E^{\prime}$ es una relación de congruencia en $Q_{1}\oplus Q_{2}$ conteniendo a E, y si $P,Q_{1}~o~Q_{2}$ es un grupo, entonces $E=E^{\prime}$. \\
3. Si P es un grupo, entonces dos elementos de $Q_{1}\oplus Q_{2}$ son congruentes módulo E si y solo si ellos pertenecen a la misma órbita de la acción de P en $Q_{1}\oplus Q_{2}$ definida por $p(q_{1},q_{2})=(q_{1}+u_{1}(p),q_{2}+u_{2}(-p)).$ \\
4. Si P y $Q_{i}$ son grupos, entonces también los es $Q_{1}\oplus_{P} Q_{2}$, la cual es justamente la suma amalgamada en la categoría de grupos abelianos.
\end{prop}

\begin{proof}
Para probar (1) observemos primero que S es estable bajo la acción de la diagonal de  $(Q_{1} \times Q_{2})\times (Q_{1}\times Q_{2})$ y contiene la diagonal. Entonces, por la Proposición~\ref{prop:cong-gen}, la relación de congruencia R generada por S es el conjunto de pares (a,b) tal que existe $(r_{0},...,r_{n})$ con $r_{0}=a$, $r_{n}=b$ y tal que cada par $(r_{i},r_{i+1})$ esta en S o en $S^{t}$. Notemos por otro lado que si $(r_{i-1},r_{i})$ y $(r_{i},r_{i+1})$ están en S o en $S^{t}$ ambas, entonces también $(r_{i-1},r_{i+1})$, por lo que la sucesión puede ser acortada. Notemos que si $(r_{0},r_{1})$ pertenece a $S^{t}$, entonces $(r_{0},r_{0},r_{1})$ satisface lo descrito en (1). Esto demuestra que el conjunto descrito en (1) es una relación de congruencia. Dado que E contiene a S, y es la relación de congruencia más pequeña que lo contiene, se tiene que $E=R$. \\
Para probar (2) notemos primero que el conjunto $E^{\prime}$ es evidentemente simétrico y reflexivo. Para probar que es transitivo diremos que un par (a,b) une un par de elementos $(q_{1},q_{2})$ y $(q_{1}^{\prime},q_{2}^{\prime})$ de $Q_{1}\oplus Q_{2}$ si $q_{1}+u_{1}(b)=q_{1}^{\prime}+u_{1}(a)$ y $q_{2}+u_{2}(a)=q_{2}^{\prime}+u_{2}(b)$. Uno puede ver fácilmente que si (a,b) une a $(q_{1},q_{2})$ y $(q_{1}^{\prime},q_{2}^{\prime})$ y $(a^{\prime},b^{\prime})$ une a $(q_{1}^{\prime},q_{2}^{\prime})$ y $(q_{1}^{\prime \prime},q_{2}^{\prime \prime})$ entonces $(a+a^{\prime},b+b^{\prime})$ une a $(q_{1},q_{2})$ y $(q_{1}^{\prime \prime},q_{2}^{\prime \prime})$. Más aún, si (a,b) une a $(q_{1},q_{2})$ y $(q_{1}^{\prime},q_{2}^{\prime })$, entonces para cualquier $(\overline{q_{1}},\overline{q_{2}})$, (a,b) une a $(q_{1}+\overline{q_{1}},q_{2}+\overline{q_{2}})$ y $(q_{1}^{\prime}+\overline{q_{1}},q_{2}^{\prime}+\overline{q_{2}})$. Entonces por la Proposición~\ref{prop:cong-gen}, $E^{\prime}$ es una relación de congruencia en $Q_{1}\oplus Q_{2}$. Más aún, si $p\in P$, (p,0) une a $(u_{1}(p),0)$ y $(0,u_{2}(p))$, y dado que E es generada por tales pares, $E\subset E^{\prime}$. Si P o $Q_{i}$ es un grupo entonces $v:=v_{i}\circ u_{i}$
se factoriza a través del grupo $Q^{\ast}$ de elementos invertibles de Q. Si (a,b) une a $(q_{1},q_{2})$ y $(q_{1}^{\prime},q_{2}^{\prime})$, tenemos que 
\begin{equation*}
\begin{split}
v_{1}(q_{1})+v_{2}(q_{2})+v(a+b) &= v_{1}(q_{1}+u_{1}(b))+v_{2}(q_{2}+u_{2}(a)) \\
&= v_{1}(q_{1}^{\prime}+u_{1}(a))+v_{2}(q_{2}^{\prime}+u_{2}(b)) \\
&= v_{1}(q_{1}^{\prime})+v_{2}(q_{2}^{\prime})+v(a+b).
\end{split}
\end{equation*}
Dado que $v(a+b)\in Q^{\ast}$, se sigue que 
$$v_{1}(q_{1})+v_{2}(q_{2})=v_{1}(q_{1}^{\prime})+v_{2}(q_{2}^{\prime}).$$
Entonces $E^{\prime}\subset E$. Esto prueba (2), y (3) y (4) son consecuencias inmediatas.
\end{proof}

\begin{exa}\label{res:2.16}
Tomando $Q_{2}=0$ en la Proposición~\ref{prop:suma-amalgamada} obtenemos el cokernel del morfismo $u_{1}:P \to Q_{1}$, o de manera equivalente, el coecualizador de $u_{1}$ y el mapeo cero de P a $Q_{1}$. Si P es un submonoide de Q, uno escribe $Q \to Q/P$ para denotar el cokernel de la inclusión $P \to Q$, y se sigue de (2) en la proposición anterior que dos elementos $q,q^{\prime}\in Q$ tienen la misma imagen en $Q/P$ si y solo si existen $p,p^{\prime}\in P$ tales que $q+p=q^{\prime}+p^{\prime}$. Por ejemplo, el cokernel del morfismo diagonal $\bb{N} \to \bb{N}\times \bb{N}$ es el homomorfismo 
$$\bb{N}\times \bb{N} \to \bb{Z}:(a,b) \mapsto a-b.$$
\end{exa}

Nota $Q/P$ puede ser cero aun si P es un submonoide propio de Q. Esto pasa, por ejemplo, si P es el submonoide de $Q:=\mathbb{N} \times \mathbb{N}$ generado por (1,0) y (1,1). Si $P^{\prime}$ es un submonoide de Q conteniendo a P, entonces $P^{\prime}/P$ es un submonoide de $Q/P$ y el mapeo natural $(Q/P)/(P^{\prime}/P) \to Q/P^{\prime}$ es un isomorfismo. \\

\subsection{Acciones de Monoides} 
\begin{obs}\label{res:2.17}
Si S es un conjunto, el conjunto de funciones de S en si misma forman un monoide (no necesariamente conmutativo) End(S) bajo la composición de funciones.
\end{obs}

\begin{Def}\label{res:2.18}
Si Q es un monoide, una acción de Q en un conjunto S es un homomorfismo de monoides $\theta_{S}:Q \to End(S)$.
\end{Def}

En este contexto escribimos la ley del monoide Q de manera multiplicativa, y si $s\in S$,$q\in Q$, escribimos $qs$ para denotar $\theta_{S}(q)(s)$.

\begin{Def}\label{res:2.19}
Un Q-conjunto es un conjunto con una acción de Q, y $Ens_{Q}$ denotará a la categoría de Q-conjuntos, con la evidente noción de homomorfismos.
\end{Def}

Si S es un Q-conjunto y $s\in S$, la imagen del mapeo $Q\to S:q \mapsto qs$ es el Q-subconjunto minimal estable de S conteniendo a s., llamado la trayectoria de s en S. \\

\begin{Def}\label{res:2.20}
Una base para un Q-conjunto S es un conjunto T y una función $i:T \to S$ tal que el mapeo inducido $\phi:Q\times T \to S:(q,t)\mapsto qi(t)$ es una biyección. Si tal base existe diremos que S es un Q-conjunto libre.
\end{Def}

Un Q-conjunto libre con base $i:T\to S$ satisface la propiedad universal de los objetos libres. Sea $S^{\prime}$ un Q-conjunto con $f:T \to S^{\prime}$. Sea $a\in S$. Como $\phi$ es biyectiva existe $(q,t\in Q\times T)$ único tal que $qi(t)=a$. Definamos $g:S \to S^{\prime}:a \to f(t)$. Por lo que claramente el diagrama conmuta. De modo que se cumple la propiedad universal de los objetos libres. \\

Si T es cualquier conjunto y si $Q\times T$ es equipado con la acción definida por $q^{\prime}(q,t)=(q^{\prime}q,t)$, entonces el mapeo $T \to Q\times T$ enviando t a (1,t) es una base. De este modo, el funtor que toma un conjunto T y lo envía al Q-conjunto libre $Q\times T$ es adjunto izquierdo del funtor olvidadizo de la categoría de monoides a la de conjuntos. Nota que si G es un grupo y S es un G-conjunto, entonces S tiene una G-base si y solo si la acción es libre en el sentido de que $gs=s$ implica que $g=1$. \\

La categoría $Ens_{Q}$ de Q-conjuntos admite límites proyectivos arbitrarios, y su formación conmuta con el funtor olvidadizo a la categoría de conjuntos. Esta es una consecuencia formal del hecho que el funtor $Ens_{Q}\to Ens $ tiene una adjunta izquierda. En particular, si S y T son Q-conjuntos, entonces Q actúa en $S\times T$ por la acción $q(s,t)=(qs,qt)$. Esto define el producto en $Ens_{Q}.$ 
Colímites en $Ens_{Q}$ también existen. La suma directa de una familia de Q-conjuntos $\{S_{i}\}$ es solo la unión disjunta con la Q-acción evidente. Para entender la construcción de los cocientes en $Ens_{Q}$, notemos que si $\pi: S \to T$ es un mapeo sobre de Q-conjuntos, la correspondiente relación de equivalencia $E\subset S \times S$ es un Q-subconjunto. 

\begin{Def}\label{res:2.21}
Sea S un Q-conjunto y E una relación de equivalencia en S. Si E es un Q-subconjunto de $S\times S$ diremos que E es una relación de congruencia en S. 
\end{Def}

Por el contrario, si E es cualquier relación de congruencia en S, entonces existe una única estructura de Q-conjunto de S/E tal que la proyección $S \to S/E$ es un morfismos de Q-conjuntos. Cuando $S=Q$ actuando regularmente en si mismo, la noción de relación de congruencia en Q como monoide coincide con la la noción de relación de congruencia como Q-conjunto, gracias a la Proposición~\ref{prop:cong-gen}. Más aún, el análogo a (2) de la Proposición~\ref{prop:cong-gen} se cumple para Q-conjuntos, y en particular la relación de equivalencia generada por un subconjunto de $S\times S$ que es estable bajo la adición de elementos en Q es una relación de congruencia. Si $u,v:S^{\prime}\to S$ son dos morfismos, entonces el coecualizador de u y v es el cociente de S con la relación de congruencia generada por el conjunto $\{(u(s),v(s))\}$. La existencia de colímites se sigue. \\

Si S y T son Q-conjuntos, el $Hom_{Q}(S,T)$  tiene una acción natural de Q, dada por $(qh)(s):=qh(s)=h(qs)$, con $h\in Hom_{Q}(S,T)$,$q\in Q$, $s\in S$. También existe una construcción para el producto tensorial de $Q-conjuntos$.

\begin{Def}\label{res:2.22}
Si $S,T,W$ son Q-conjuntos entonces un Q-bimorfismo es una función $\beta:S\times T \to W$ tal que $\beta(qs,t)=\beta(s,qt)=q\beta(s,t)$.
\end{Def}

\begin{prop}\label{res:2.23}
Sean S y T dos monoides. Entonces el producto tensorial de S y T existe. 
\end{prop}

\begin{proof}
 Para construirlo consideramos al producto $S\times T$ con la acción $q(s,t)=(qs,t)$, y consideramos la relación de equivalencia R en $S\times T$ generada por el conjunto de pares 
$$((qs,t),(s,qt)),q\in Q,s\in S,t\in T.$$
Notemos que este conjunto de pares es estable bajo la acción de Q, ya que si $q^{\prime}\in Q$ entonces
$$((q^{\prime}qs,t),(q^{\prime}s,qt)) = ((q(q^{\prime}s),t),(q^{\prime}s,qt))$$
Se sigue que la relación de equivalencia R es una relación de congruencia. Entonces la proyección $\pi:S\times T \to (S\times T)/R$ es un Q-bimorfismo y satisface la propiedad universal del producto tensorial.
\end{proof}

Si Q es un grupo conmutativo, entonces $S\otimes_{Q}T$ puede ser construido en la forma usual como el espacio de órbitas de la acción de Q en $S\times T$ dada por $q(s,t)=(qs,q^{-1}t)$.
En general uno tiene el isomorfismo natural de Q-conjuntos,
\begin{equation*}
\begin{split}
Hom_{Q}(S\otimes_{Q}T,W) &\cong Hom_{Q}(S,Hom_{Q}(T,W)) \\
\beta:S\otimes_{Q}T \to W &\mapsto \gamma:S \to Hom(T,W):s \mapsto \beta(s,_)
\end{split}
\end{equation*}
Para el reverso si tenemos $\gamma:S \to Hom(T,W)$ definimos el Q-bimorfismo $f:S \times T \to W:(s,t)\to \gamma(s)(t)$. Por la propiedad universal del producto tensorial existe $g:S\otimes_{Q}T\to W$ único tal que el diagrama del producto tensorial conmuta. Por lo que,

\begin{equation*}
\begin{split}
g:S\otimes_{Q}T \to W &\mapsfrom \gamma:S \to Hom(T,W) 
\end{split}
\end{equation*}

Se sigue que para un T fijo, el funtor $S \mapsto S\otimes_{Q}T$ conmuta con colímites. \\
Sea $\theta:P \to Q$ un homomorfismo de monoides. Entonces $\theta$ define una acción de P en Q dada por $pq=\theta(p)q$. Si S es un P-conjunto, el producto tensorial $Q\otimes_{P}S$ tiene la acción natural de Q-conjunto dada por $q(q^{\prime}\otimes s)= (qq^{\prime}\otimes s)$, y el mapeo $S \to Q\otimes_{Q}S$ enviando s a $1\otimes s$ es un morfismo de P-conjuntos. Si $\theta_{i}:P \to Q_{i}$ es un par de homomorfismos, entonces existe una estructura única de monoide para $Q_{1}\otimes_{P}Q_{2}$ tal que
$$ (q_{1}\otimes q_{2})(q_{1}^{\prime}\otimes q_{2}^{\prime})= (q_{1}q_{1}^{\prime}\otimes q_{2}q_{2}^{\prime})$$
y esta es también la única estructura monoidal para la cual los mapeos naturales $Q_{i} \to Q_{1}\otimes_{P}q_{2}$ son homomorfismos. Se puede ver que esta estructura convierte a $Q_{1}\otimes_{P}Q_{2}$ es una suma amalgamada de $Q_{1}$ y $Q_{2}$ a través de P. \\

\begin{obs}\label{res:2.24}
Hemos visto que para un Q-conjunto fijo T,el funtor $S \mapsto S\otimes_{Q}T$ conmuta con colímites. Quizás no es sorpresa que no conmuta con límites en general. Queremos enfatizar que este funtor no conmuta con productos finitos áún cuando T es libre. De hecho, si T tiene base $\omega$, entonces $S\otimes_{Q}T \cong S\times \omega$, y si la cardinalidad de $\Omega$ es más grande que uno, el funtor $S \mapsto S\times \Omega$ no conmuta con productos. Este hecho complica el cálculo de productos tensoriales por medio de sus generadores y relaciones. De hecho, supón que $F \to S$ es un homomorfismo sobre y $E$ es la correspondiente relación de equivalencia en F, con F libre. Entonces $F \to S$ es el coecualizador de los dos mapeos $E\rightrightarrows F$, y da dado que $\otimes{Q}T$ conmuta con colímites, se sigue que $F\otimes_{Q}T \to S\otimes_{Q}T$ es el coecualizador de $E\otimes_{Q}T \rightrightarrows F\otimes_{Q}T$. Por otro lado, el mapeo natural 

$$(F\times F)\otimes_{Q}T \to (F\otimes_{Q}T)\times (F\otimes_{Q}T)$$

no necesariamente es un isomorfismo, y la imagen de $E\otimes_{Q}T$ en $(F\otimes_{Q}T)\times (F\otimes_{Q}T)$ podría no ser una relación de equivalencia. Por lo que uno se queda con el problema de calcular la relación de congruencia que genera.  
\end{obs}

\begin{Def}\label{res:2.25}
Sea Q un monoide y S un Q-conjunto. El transportador de S es la categoría $\fk{T}_{Q}S$ cuyos objetos son los elementos de S y los morfismos que van de un objeto s a uno t son los elementos $q\in Q$
tal que $qs=t$. El transportador de un monoide Q es el transportador de Q considerado como Q-conjunto, y es denotado simplemente por $\fk{T}S$. 
\end{Def}

Asociado con la categoría $\fk{T}_{Q}S$ existe un conjunto parcialmente ordenado que vale la pena hacer explícito. \\

\begin{Def}\label{res:2.26}
Sea Q un monoide y S un Q-conjunto. Si $s,t\in S$, diremos que $s\leq t$ si existe $q\in Q$ tal que $qs=t$, y que $s\sim t$ si $s\leq t$ y $t\leq s$. 
\end{Def}

Supongamos que $s\leq t$ y $t\leq w$. Entonces existe $q,p\in Q$ tales que $sq=t$ y $tp=w$. Tenemos que $qp\in Q$, y así
$$s(qp)=tp=w$$
De modo que $s\leq w$. Sea $s\in S$. Entonces,
$$1 s = s $$
Por lo que $s\leq s$. Por tanto esta relación define un preorden en S. Veamos que $\sim$ es una relación de congruencia en S.  Claramente es reflexiva y simétrica. Supongamos que $s\sim t$ y $t\sim w$. Así que $s\leq t$, $t\leq s$, $t\leq w$ y $w\leq t$. Dado que $\leq$ es un preorden en S, tenemos que $s\leq w$ y $w\leq s$, de modo que $s\sim w$. Supongamos que $q\in Q$ y $a\sim b$. Entonces, $a\leq b$ y $b\leq a$. Por lo que existen $t,r\in Q$ tales que $ta=b$ y $rb=a$. Así que utilizando la conmutatividad $tqa=qb$ y $rqb=qa$. De modo que $qa\leq qb$ y $qb \leq qa$. Así, $qa \sim qb$, y en conclusión $\sim $ es una relación de congruencia en S. Además es fácil ver que $\leq$ es un orden parcial en $S/ \sim$. \\
Utilizaremos esta notación especialmente cuando $S=Q$ con la representación regular. Dado que $\sim$ es una relación de congruencia, se sigue de la Proposición~\ref{prop:cong-gen} que $Q/\sim$ hereda una estructura de monoide.

\subsection{Monoides \textit{integral}es, finos y saturados}

\begin{prop}\label{res:2.27}
Sea M un monoide conmutativo. Entonces existe un grupo $M^{gp}$ y un homomorfismo universal $\lambda_{M}:M \to M^{gp}$ tal que dado cualquier otro homomorfismos $f:M \to D$ existe un único homomorfismo $g:M^{gp} \to D$ con  $g\circ \lambda_{M} = f$.  
\end{prop}

\begin{proof}
Uno puede construir $M^{gp}$ como el conjunto de clases de equivalencia de pares $(x,y)\in M\times M$, donde $(x,y)$ es equivalente a $(x^{\prime},y^{\prime})$ si y solo si existe $z\in M$ tal que $x+y^{\prime}+z=x^{\prime}+y+z$. Uno escribe $x-y$ para denotar la clase de equivalencia de $(x,y)$, y entonces $(x-y)+(x^{\prime}-y^{\prime})= (x+x^{\prime})-(y+y^{\prime}).$ Y considerar 
$$\lambda_{M}:M \to M^{gp}:m \mapsto [(m,0)].$$
\end{proof}

Entonces
$$Hom(M^{gp},G) \cong Hom(M,G),$$
es decir, el funtor $M \mapsto M^{gp}$ es adjunto izquierdo del funtor inclusión que va de la categoría de grupos a la categoría de monoides; dado que tiene adjunta derecha, automáticamente conmuta con la formación de límites directos.   

\begin{Def}\label{res:2.28}
Si M es un monoide, definimos $$M^{\ast}:=\set{m\in M:\exists n\in M,m+n=0}.$$
\end{Def}

Notemos que $M^{\ast}$ es un submonoide de M. Es de hecho un grupo. Este grupo es llamado el grupo de unidades de M. Cualquier homomorfismo de un grupo a M se factoriza de manera única a través de $M^{\ast}$. Así que 
$$Hom(M^{\ast},G)=Hom(M,G),$$
es decir, $M\mapsto M^{\ast}$ es adjunto derecho del funtor inclusión que va de la categoría de grupos a la categoría de monoides. El grupo $M^{\ast}$ actúa naturalmente en M por traslación. En particular definimos, $\overline{M}= M/M^{\ast}$

\begin{Def}\label{res:2.29}
Un monoide conmutativo se dice que es: \\
1. \textit{Sharp} si $M^{\ast}=\set{0}$ \\
2. \textit{Dull} si $M^{\ast}=M$, es decir, si M es un grupo. \\
3. u-\textit{integral} si $m\in M$,$u\in M^{\ast}$ y $m+u=m$ implica que $u=0$.\\
4. Cuasi-\textit{integral} si $m,n\in M$ y $m+n=m$ implica que $n=0$. \\
5. \textit{Integral} si $m,n,p\in M$ y $m+n=p+n$ implica que $m=p$.
\end{Def}

Evidentemente todo monoide \textit{integral} es cuasi \textit{integral}, y todo cuasi \textit{integral} es u-\textit{integral}. Veamos que si M es es u-\textit{integral}, entonces $M^{\ast}$ actúa libremente en M. Supongamos que $p+m=m$, con $p\in M^{\ast}$ y $m\in M$. Entonces por ser u-\textit{integral}, se sigue que $p=0$. El mapeo universal $\lambda_{M}$ es inyectivo si y solo si M es \textit{integral}. Supongamos que es inyectivo y $m+n=p+n$ con $m,n,p\in M$. Entonces, $[(m,0)]=[(p,0)]$, y al ser $\lambda_{M}$ inyectiva tenemos que $m=p$. Supongamos que M es \textit{integral}. Si $\lambda_{M}(m)=\lambda_{M}(p)$, tenemos que $[(m,0)]=[(p,0)]$, es decir, existe $n\in M$ tal que $m+n=p+n$. Luego, al ser M \textit{integral}, se sigue que $m=p$, por lo que $\lambda_{M}$ es inyectiva. De manera completamente análoga podemos demostrar que el mapeo inducido $M^{\ast} \to M^{gp}$ es inyectivo si y solo si M es u-\textit{integral}. \\
Notemos que para cualquier monoide M el cociente $\overline{M}$ es \textit{\textit{Sharp}}, ya que si $[m]\in \overline{M}$ es una unidad, existe $[n]$ tal que $[m+n]=[m]+[n]=[0]$. Lo cual implica que existe $u\in M^{\ast}$ tal que $m+n=u$. Así, $m+(n-u)=0$, por lo que $m\in M^{\ast}$, es decir, $[m]=[0]$. De este modo el mapeo $M \to \overline{M}$ es el homomorfismo universal que va de M a un monoide \textit{\textit{Sharp}}. \\
Para cualquier monoide M, el monoide $M/ \sim$ es \textit{Sharp}. Sea $[m]\in M/ \sim$ unidad, entonces existe $[n]\in M/ \sim$ tal que $[m+n][m]+[n]=[0]$. Por lo que existe $u\in M$ tal que $u+m+n=0$, así $(u+n)+m=0$. De modo que $[m]=[0]$.  \\
Veamos que si M es quasi-\textit{integral}, entonces 
\begin{equation*}
\begin{split}
\phi:M/M^{\ast} &\to M/ \sim \\
\overline{m} &\mapsto [m]
\end{split}
\end{equation*}
 
es un isomorfismo. Supongamos que $\phi(\overline{m})=\phi(\overline{n})$. De modo que $[m]=[n]$, y existe entonces $f,g\in M$ tales que $f+m=n$ y $g+n=m$. Luego, $f+g+m=g+n=m$, y al ser M quasi-\textit{integral}, tenemos que $f+g=0$, por lo cual $f,g\in M^{\ast}$. De este modo $\overline{m}=\overline{n}$, y así $\phi$ es inyectiva. Claramente es sobre. Por lo tanto es un isomorfismo. 

\begin{obs}\label{res:2.30}
La formación del grupo $M^{gp}$ conmuta con límites directos pero no con productos fibrados en general. Por ejemplo, sea $s:\bb{N}^{2} \to \bb{N}:(a,b)\to a+b$ y $t:\bb{N}^{2} \to \bb{N}:(a,b) \to 0$. Entonces el ecualizador de s y t es cero. Sin embargo, el mapeo asociado en grupos $\bb{Z}^{2} \to \bb{Z}$ es la antidiagonal $\bb{Z} \to \bb{Z}^{2}$, enviando c a $(c,-c)$. \\

Por otro lado, es cierto que si $\theta:P \to Q$ es inyectiva y Q \textit{integral} entonces $\theta^{gp}:P^{gp} \to Q^{gp}$ es también inyectiva. Supongamos que $\theta^{gp}([a,b])=\theta([c,d])$. Entonces $[\theta(a),\theta(b)]=[\theta(c),\theta(d)]$. Así que existe $z\in Q$ tal que $\theta(a)+\theta(d)+z=\theta(c)+\theta(b)+z$. Como $\theta$ es un homomorfismo $\theta(a+d)+z=\theta(c+b)+z$, y al ser Q \textit{integral} se sigue que $\theta(a+d)=\theta(c+b)$. Como $\theta$ es inyectiva, se sigue que $a+d=c+b$ por lo cual $[a,b]=[c,d]$, y así $\theta^{gp}$ es inyectiva también. 
\end{obs}

\begin{prop}\label{res:2.31}\label{prop:Q-sobre-P}
Si Q es un monoide \textit{integral} y P es un submonoide, el mapeo natural $Q/P \to Q^{gp}/P^{gp}$ es inyectivo. Entonces $Q/P$ puede ser identificado con la imagen de Q en $Q^{gp}/P^{gp}$. Un monoide Q es \textit{integral} si y solo si es u-\textit{integral} y $\overline{Q}$ es \textit{integral}.  
\end{prop}

\begin{proof}
Si $q,q^{\prime}$ son dos elementos de Q con la misma imagen en $Q^{gp}/P^{gp}$, entonces existen $p,p^{\prime}$ tales que $q-q^{\prime}=p-p^{\prime}$ en $Q^{gp}$. Dado que Q es \textit{integral}, $q+p^{\prime}=p+q^{\prime}$ en Q. Se sigue que q y $q^{\prime}$ tienen la misma imagen en $Q/P$. Entonces $Q/P \to Q^{gp}/P^{gp}$ es inyectivo, y $Q/P$ es \textit{integral}. En particular si Q es \textit{integral} también lo es $\overline{Q}$. Por el contrario, supongamos que Q es u-\textit{integral} y $\overline{Q}$ es \textit{integral}, y $q,q^{\prime},p$ son elementos de Q tales que $p+q=p+q^{\prime}$. Dado que $\overline{Q}$ es \textit{integral}, existe una unidad u tal que $q^{\prime}=q+u$.Entonces, $q+p=q+p+u$, Dado que Q es u-\textit{integral}, $u=0$ y entonces $q=q^{\prime}$. Esto demuestra que Q es \textit{integral}. 
\end{proof}

\begin{Def}\label{res:2.32}
Sea $Mon^{int}$ la subcategoría completa de Mon cuyos objetos son los monoides \textit{integral}es. 
\end{Def}

Para cualquier monoide M sea $M^{int}$ la imagen de $\lambda_{M}$. Entonces
$$Hom(M^{int},N)=Hom(M,N),$$
es decir, el funtor $M \mapsto M^{int}$ es adjunto izquierdo del funtor inclusión $Mon^{int} \to Mon.$

\begin{prop}\label{res:2.33}
Sea Q la suma amalgamada de dos homomorfismos $u_{i}:P \to Q_{i}$ en la categoría de monoides. Entonces $Q^{int}$ es la suma amalgamada de $u_{i}^{int}:P^{int} \to Q_{i}^{int}$ en la categoría $Mon^{int}$, y puede ser identificada naturalmente con la imagen de Q en $Q_{1}^{gp } \oplus_{P^{gp}}Q_{2}^{gp}$ Si $P,Q_{1},Q_{2}$ son \textit{integral}es y cualquiera de estos monoides es un grupo, entonces Q es \textit{integral}. 
\end{prop}

\begin{proof}
El hecho de que $Q^{int}$ es la suma amalgamada de $u_{i}^{int}$ en $Mon^{int}$ es una consecuencia formal del hecho de que $M \mapsto M^{int}$ preserva colímites. Más aún, dado que $M \mapsto M^{gp}$ también preserva colímites, $Q^{gp} \cong Q_{1}^{gp} \oplus_{P^{gp}} Q_{2}^{gp}$. Se sigue que $Q^{int}$ es la imagen de Q en $Q^{gp} \cong Q_{1}^{gp} \oplus_{P^{gp}} Q_{2}^{gp}$. Supongamos ahora que cualquiera de $P,Q_{1},Q_{2}$ es un grupo y que $q,q^{\prime}$ son dos elementos de Q con la misma imagen en $Q^{gp}$. Elige $(q_{1},q_{2})$ y $(q_{1}^{\prime},q_{2}^{\prime})$ en $Q_{1}\oplus Q_{2}$ mapeando a $q$ y $q^{\prime}$ respectivamente. Entonces $v_{1}(q_{1})+v_{2}(q_{2})=v_{1}(q_{1}^{\prime})+v_{2}(q_{2}^{\prime})$ en $Q^{gp}$, y entonces existen elementos a y b en P tales que $(q_{1}^{\prime}-q_{1},q_{2}^{\prime}-q_{2})= (u_{1}(a-b),u_{2}(b-a))$. Entonces $q_{1}^{\prime}+u_{1}(b)=q_{1}+u_{1}(a)$ y $q_{2}^{\prime}+u_{2}(a)=q_{2}+u_{2}(b)$. Se sigue de (2) de la Proposición~\ref{prop:suma-amalgamada} que $v_{1}(q_{1})+v_{2}(q_{2})= v_{1}(q_{1}^{\prime})+v_{2}(q_{2}^{\prime})$ en Q, es decir, que $q=q^{\prime}$. Por tanto, el mapeo $Q\to Q_{1}^{gp}\oplus_{P^{gp}}Q_{2}^{gp}$ es inyectivo y Q es \textit{integral}. 
\end{proof}

\begin{Def}\label{res:2.34}
Un monoide que es \textit{integral} y finitamente generado es llamado un monoide fino.
\end{Def}

\begin{Def}\label{res:2.35}
Un monoide Q es llamado saturado si es \textit{integral} y además dado $mq\in Q^{gp}$, con $m\in \bb{Z}^{+}$, entonces $q\in Q$.
\end{Def}

Por ejemplo el conjunto de los números naturales $\bb{N}$ es saturado pero el conjunto de los enteros más grandes que un número natural $d>1$ no es saturado. 

\begin{prop}\label{res:2.36}
Sea Q un monoide \textit{integral}. \\
1. El homomorfismo natural $Q^{gp}/Q^{\ast} \to \overline{Q}^{gp}$ es un isomorfismo. \\
2. Si Q es saturado, entonces $\overline{Q}^{gp}$ es libre de torsión. \\
3. $Q^{sat}$ es el conjunto de elementos x en $Q^{gp}$ tales que existe $n\in \bb{Z}^{+}$ con $nx\in Q$ es un submonoide saturado de $Q^{gp}$, y el funtor $Q \mapsto Q^{sat}$ es adjunto izquierdo del funtor inclusión que va de la categoría $Mon^{sat}$ de monoides saturados a $Mon^{int}.$ \\
4. Q es saturado si y solo si $\overline{Q}$ es saturado. Un elemento de Q es una unidad si y solo si su imagen en $Q^{sat}$ es una unidad.\\
5. El mapeo natural $Q^{sat}/Q^{\ast} \to \overline{Q}^{sat}$ es un isomorfismo. Más aún, toda unidad de $\overline{Q}^{sat}$ es torsión, y el mapeo natural
$$\overline{Q^{sat}} \to \overline{\overline{Q}^{sat}}$$
es un isomorfismo. 
\end{prop} 

\begin{proof}
Supón que $q_{1},q_{2}\in Q$ y $q_{2}-q_{1}$ se mapea al cero en $\overline{Q}^{gp}$. Dado que $\overline{Q} \subset \overline{Q}^{gp}$, $\overline{q_{1}}=\overline{q_{2}}$ en $\overline{Q}$, entonces existe $u\in Q^{\ast}$ con $q_{2}=q_{1}+u$. Entonces $q_{2}-q_{1}=u\in Q^{\ast}$. Esto prueba (1). Supongamos que Q es saturado y $q\in Q^{gp}$ se mapea a un elemento de torsión $\overline{q}$ de $\overline{Q}^{gp}$. Entonces $nq\in Q^{\ast}$ para algún $n\in \bb{Z}^{+}$, y dado que Q es saturado, tenemos que $q\in Q$.  El hecho que nq pertenece a $Q^{\ast}$ implica ahora que q pertence a $Q^{\ast}$, así $\overline{q}=0\in \overline{Q}$. Entonces $\overline{Q}^{gp}$ es libre de torsión. Si p y q son elementos de $Q^{gp}$ con $mq,np\in Q$, entonces $mn(p+q)\in Q$, por lo que se sigue que $Q^{sat}$ es un submonoide de $Q^{gp}$. Entonces $(Q^{sat})^{gp}=Q^{gp}$ y si $q\in Q^{gp}$ y $nq\in Q^{sat}$, entonces existe un $m\in \bb{Z}^{+}$ tal que $mnq\in Q$. Se sigue que $q\in Q^{sat}$, así que $Q^{sat}$ es saturado. La verificación de la adjunción del funtor $Q \mapsto Q^{sat}$ es inmediata como la de (4). \\
Si $q\in Q^{sat}$ y $\overline{q}$ es una unidad de $\overline{Q}^{sat}$, entonces existe un elemento p de $Q^{sat}$ tal que $q+p\in Q^{\ast}$. Entonces existen $m,n\in \bb{Z}^{+}$ tal que mq y np pertenecen a Q. Pero entonces $mnp+mnq\in Q^{\ast}$, y así mnq es una unidad de Q. Esto muestra que $\overline{q}$ es un elemento de torsión de $\overline{Q}^{sat}$. Es claro que el mapeo de (5) es sobre. \\
Supongamos que p y q son dos elementos de $Q^{sat}$ con la misma imagen en $\overline{\overline{Q}^{sat}}$ Entonces $q-p\in Q^{gp}$ se mapea a una unidad de $\overline{Q}^{sat}$, y entonces a un elemento de torsión de $\overline{Q}^{sat} \subset \overline{Q}^{gp}$. De modo que $mp-mq\in (Q^{sat})^{\ast}$  para algún m. Entonces también $mp-mq\in Q^{\ast}$, así que $q-p$ es una unidad de $Q^{sat}$ y p y q tienen la misma imagen en $\overline{Q}^{sat}$, lo cual prueba la inyectividad.    
\end{proof}

\begin{Def}\label{res:2.37}
A los monoides que son finos y saturados los llamaremos fs-monoides.
\end{Def}

Este tipo de monoides son de vital importancia en la geometría algebraica logarítmica. 

\begin{Def}\label{res:2.38}
Un monoide P se dice tórico si es fino, saturado y además $P^{gp}$ es libre de torsión. 
\end{Def}

En este caso $P^{gp}$ puede ser visto como el grupo de caracteres de un toro algebraico. Los esquemas que surgen de los monoides tóricos forman los componentes básicos de la geometría tórica.

\begin{prop}\label{res:2.39}
Sea $\set{M_{i}:i\in I}$ un sistema dirigido de monoides donde cada uno satisface una de las siguientes propiedades P: \textit{integral}, saturado, \textit{Dull}. Entonces el límite directo M satisface P. 
\end{prop}

\begin{proof}
Supón que cada $M_{i}$ es \textit{integral} y sea $m\in M$.Entonces existe $i\in I$ y $m_{i}\in M_{i}$ tal que $m_{i}$ se mapea a m. Para cada $i \to j$, sea $m_{ji}$ la imagen de $m_{i}$ en $M_{j}$. Entonces la multiplicación $m_{ji}: M_{j} \to M_{j}$ es inyectiva. Se sigue que el límite de estos mapeos, es decir, la multiplicación por m, es también inyectiva. Entonces M es \textit{integral}. \\
Supongamos que cada $M_{i}$ es saturado y $x\in M^{gp}$ con $nx\in M$ para algún $n>0$. Dado que la formación de $M^{gp}$ conmuta con límites directos, existe $i\in I$ y $x_{i}\in M_{i}^{gp}$ mapeando a x. Reemplazando i por algún elemento al cual mapea, podemos asumir que existe $m_{i}\in M_{i}$ mapeando a nx. De nuevo reemplazando i, podemos asumir que $nx_{i}=m_{i}$ en $M_{i}$. Dado que $M_{i}$ es saturado se sigue entonces que $x_{i}\in M_{i}$ y entonces $x\in M$. Por tanto M es saturado. Dado que la formación de límites directos es la misma en la categoría de monoides conmutativos y grupos, el límite directo de monoides \textit{Dull} es \textit{Dull}. 
\end{proof}

\begin{Def}\label{res:2.40}
Un monoide M se dice Valuativo si es \textit{integral} y para todo $x\in M^{gp}$ ocurre que x o $-x$ pertenecen a M. 
\end{Def}

Esto es equivalente a decir que el preorden en $M^{gp}$ definida por la acción de M es un preorden total. El monoide $\bb{N}$ es valuativo y si V es un anillo de valuación, el submonoide $V^{\prime}$ de elementos no cero de V es valuativo. Todo monoide valuativo es saturado. 
Si R es un anillo conmutativo, su monoide multiplicativo $(R,\cdot,1)$ no es cuasi \textit{integral} a menos de que sea el monoide cero ya que $0\cdot 0 = 0 \cdot 1$. Por otro lado, el conjunto $R^{\prime}
$ de divisores no cero de R forma un submonoide \textit{integral} del monoide multiplicativo R. Por ejemplo, $\bb{Z}^{\prime}$ es \textit{integral}, y $\overline{\bb{Z}}^{\prime} = \bb{Z}^{\prime}/(\pm)$ es un monoide libre, generado por los números primos. Si R es un anillo de valuación discreta $\overline{R}^{\prime}=R^{\prime}/ R^{\ast}$ es generado libremente por la imagen de un uniforme de $R^{\prime}$. Aunque existe un único isomorfismo de monoides $R^{\prime}/ R^{\ast} \cong \bb{N}$, este isomorfismo no es funtorial: si $R \to S$ es una extensión finita de anillos de valuación con índice de ramificación e, el mapeo inducido $\overline{R}^{\prime} \to \overline{S}^{\prime}$ envía el generador único de $\overline{R}^{\prime}$ a e veces el de $\overline{S}^{\prime}$. \\
Si Q es un monoide conmutativo \textit{Sharp}, los Q-conjuntos libres son muy rígidos, como la siguiente observación muestra. 

\begin{prop}\label{res:2.41}
Sea Q un monoide conmutativo \textit{Sharp}, y sea S un Q-conjunto libre. Entonces cualquier base para S es única bajo isomorfismos. Explícitamente, toda base $i:T \to S$ induce una biyección entre T y $S\setminus Q^{+}S.$
\end{prop}

\begin{proof}
Sea $i:T \to S$ una base de S. Dado que el mapeo inducido $Q\times T \to S$ es biyectivo, i debe ser inyectiva. Verifiquemos que si $t\in T$, entonces $i(t)\in S\setminus Q^{+}S.$ Supongamos que $i(t)=qs$ con $q\in Q$ y $s\in S$. Entonces existe un único $(q^{\prime},t^{\prime})\in Q\times T$ tal que $s=q^{\prime}i(t^{\prime})$, y también $i(t)=qq^{\prime}i(t^{\prime})$. Entonces $(1,t)$ y $(qq^{\prime},t^{\prime})$ son dos elementos de $Q\times T$ con la misma imagen en S, por  lo que $qq^{\prime}=1$. Dado que Q es conmutativo, $q\in Q^{\ast}$ y al ser Q \textit{Sharp}, $q=1$. Como q era arbitrario, $i(t)\notin Q^{+}S$. Por otro lado, supón que $s\in S\setminus Q^{+}S$. Dado que i es una base para S, existe algún $(q,t)\in Q\times T$ con $qi(t)=s$; como $s\notin Q^{+}S$, $q=1$ y $s=i(t)$. Entonces el mapeo inducido que va de T a $S \setminus Q^{+}S$ es también sobre.  
\end{proof}

\subsection{Ideales, caras y localización}
\begin{Def}\label{res:2.42}
Un ideal de un monoide M es un subconjunto I tal que para todo $q\in I$ y $k\in M$ tenemos que $q+k\in I$. 
\end{Def}

\begin{Def}\label{res:2.43}
Un ideal I de M es primo si $I\neq M$ y si $q+k\in I$ entonces $q\in I$ o $k\in I$.  
\end{Def}

\begin{Def}\label{res:2.44}
Una cara F de M es un submonoide tal que $p+q\in F$ implica que $p\in F$ y $q\in F$. 
\end{Def}

Por ejemplo el conjunto vacío es un ideal primo de M, como también el conjunto $M^{+}$ de no unidades de M. Notemos que cualquier ideal I conteniendo una unidad u debe ser todo M, ya que de este modo $0=u-u \in I$, y así para todo $m\in M$ se tendría que $m=0+m\in I$. Por tanto todo ideal propio debe estar contenido en $M^{+}$. Entonces $M^{+}$ es el único ideal propio maximal de M. Más aún, como el vacío es subconjunto de todo conjunto. tenemos que $\emptyset$ es el único ideal minimal de M. En muchos aspectos un monoide es análogo a un anillo local. 

\begin{Def}\label{res:2.45}
Un homomorfismo de monoides $\theta:P \to Q$ se dice local si $\theta^{-1}(Q^{+})=P^{+}$, o de manera equivalente, $\theta^{-1}(Q^{\ast})=P^{\ast}$. 
\end{Def}

Observemos que dado un monoide M una cara F es justamente un submonoide cuyo complemento es un ideal primo, y un ideal primo P es un ideal cuyo complemento es una cara. Entonces, 
\begin{equation*}
\begin{split}
\fk{p} &\mapsto M\setminus \fk{p} \\
M\setminus F &\mapsfrom F
\end{split}
\end{equation*}

define una biyección entre las caras y los ideales primos de M. \\
El conjunto de unidades $M^{\ast}$ es la cara más pequeña de M ya que su complemento $M^{+}$ es el ideal maximal. Además M es la cara más grande ya que su complemento es $\emptyset$ el cual es el ideal minimal.  La noción de una cara de un monoide corresponde a la noción de un subconjunto multiplicativo saturado de un anillo. \\
Como en el caso de anillos, la intersección de una familia de ideales es un ideal, pero para monoides la unión de una familia de ideales también es un ideal.  Más aún la unión de una familia de ideales primos es un ideal primo y la intersección de una familia de caras es una cara. 

\begin{Def}\label{res:2.46}
Sea M un monoide y T un subconjunto de M. Entonces la cara generada por T, denota por $<T>$, es la intersección de todas las caras que contienen a T.  
\end{Def}

Esto es análogo al conjunto saturado multiplicativo generado por un subconjunto de un anillo. 

\begin{Def}\label{res:2.47}
Sea M un monoide. Definimos el ideal interior $I_{M}$ como la intersección de todos los ideales primos no vacíos de M. 
\end{Def}

\begin{Def}\label{res:2.48}
Sea M un monoide. Denotamos por $\Spec{M}$ al conjunto de ideales primos de M.
\end{Def}

\begin{Def}\label{res:2.49}
Sea M un monoide. Para cada ideal I denotamos por Z(I) al conjunto de ideales primos de M que contienen a I. 
\end{Def}

Entonces si M es un monoide y $(I_{\lambda})$ es alguna familia de ideales, $\cup I_{\lambda}$ es un ideal de M y $Z(\cup I_{\lambda}) = \cap Z(I_{\lambda})$.   También si I y J son ideales, también lo es IJ el conjunto de todos los elementos de la forma i+j con $i\in I$ y $j\in J$, y $Z(I)\cap Z(J)=Z(IJ)= Z(I\cup J)$. Por lo tanto el conjunto de todos los subconjuntos de $\Spec{M}$ de la forma Z(I) con I un ideal es cerrado bajo intersecciones y uniones finitas, y entonces define una topología en $\Spec{M}$ llamada la Topología de Zariski. Dado que M tienen un único ideal primo minimal, $\Spec{M}$ tiene un único punto genérico, es decir,
\begin{equation*}
\begin{split}
\overline{\set{\emptyset}} &= \cap_{\set{\emptyset} \subset Z(I)} Z(I) \\
&= \cap_{I \subset \emptyset} Z(I) \\
&= Z(\emptyset) \\
&= \Spec{M}
\end{split}
\end{equation*}

y en particular es irreducible. Dado que M tiene un único ideal maximal este tiene un único punto cerrado, es decir, 

\begin{equation*}
\begin{split}
\overline{\set{M^{+}}} &= \cap_{\set{M^{+}} \subset Z(I)} Z(I) \\
&= \cap_{I \subset M^{+}} Z(I) \\
&= \cap_{I\in \Spec{M}} Z(I) \\
&= \set{M^{+}}.
\end{split}
\end{equation*}

\begin{Def}\label{res:2.50}
Sea $f\in M$ y F la cara generada por dicho elemento. Definimos,
$$D(f)=S_{f}:=\set{\fk{p}:f\notin \fk{p}} = \set{\fk{p}:\fk{p} \cap F=\emptyset}$$
\end{Def}

Es fácil ver que $S_{f}$ es un abierto en $\Spec{M}$, y que el conjunto de este tipo de subconjuntos forman una base para la topología en $\Spec{M}$. \\
Notemos que $\Spec{M}$ nunca es vacío ya que siempre contiene al ideal primo vacío y al ideal maximal $M^{+}$. Estos dos puntos coinciden si y solo si todo elemento de M es una unidad, es decir, si y solo si M es un grupo.  \\ 
Si $\theta:P\to Q$ es un homomorfismos de monoides, tenemos que la imagen inversa de un ideal es un ideal, la imagen inversa de un ideal primo es un ideal primo y la imagen inversa de una cara es una cara. Entonces $\theta$ induce un mapeo continuo,
$$\theta^{\ast}:\Spec{Q} \to \Spec{P}: \fk{p} \mapsto \theta^{-1}(\fk{p})$$
La relación de preorden es usado para describir los ideales y las caras de un monoide, como muestra la siguiente proposición. 

\begin{prop}\label{res:2.51}\label{prop:caras-ideales}
Sea S un subconjunto de un monoide Q y P el submonoide de Q generado por S. \\
1. El ideal (S) de Q generado por S es el conjunto de todos los $q\in Q$ tal que $s\leq q$, para algún $s\in S$. \\ 
2. La cara $<S>$ de Q generada por S es el conjunto $P^{\prime}$ de elementos $q\in Q$ para los cuales existe $p\in P$ tal que $q\leq p$. En particular, la cara generada por un elemento $p\in Q$ es el conjunto de todos los elementos $q\in Q$ tal que $q\leq np$ para algún $n\in \bb{N}$. \\ 
3. Si Q es \textit{integral}, entonces $Q/P$ es \textit{Sharp} si y solo si $P^{gp}\cap Q$ es una cara de Q. En particular si F es una cara de Q, entonces $Q/F$ es \textit{Sharp}. 
\end{prop}

\begin{proof}
1) Sea $T=\set{q\in Q:s\leq q,s\in S}$. Consideremos $q\in T$ y $k\in Q$. Entonces existe $s\in S$ tal que $s\leq q$. De modo que $s\leq q+k$ y así $q+k\in T$. Por tanto T es un ideal de Q. Además claramente $S\subset T$ ya que $s\leq s$. Por otro lado, sea $q\in T$, por lo que existe $s\in S$ con $s\leq q$. Más aún, existe $f\in Q$ tal que $f+s=q$. Sea $S\subset I$ ideal de Q. Entonces $q=f+s\in I$ ya que $s\in S\subset I$. Al ser I arbitrario tenemos que $q\in (S)$. Concluimos por doble contención que $T=(S)$. \\ 

2) Notemos que un submonoide F de Q es una cara si y solo si F contiene a q siempre que $q\leq f$ para algún  $f\in F$. Entonces $<S>$ contiene $P^{\prime}$. Dado que de hecho $P^{\prime}$ es necesariamente una cara de Q conteniendo S, se sigue que $P^{\prime}=<S>.$ \\ 

3) Supongamos que Q es \textit{integral}. Entonces Q/P puede ser identificado con la imagen de Q en $Q^{gp}/P^{gp}$, por la Proposición~\ref{prop:Q-sobre-P}. De modo que un elemento $q\in Q$ mapea a 0 en $Q/P$ si y solo si $q\in Q\cap P^{gp}$, y q mapea a una unidad en Q/P si y solo si existe un elemento $q^{\prime}\in Q$ tal que $q+q^{\prime}\in P^{gp}$, es decir, si y solo si $q\in <Q\cap P^{gp}>$. Esto muestra que Q/P es \textit{Sharp} si y solo si $Q\cap P^{gp}$ es una cara de Q. Finalmente notemos que si F es una cara de Q, y $q\in Q\cap F^{gp}$, entonces $q+f\in F$, para algún $f\in F$; entonces $q\in F$.  
\end{proof}

\begin{cor}\label{res:2.52}
Sea K un ideal de un monoide Q y 
$$\sqrt{K}=\set{q\in Q:nq\in K,n\in \bb{N}}$$
su radical. \\ 
1. $\sqrt{K}$ es un ideal radical, es decir, $\sqrt{\sqrt{K}}=\sqrt{K}$ \\ 
2. $\sqrt{K}$ es la intersección de todos los ideales primos que contienen a K. \\ 
3. El mapeo $I \to Z(I)$ induce una biyección entre los ideales radicales de Q y los subconjuntos cerrados de $\Spec{Q}$, con inversa $S \mapsto \cap \set{p:p\in S}.$ \\
4. Un subconjunto cerrado S de Q es irreducible si y solo si el correspondiente ideal radical es primo. 
\end{cor}

\begin{proof}
1) Procedamos por doble contención. Sea $q\in \sqrt{\sqrt{K}}$, así que existe $n\in \bb{N}$ tal que $nq \in \sqrt{K}$. Nuevamente, existe $m\in \bb{N}$ tal que $m(nq)\in K$. De modo que $q\in \sqrt{K}$. Así $\sqrt{\sqrt{K}} \subset \sqrt{K}$. La otra inclusión es clara. \\ 
2) Es claro que $\sqrt{K}$ está contenido en todo ideal primo conteniendo a K. Por el contrario, supongamos que $q\in Q\setminus  \sqrt{K}$ y sea f un elemento de la cara F de Q generada por q. Por la Proposición~\ref{prop:caras-ideales}, existe n y $q^{\prime}$ tal que $nq=f+q^{\prime}$. Dado que $nq\notin K$, lo mismo es cierto para f. Esto muestra que $F\cap K = \emptyset$, y entonces $p= Q \setminus F$ es un ideal primo de Q conteniendo a K pero no a q. \\
3) Es claro que $Z(J) \subset Z(I)$ si $I\subset J$. Si S es cualquier subconjunto de $\Spec{Q}$, entonces $\cap S:= \cap \set{p:p\in S}$ es claramente un ideal radical de Q, y $S\subset Z(\cap S)$, dado que para todo primo $p\in S$, $\cap S \subset p$. Más aún, si I es un ideal cualquiera de Q y $S\subset Z(I)$ entonces $I\subset p$ para todo $p\in S$, y entonces $I \subset \cap S$ y $Z(\cap S) \subset Z(I)$. Entonces $Z(\cap S)$ es la clausura de S. En particular si S es cerrado, $S=Z(\cap S)$. Por otro lado, si K es un ideal radical, (2) establece que $K=\cap Z(K)$. Esto completa la demostración de (3). \\
4)  Si S es cerrado y Z(S) es un ideal primo p, entonces $p\in S$ y S es la clausura de $\set{p}$ y entonces es irreducible. Por el contrario, si S es irreducible $a+b \in \cap S$, entonces $a+b$ pertenece a todo $p\in S$; entonces todo p contiene a o b, y también $S \subset Z(a) \cap Z(b)$. Ya que  S es irreducible, este contiene a Z(a) o a Z(b). Si por ejemplo $S \subset Z(a)$, se sigue que $\sqrt{a} \subset \cap S$,y dado que $\cap S$ es un ideal radical se tiene que $a \in \cap S$. Por tanto $\cap S$ es primo.     
\end{proof}

\begin{prop}\label{res:2.53}
Sea M un monoide, S un subconjunto de M y E un M-conjunto. Entonces existe un M-conjunto, denotado por $S^{-1}E$, en el cual los elementos de S actúan biyectivamente y un mapeo $\lambda_{S}:E\to S^{-1}E$ con la siguiente propiedad universal: para cualquier homomorfismo de M-conjuntos $f:E \to E^{\prime}$ tal que los elementos de S actúan biyectivamente en él, existe un único M-mapeo $g:S^{-1}E \to E^{\prime}$ tal que 
$$\xymatrix{
E \ar[dr]_f \ar[r]^\lambda & S^{-1}E \ar[d]^g \\
& E^{\prime}
}
$$

el diagrama conmuta. El morfismo $\lambda_{S}$ es llamado la localización de E en S. Un morfismo de M conjuntos $\phi:E\to E^{\prime}$ induce un morfismo $\phi_{S}:S^{-1}E \to S^{-1}E^{\prime}$; si $\phi$ es inyectivo (suprayectivo) entonces $\phi_{S}$ es inyectivo (suprayectivo).  
\end{prop}

\begin{proof}
Sea T el submonoide de M generado por S. El conjunto $S^{-1}E$ puede ser construido como el conjunto de clases de equivalencia de pares $(e,t)\in E\times T$, donde $(e,t)\equiv (e^{\prime},t^{\prime})$ si, y solo si $ett^{\prime \prime}=e^{\prime}t^{\prime}t^{\prime \prime}$ para algún $t^{\prime \prime}\in T$. Y definimos $\lambda_{S}$ como
$$\lambda_{S}: E \to S^{-1}E:e \to [(e,1)]$$
Y la acción de M en $S^{-1}E$ se define como $m[e,t]=[me,t]$. Veamos ahora que dichas definiciones satisfacen la propiedad universal. \\
Sea $\theta: M \to End(S^{-1}E)$ el homomorfismo de la acción definida anteriormente y $s\in S$. Probemos que $\theta(s):S^{-1}E \to S^{-1}E$ es biyectiva. Supongamos que $\theta(s)([a,b])=\theta(s)([c,d])$, entonces $[sa,b]=[sc,d]$. Por lo que existe $t\in T$ tal que $sadt=scbt$. Se sigue que $[a,b]=[c,d]$ por o que $\theta(s)$ es inyectiva. \\
Sea $[e,t]\in S^{-1}E$. Notemos que $st\in T$ ya que $s\in S$. Por lo que claramente $\theta(s)([e,st])=[se,st]=[e,t]$. En conclusión $\theta(s)$ es biyectiva. Sea $f:E\to E^{\prime}$ con S actuando biyectivamente en $E^{\prime}$. Sea $[e,t]\in S^{-1}E$ definimos $g:E \to E^{\prime}:[e,t]\mapsto f(te)$. Entonces,
$$g \circ \lambda_{S}(e)=g([e,1])=f(1\cdot e)=f(e)$$
Por lo tanto nuestras definiciones satisfacen la propiedad universal. \\ 
Supongamos que $\phi:E \to E^{\prime}$ es inyectiva y que $\phi_{S}([e,t])=\phi_{S}([a,b])$. Entonces, $[\phi(e),t]=[\phi(a),b]$, y así que existe $w\in T$ tal que $\phi(e)bw=\phi(a)tw$. Como $\phi$ es un M-morfismo se sigue que $\phi(bwe)=\phi(twa)$. Al ser $\phi$ inyectiva tendremos que $bwe=twa$, por lo cual $[e,t]=[a,b]$. De modo que $\phi_{S}$ es inyectiva.   
\end{proof}

Sea M un monoide y F una cara de M. Si $\fk{p} = M \setminus F$ es el ideal primo correspondiente, escribiremos $E_{p}$ en lugar de $S^{-1}E$. 

\begin{Def}\label{res:2.54}
Un M-conjunto E es llamado M-regular si los elementos de M actúan inyectivamente en E. 
\end{Def}

Si este es el caso, el morfismo localización $\lambda_{S}:E \to S^{-1}E$ es inyectivo, para todo subconjunto S de M. \\

El caso más importante de la definición anterior es la representación regular, cuando M actúa sobre $E=M$ por traslación. Entonces $M_{S}:=S^{-1}M$ tiene una estructura única de monoide dada por
$$[a,b][c,d]=[ac,bd]$$
y para la cual $\lambda_{S}$ es compatible con las M-acciones. El morfismo $\lambda_{S}:M \to M_{S}$ está también caracterizado por una propiedad universal: cualquier homomorfismo $\lambda:M \to N$ con la propiedad $\lambda(s)\in N^{\ast}$ para cada $s\in S$ se factoriza de manera única a través de $M_{S}$. De hecho, todo elemento de la cara $<S>$ generada por S mapea a una unidad en $S^{-1}M$, y $\lambda_{S}^{-1}(M_{S}^{\ast}) = <S>$. De hecho, si $m\in <S>$, se sigue por la parte (2) de la Proposición~\ref{prop:caras-ideales} existe $m^{\prime}\in M$ tal que $mm^{\prime}$ pertenece al submonoide T de M generado por S. Entonces $\lambda_{S}(mm^{\prime})$ es una unidad en $M_{S}$, por lo que también lo es $m^{\prime}$. Por el contrario, si $\lambda_{S}(m)$ es una unidad de $M_{S}$, entonces existe $m^{\prime}\in M$ y $t^{\prime}\in T$ tal que $(m,1)(m^{\prime},t^{\prime})\equiv (1,1)$. Esto significa que existe $t\in T$ tal que $mm^{\prime}t=t^{\prime}t$. Dado que $tt^{\prime}\in <S>$, se sigue que $m\in <S>$. Si M es \textit{integral}, entonces el mapeo natural $S^{-1}M \to M^{gp}$ es inyectivo, y $S^{-1}M$ puede ser identificado con el conjunto de elementos en $M^{gp}$ de la forma $m-t$ con $m\in M$ y t en la cara de M generada por S. Si $\theta: M \to N$ es un morfismo de monoides y S es un subconjunto de M, escribimos $S^{-1}N$ para denotar la localización de N por la imagen de S bajo $\theta$, cuando no hay riesgo de confusión. Notemos que si E es un M-conjunto, entonces el monoide localizado $M_{S}$ actúa naturalmente en $E_{S}$, y de hecho el mapeo natural $M_{S}\otimes_{M}E \to E_{S}$ es un isomorfismo.

\begin{Def}\label{res:2.55}
Sea Q un monoide. \\
1. La dimensión de Q es la longitud máxima d de todas las cadenas de ideales primos 
$$\emptyset =\fk{p}_{0} \subset \fk{p}_{1} \subset \cdots \subset \fk{p}_{d} = Q^{+},$$
es decir la dimensión de krull del espacio topológico $\Spec{Q}.$ \\ 
2. Si $\fk{p}\in \Spec{Q}$, la altura de $\fk{p}$, denotada por $ht(\fk{p})$, es la longitud máxima de una cadena de ideales primos 
$$\fk{p}_{h} \subset \cdots \subset \fk{p}_{1} \subset \fk{p}_{0}=\fk{p}.$$
\end{Def}

Si $\fk{p}$ es un ideal primo de Q, el mapeo $\Spec{Q_{\fk{p}}} \to \Spec{Q}$ inducido por el mapeo localización $\lambda:Q \to Q_{\fk{p}}$ es inyectivo e identifica $\Spec{Q_{\fk{p}}}$  con el subconjunto de $\Spec{Q}$ que consiste de aquellos ideales primos que contienen a p. De manera equivalente, $F \mapsto \lambda^{-1}(F)$ es una biyección entre el conjunto de caras de $Q_{\fk{p}}$ y el conjunto de caras de Q conteniendo a $Q\setminus \fk{p}$. Esta biyección preserva la topología y el orden. En particular, todo ideal de $Q_{\fk{p}}$ es inducido por un ideal de Q. Más aún, tenemos que $ht(\fk{p})=dim(Q_{\fk{p}})$. Si Q es fino, $\Spec{Q}$ es un espacio topológico finito. La siguiente proposición será demostrada hasta después:

\begin{prop}\label{res:2.56}
Sea Q un monoide \textit{integral}. \\
1. $\Spec{Q}$ es un conjunto finito si Q es finitamente generado. \\
2. $dim(Q)\leq rank(\overline{Q}^{gp})$, con igualdad si Q es fino. \\ 
3. Si Q es fino, todo cadena maximal de ideales primos tiene longitud $dim(Q)$.     
\end{prop}

\subsection{Monoides idealizados}
\begin{Def}\label{res:2.57}
Un monoide idealizado es un par $(M,K)$, donde M es un monoide y K es un ideal de M. 
\end{Def}

\begin{Def}\label{res:2.58}
Un homomorfismo de monoides idealizados 
$$\theta:(P,I) \to (Q,J)$$
es un homomorfismo de monoides $P \to Q$ que envía I a J. 
\end{Def}

\begin{Def}\label{res:2.59}
Un ideal primo de un monoide idealizado (M,K) es un ideal primo P de M tal que $K \subset P$.
\end{Def}

\begin{Def}\label{res:2.60}
Una cara de un monoide idealizado (M,K) es una cara F de M tal que $F \cap K= \emptyset$.    
\end{Def}

\begin{Def}\label{res:2.61}
$\Spec{(M,K)}$ es el conjunto de ideales primos de $(M,K)$.
\end{Def}

Notemos que $\Spec{(M,K)}$ es vacío si y solo si $K=M$. Si $\Spec{(M,K)}=\emptyset$, entonces no hay ideales primos de M que contengan a K, por lo que $K=M$. Por otro lado, si $K=M$ y dado que todo ideal primo es propio, tenemos que no existe $\fk{p}\in \Spec{M}$ tal que $M\subset \fk{p}$. Por lo que definiremos lo siguiente.

\begin{Def}\label{res:2.62}
Diremos que un monoide idealizado (M,K) es aceptable si K es un ideal propio de M o M es el monoide cero.  
\end{Def}

Escribiremos Moni para denotar la categoría de monoides idealizados aceptables. El funtor $Moni \to Mon: (Q,I) \to Q$ tiene un adjunta izquierda completamente fiel, la cual toma un monoide P y lo envía al monoide $(P,\emptyset)$. Entonces, podemos ver a Mon como una subcategoría completa de Moni.

\begin{obs}\label{res:2.63}
La dimensión de Krull de (M,K) es la longitud máxima de las cadenas de ideales primos de (M,K), o de manera equivalente, de una cadena de caras de (M,K). Si C es una cadena maximal de caras, entonces $F=\cup C$ es otra cara de (M,K) y entonces pertenece a C; más aún, cada miembro de C es una cara de F. Entonces el conjunto de caras de (M,K) admite elementos maximales, y la dimensión de (M,K) es igual al máximo de la dimensión de sus caras. 
\end{obs}

Límites y colímites existen en la categoría de monoides idealizados, y son compatibles con el funtor olvidadizo a la categoría de monoides. 

\begin{obs}\label{res:2.64}
El pushout de un par de homomorfismos $u_{i}:(P,I) \to (Q_{i},J_{i})$ en la categoría de monoides idealizados está dado por los mapeos obvios $v_{i}:(Q_{i},J_{i}) \to (Q,J)$ donde $Q_{i} \to Q$ es el pushout como monoides y J es el ideal de Q generado por las imágenes de $J_{i}$. 
\end{obs}
\section{Finitud, convexidad y dualidad}

\subsection{Finitud}
\begin{prop}\label{res:3.1}
Un monoide M es finitamente generado si y solo si $M^{\ast}$ es finitamente generado como grupo y $\overline{M}$ es finitamente generado como monoide. 
\end{prop}

\begin{proof}
Si S es un conjunto finito de generadores para M, entonces cualquier elemento no cero m de M se puede escribir como una suma $\sum n_{i} s_{i}$ con $n_{i}>0$. Si m es una unidad, entonces cada $s_{i}$ lo es, y se sigue que $M^{\ast}$ es generada por el conjunto finito $S \cap M^{\ast}$. Ya que $M\to \overline{M}$ es sobre, $\overline{M}$ es finitamente generado como monoide. Por el contrario, supongamos que $\set{s_{i}}$ es un conjunto finito de generadores del grupo $M^{\ast}$ y $\set{t_{j}}$ es un subconjunto finito de M cuya imagen en $\overline{M}$ genera a $\overline{M}$ como monoide. Entonces el conjunto, $\set{s_{i},-s_{i},t_{i}}$ genera a M como monoide.    
\end{proof}

\begin{Def}\label{res:3.2}
Sea M un monoide, X un M-conjunto y $S\subset X$. Decimos que $s\in S$ es un elemento M-minimal de S si, para cualquier $s^{\prime} \in S$ con $s^{\prime} \leq s$ se tiene que $s\leq s^{\prime}$.
\end{Def}
Si no hay peligro de confusión solo diremos minimal en lugar de M-minimal. Por ejemplo con la acción regular de M en si mismo, las unidades son los elementos minimales de M. Si M es cuasi \textit{integral} los elementos minimales del ideal maximal $M^{+}$ son llamados irreducibles. Entonces, un elemento $c\in M$ es irreducible si y solo si $c\in M^{+}$ y si $c=a+b$, entonces a o b es una unidad de M. 

\begin{prop}\label{res:3.3}
Sea M un monoide cuasi \textit{integral} \textit{sharp}. Entonces todo conjunto de generadores de M contiene a todos los elementos irreducibles de M. Si en adición M es finitamente generado, entonces el conjunto de elementos irreducibles de M es finito y genera a M. 
\end{prop}

\begin{proof}
 El primer enunciado es claro. Supongamos ahora que M es finitamente generado. Es claro que todo conjunto finito de generadores contiene un conjunto minimal de generadores. Sea S tal conjunto minimal. Veamos que todo elemento $x\in S$ es irreducible. Si $x=y+z$ con $y,z\in M$, podemos escribir a y y z como suma de elementos de S, digamos $y=\sum a_{s}s$ y $z=\sum b_{s}s$, donde $a_{s},b_{s}\in \bb{N}$ para todo $s\in S$. Entonces $x=\sum c_{s}s$ con $c_{s}=a_{s}+b_{s}$. Sea $S^{\prime} = S\setminus \set{x}$, así que $x=c_{x}x+m^{\prime}$ donde $m^{\prime}$ está en el submonoide $M^{\prime}$ de M generado por $S^{\prime}$. Si $c_{x}=0$, entonces $x\in M^{\prime}$ y $S^{\prime}$ genera a M, contradiciendo la minimalidad de S. Entonces $c_{x}\geq 1$ y podemos escribir $(c_{x}-1)x+m^{\prime}=0$, y al ser \textit{sharp} implica que $c_{x}=1$ y $m^{\prime}=0$. Entonces $y=a_{x}x$ y $z=b_{x}x$, donde $a_{x}+b_{x}=1$. Entonces exactamente uno de los dos de y y z son cero, por lo que x es irreducible. Ya que S genera a M por hipótesis y al ser los elementos de S irreducibles, M es generado por su conjunto de elementos irreducibles. Ya que S contiene a todos los elementos irreducibles de M y es finito, estos solo pueden ser un número finito.  
\end{proof}

\begin{cor}\label{res:3.4}
El grupo de automorfismos de un monoide \textit{sharp} fino es finito y está contenido en el grupo de permutaciones del conjunto de sus elementos irreducibles. 
\end{cor}

\begin{proof}
Sea M un monoide \textit{sharp} y fino. Por la proposición~\ref{res:3.3} tenemos que el conjunto de elementos irreducibles S de M es finito y genera a M. Dado que los automorfismos de M están determinados por sus generadores y estos son finitos tendremos que el número de automorfismos también será finito. Sea X el conjunto de sus elementos irreducibles y $\phi$ un automorfismo de M. Veamos que $s\in M$ es unidad si y solo si $\phi(s)$ es unidad. Supongamos s es unidad, entonces existe $t\in M$ tal que $s+t=0$. De modo que $\phi(s)+\phi(t)=\phi(s+t)=0$. Por otro lado si $\phi(s)$ es unidad tenemos que $\phi(s)+m=0$. Al ser $\phi$ un automorfismo existe $t\in M$ tal que $\phi(s+t)=\phi(s)+\phi(t)=0$ y por tanto $s+t=0$. Sea $s\in X$ entonces por lo anterior $\phi(s)$ no es una unidad. Supongamos que $\phi(s)=a+b$. Luego, por ser $\phi$ un automorfismo existen $x,y\in M$ tales que $\phi(s)=\phi(x)+\phi(y)=\phi(x+y)$ y así $s=x+y$. Dado que $s\in X$ tenemos que x o y son unidades, supongamos que x lo es. Entonces por la observación inicial $\phi(x)=a$ es una unidad. De modo que $\phi(s)\in X$. Supongamos ahora que $\phi(s)=t\in X$. Por la observación inicial tenemos que s no es unidad. Supongamos que $s=a+b$. Entonces $t=\phi(s)=\phi(a)+\phi(b)$, por lo que $\phi(a)$ o $\phi(b)$ es unidad, supongamos que $\phi(a)$ lo es. Por tanto a es unidad y así $s\in X$. Concluimos que el automorfismo $\phi$ permuta los elementos irreducibles de M y por lo tanto es isomorfo a un subgrupo del grupo de permutaciones de elementos irreducibles. 

\end{proof}

\begin{Def}\label{res:3.5}
Sea M un monoide y S un M-conjunto. Diremos que S es noetheriano si todo M-conjunto de S es finitamente generado.  
\end{Def}

Diremos que M es noetheriano si lo es con la representación regular. Es claro que una unión finita de M-conjuntos noetherianos es noetheriano, que un M-subconjunto de un M-conjunto noetheriano es nuevamente noetheriano, y que la imagen de un M-conjunto noetheriano es noetheriano. Se sigue que si M es un monoide noetheriano, entonces un M-conjunto es noetheriano si y solo si es finitamente generado como M-conjunto.

\begin{prop}\label{res:3.6}
Sea M un monoide y S un M-conjunto. Entonces las siguientes condiciones son equivalentes: \\ 
1. Todo M-subconjunto de S es finitamente generado, es decir, S es noetheriano. \\ 
2. Toda cadena ascendente $S_{1} \subset S_{2} \subset ...$ de M-conjuntos de S es eventualmente constante. \\ 
3. Toda sucesión $(s_{1},s_{2},...)$ en S contiene una subsucesión creciente. \\ 
4. Todo subconjunto no vacío de S contiene un elemento minimal, y solo hay finitas clases de equivalencia de tales elementos para la relación de equivalencia $\sim$.  \\ 
5. Todo conjunto no vacío de M-subconjuntos de S tiene un elemento maximal.  \\ 
6. El cociente de S por la relación de congruencia $\sim$ es noetheriano.  
\end{prop}

\begin{proof}
La equivalencia de (1),(2) y (5) se puede probar de la misma forma que en el caso de módulos sobre anillos. Supongamos que se cumple (2) y sea $(s_{1},s_{2},...)$ una sucesión en S. Para cada i sea $(s_{i})$ el M-subconjunto de S generado por $s_{i}$ y $S_{i}=(s_{1})\cup ... \cup (s_{i})$. Entonces $S_{1} \subset S_{2}  \subset ...$, así que por (2), existe algún $N\in \bb{N}$ tal que $S_{j}=S_{N}$ para todo $j\geq N$. Entonces para todo $j \geq N$ existe $i\leq N$ tal que $s_{j}\geq s_{i}$. Ya que hay infinitos j y finitos i, tiene que existir un $i\leq N$ y una sucesión $j_{1} < j_{2} < ...$ tal que $s_{i} \leq s_{j_{k}}$ para todo k. Entonces, reemplazando $(s_{1},s_{2},...)$ por la subsucesión $(s_{i},s_{j_{1}},s_{j_{2}},...)$, podemos asumir que $s_{i}\leq s_{j}$ para $j>1$. Repitiendo este proceso, podemos asegurar que $s_{2} \leq s_{j}$ para $j>2$ y así sucesivamente. \\ 
Para probar que (3) implica (4) primero observemos que cualquier sucesión decreciente en S es eventualmente una sola clase de equivalencia bajo $\sim$. De hecho, (3) implica que existe una sucesión creciente $(i_{j}:j\in \bb{Z}^{+})$ tal que $(s_{i_{1}},s_{i_{2}},...)$ es creciente. Entonces todo los $s_{i}$ con $i\geq i_{1}$ son equivalentes. De hecho, si $i\geq i_{1}$, elegimos j tal que $i_{j}\geq i$, y entonces
$$s_{i_{1}} \geq s_{i} \geq s_{i_{j}} \geq s_{i_{1}}.$$
Ahora si T es un subconjunto no vacío de S, tomamos un elemento $t_{1}$ de T. Si $t_{1}$ es M-minimal, no hay nada que hacer; si no, entonces existe $t_{2}\in T$ tal que $t_{2}\leq t_{1}$ y $t_{2} \ngeq t_{1}$. Si $t_{2}$ es M-minimal, está hecho, y si no, existe $t_{3}$ con $t_{3}\leq t_{2}$ y $t_{3} \ngeq t_{2}$. Continuando de esta forma, podemos encontrar una sucesión $(t_{1},t_{2},...)$ de elementos de T tal que $t_{i} \ngeq t_{i-1}$ para todo i. Como hemos visto, tal sucesión debe terminar, por lo que podemos encontrar un elemento M-minimal de T. Si hubiera un número infinito de clases de equivalencia de tales elementos minimales, podríamos encontrar una sucesión infinita $(s_{1},s_{2},...)$ de elementos todos pertenecientes a distintas clases de equivalencia, y por (3) tal sucesión debería contener una subsucesión creciente. Pero entonces $s_{1}\leq s_{2}$ y $s_{1} \nsim s_{2}$, contradiciendo la minimalidad de $s_{2}$, lo cual prueba (4). \\
Supongamos que (4) se cumple y sea T un M-subconjunto de S. Por (4), existe un conjunto finito $T^{\prime}$ de elementos minimales de T tal que todo elemento minimal es equivalente a algún elemento de $T^{\prime}$. Sea t un elemento arbitrario de T y $T_{t}:=\set{s\in T:s\leq t}$. Entonces $T_{t}$ no es vacío y así por (4) contiene un elemento minimal $t^{\prime \prime}$. Notemos que si $t^{\prime \prime \prime}\in T$ y $t^{\prime \prime \prime}\leq t^{\prime \prime}$, entonces $t^{\prime \prime \prime}\in T_{t}$ y entonces $t^{\prime \prime \prime} \sim t^{\prime \prime}$, por la minimalidad de $t^{\prime \prime}$. Entonces $t^{\prime \prime}$ es de hecho un elemento minimal de T, y entonces es equivalente a algún elemento $t^{\prime}\in T^{\prime}$. Ya que $t^{\prime \prime} \in T_{t}$, $t=m+t^{\prime \prime}= m^{\prime}+t^{\prime}$ para algunos $m,m^{\prime}\in M$. Entonces el conjunto finito $T^{\prime}$ genera a T. Esto prueba (1). \\
Sea $\pi:S \to S/ \sim$ la proyección natural y $S^{\prime}$ un M-subconjunto de S. Notemos que si $s^{\prime}\in S^{\prime}$ y $s\in S$ y $s\sim s^{\prime}$, entonces $s\in S^{\prime}$. De modo que $S^{\prime}=\pi^{-1}(\pi(S))$, por lo que $\pi$ induce una biyección entre la familia de M-subconjuntos de S y los de $S/ \sim$. Entonces (5) se cumple para S si y solo si se cumple para $S/ \sim$.  
\end{proof}

\begin{cor}\label{res:3.7}
Todo monoide \textit{dull} es noetheriano, y un monoide es noetheriano si y solo si su \textit{sharp}erización es noetheriano.  
\end{cor}

\begin{lem}\label{res:3.8}
Si P y Q son noetherianos, entonces $P \oplus Q$ es noetheriano. En particular si Q es noetheriano, entonces $Q \oplus \bb{N}$ es también noetheriano. 
\end{lem}

\begin{proof}
Sea $(p_{i},q_{i})$ una sucesión en $P\oplus Q$. Dado que P es noetheriano, existe una sucesión creciente $(n)$ en $\bb{N}$ tal que la subsucesión $(p_{i_{n}})$ es creciente. De manera análoga si Q es noetheriano tenemos que existe $(m)$ en $\bb{N}$ tal que $(q_{i_{n_{m}}})$ es creciente, por lo que $(p_{i_{n_{m}}},q_{i_{n_{m}}})$ también es creciente. 
\end{proof}

\begin{lem}\label{res:3.9}
Toda relación de congruencia en un monoide finitamente generado es finitamente generado (como relación de congruencia). 
\end{lem}

\begin{proof}
La demostración de este lema se puede consultar en ~\cite{Ogus:2018aa}. 
\end{proof}

\begin{Theorem}\label{res:3.10}
Un monoide finitamente generado es noetheriano y finitamente presentado. Por el contrario, un monoide noetheriano, \textit{sharp} y cuasi \textit{integral} es finitamente generado.   
\end{Theorem}

\begin{proof}
Es inmediato ver que si $M \to M^{\prime}$ es sobre y M es noetheriano, entonces $M^{\prime}$ es noetheriano. Entonces se sigue del Lema anterior y por inducción que todo monoide finitamente generado es noetheriano. El hecho de que un monoide finitamente generado es finitamente representado se sigue del lema anterior. \\
Por otro lado, supongamos que M es \textit{sharp}, cuasi \textit{integral} y noetheriano. Podemos asumir que M no es el monoide cero. Por (4) de la proposición~\ref{res:3.6} al subconjunto no vacío $M^{+}$ de M, podemos ver que el conjunto S de elementos minimales de $M^{+}$ es finito y no vacío. Más aún, todo elemento de M puede ser escrito como suma de elementos de S. Entonces M es finitamente generado. 
\end{proof}

\begin{cor}\label{res:3.11}
Un monoide M cuasi \textit{integral} es noetheriano si y solo si $\overline{M}$ es finitamente generado. 
\end{cor}

\begin{proof}
Supongamos que M es noetheriano. Entonces por el corolario~\ref{res:3.7} tenemos que $\overline{M}$ es noetheriano. Dado que M es cuasi \textit{integral} también $\overline{M}$ lo es. De igual forma recordemos que $\overline{M}$ siempre es \textit{sharp}. Por el Teorema~\ref{res:3.10} se sigue que $\overline{M}$ es finitamente generado. Supongamos ahora que $\overline{M}$ es finitamente generado. Por el Teorema~\ref{res:3.10} se tiene que $\overline{M}$ es noetheriano y finitamente representado. Por el corolario~\ref{res:3.7} se sigue que M es noetheriano. 
\end{proof}

\begin{cor}\label{res:3.12}
Si M es un monoide finitamente generado, cualquier M-subconjunto de un M-conjunto finitamente generado es finitamente generado, y de hecho es generado por un conjunto finito de elementos minimales. 
\end{cor}

\begin{prop}\label{res:3.13}
Sea Q un monoide y S un Q-conjunto finitamente generado. Supongamos que $qs \neq s$ para todo $s\in S$ y para todo $q\in Q^{+}$. Entonces S es generado por $S \setminus Q^{+}S$ . En particular, si $Q^{+}S=S$, entonces de hecho $S=\emptyset$. 
\end{prop}

\begin{proof}
Sea $\set{s_{1},...,s_{n}}$ un conjunto finito de generadores para S. Después de omitir algunos elementos, podemos suponer que ningún subconjunto propio genera a S. Entonces es suficiente probar que cada $s_{i}\notin Q^{+}S$. Si $s_{n}\in Q^{+}S$, entonces existe $q\in Q^{+}$ y $s\in S$ tal que $s_{n}=sq$, y dado que $\set{s_{1},...,s_{n}}$ genera a S, existen $i\leq n$ y $q^{\prime}\in Q$ tal que $s=q^{\prime}s_{i}$. Entonces $s_{n}=qq^{\prime}s_{i}$, y si $i<n$ el conjunto $\set{s_{1},...,s_{n-1}}$ genera a S, contradiciendo la minimalidad de S. Se sigue que $i=n$ así que $s_{n}=qq^{\prime}s_{n}$. La hipótesis entonces implica que $qq^{\prime}$ es unidad contradiciendo que $q\in Q^{+}$.  
\end{proof}

\begin{Def}\label{res:3.14}
Un ideal propio $\fk{q} \subset M$ es primario si $ax\notin \fk{q}$ para cualquier $a\in M\setminus \sqrt{\fk{q}}$ y $x\notin \fk{q}$.  
\end{Def}

\begin{prop}\label{res:3.15}
Sea M un monoide noetheriano. \\
1. Si $\fk{q}$ es un ideal primario de M, entonces su radical $\sqrt{\fk{q}}$ es primo, y existe un elemento x de M tal que $\sqrt{\fk{q}}=\set{a:ax\in \fk{q}}$. \\ 
2. Todo ideal de M puede ser escrito como la intersección de un número finito de ideales primarios. \\ 
3. Si $K=\fk{q}_{1} \cap ... \cap \fk{q}_{n}$ donde cada $\fk{q}_{i}$ es un ideal primario de M, sea $F_{i}:=M \setminus \sqrt{\fk{q}_{i}}$ y $F=F_{1} \cap ... \cap F_{n}$. Entonces $M\setminus K$ es estable bajo la multiplicación por F. 
\end{prop}

\begin{proof}
Supongamos que $\fk{p}$ es primario, y $a\notin \sqrt{\fk{p}}$ pero $ax\in \sqrt{\fk{p}}$. Entonces existe $n\in \bb{N}$ tal que $a^{n}x^{n}\in \fk{p}$, y como $\fk{q}$ es primario, se sigue que $x^{n}\in \fk{q}$ y por tanto $x\in \sqrt{\fk{q}}$. Entonces $\sqrt{\fk{q}}$ es primo. Para cada $x\in M\setminus \fk{q}$, sea $K_{x}:=\set{a:ax\in \fk{q}}$. Entonces $K_{x} \subset \sqrt{\fk{q}}$. Más aún, como M es noetheriano, existe $x\in M\setminus \fk{q}$ tal que $K_{x}$ no está propiamente contenido en en cualquier $K_{x^{\prime}}$. Veamos que de hecho $K_{x}=\sqrt{\fk{q}}$. De hecho, si $b\in \sqrt{\fk{q}}$, entonces alguna potencia de b pertenece a $\sqrt{\fk{q}}$, y entonces existe un $n\in \bb{N}$ tal que $b^{n}x\notin \fk{q}$ pero $b^{n+1}x\in \fk{q}$. Sea $x^{\prime}=b^{n}x$. Entonces $K_{x} \subset K_{x^{\prime}}$, y por maximalidad $K_{x}=K_{x^{\prime}}$. Ya que $b\in K_{x^{\prime}}$, se sigue que $b\in K_{x}$, completando la demostración de (1). \\
Si (2) no se cumpliera, podríamos encontrar un ideal de M que es maximal entre todos los ideales que no admiten una descomposición primaria. Tal ideal K debería ser necesariamente propio. Para $a\in  M \setminus \sqrt{K}$, sea $K_{n}=\set{x:a^{n}x\in K}$. Entonces $K \subset K_{1} \subset K_{2} \subset ...$, por lo que existe un $n\in \bb{N}$ tal que $K_{n}=K_{n+1}$. Notemos que si $x\in K_{1} \cap (a^{n})$, podemos escribir $x=a^{n}y\in K$ para algún $y\in M$, y $ax=a^{n+1}y$. Entonces $y\in K_{n+1}=K_{n}$, por lo que $x=a^{n}y\in K$. Esto implica que $K_{1} \cap (a^{n}) \subset K$, y entonces $K=K_{1} \cap K^{\prime}$ con $K^{\prime}=K \cup (a^{n})$. Ya que $a\notin \sqrt{K}$, $K^{\prime}$ estrictamente contiene a K y entonces admite una descomposición primaria. Entonces $K_{1}$ no admite una descomposición, y entonces $K_{1}=K$. Esto significa que $x\in K$ para cualquier $ax\in K$, suponiendo como antes que $a\notin \sqrt{K}$. Hemos demostrado de hecho que K es primario, otra contradicción. \\
Para (3) supongamos que $fm\in K$, donde $f\in F$ y $m\in K$. Entonces $fm\in \fk{q}_{i}$ para todo i y ya que $f\notin \sqrt{\fk{q}_{i}}$, se sigue que $m\in \fk{q}_{i}$ para todo i, por lo que $m\in K$. 
 \end{proof}

\begin{Def}\label{res:3.16}
Un homomorfismo de monoides es exacto si el siguiente diagrama
$$\xymatrix{
P \ar[d] \ar[r]^{\theta} & Q \ar[d] \\
P^{gp} \ar[r]_{\theta^{gp}} & Q^{gp}
}
$$
es cartesiano. 
\end{Def}

\begin{prop}\label{res:3.17}
Los siguientes enunciados se cumplen. \\ 
1. Si M es un monoide, el mapeo diagonal $\Delta_{M}:M \to M\times M$ es exacto si y solo si M es \textit{integral}. \\ 
2. Un homomorfismo exacto $\theta:P \to Q$ es local, y es inyectivo si P es \textit{sharp}. \\ 
3. Sea $\theta:P \to Q$ un homomorfismo de monoides \textit{integral}es. Entonces $\theta$ es exacto si y solo si para cualquier $p_{1},p_{2} \in P$, $\theta(p_{2}) \geq \theta(p_{1})$ implica que $p_{2}\geq p_{1}$. \\ 
4. En la categoría de monoides \textit{integral}es, el pullback de un homomorfismo exacto es exacto. \\ 
5. Sea P un submonoide de un monoide \textit{integral} Q. Entonces la inclusión $i:P \to Q$ exacta si y solo si $Q \setminus P$ es estable bajo la acción de P en Q. 
\end{prop}

\begin{proof}
Los detalles de la demostración se pueden encontrar en ~\cite{Ogus:2018aa}. 
\end{proof}

\begin{lem}\label{res:3.18}
Sea P un submonoide exacto de un monoide fino \textit{sharp} Q y $\Sigma$ un subconjunto no vacío de $P^{gp}$. Entonces un elemento s de $\Sigma$ es P-minimal in $\Sigma$ si y solo si $\theta^{gp}(s)$ es Q-minimal en $\theta^{gp}(\Sigma)$.    
\end{lem}

\begin{proof}
Podemos suponer sin pérdida de generalidad que P y Q son \textit{sharp}. Entonces $\theta^{gp}$ es inyectivo, ya que $\theta$ es exacto. Sea s un elemento de $\Sigma$. Es claro que s es P-minimal en $\Sigma$ si $\theta(s)$ es Q-minimal en $\theta^{gp}(\Sigma)$. Por otro lado, si $s^{\prime}\in \Sigma$ y $\theta^{gp}(s^{\prime})< \theta^{gp}(s)$ entonces $q:=\theta^{gp}(s)-\theta^{gp}(s^{\prime})\in Q^{+}$ y se sigue de la exactitud de $\theta$ que $p=s-s^{\prime}\in P^{+}$. Entonces S no puede ser P-minimal en S.
\end{proof}

\begin{lem}\label{res:3.19}
Si $\theta: P \to Q$ es un homomorfismo exacto de monoides \textit{integral}es y Q es saturado, entonces P también es saturado.  
\end{lem}

\begin{proof}
Supongamos que $x\in P^{gp}$ y $nx\in P$. Entonces $n\theta(x)=\theta(nx)\in Q$, y ya que Q es saturado, se sigue que $\theta(x)\in Q$. Como Q es exacta concluimos que $x\in P$, es decir, que P es saturado.  
\end{proof}

\begin{Theorem}\label{res:3.20}
1. Sea $\theta:P \to Q$ un homomorfismo exacto de monoides \textit{integral}es. Supón que S es un P-subconjunto de $P^{gp}$ cuya imagen en $Q^{gp}$ está contenida en un Q-conjunto noetheriano. Entonces S es un P-conjunto noetheriano. \\ 
2. Todo submonoide exacto de un monoide fino (respectivamente saturado,tórico) es fino (respectivamente saturado, tórico). \\ 
3. Una cara de un monoide \textit{integral} es un submonoide exacto. Toda cara de un monoide afín es finitamente generado. \\ 
4. Toda localización de un monoide fino es saturado. \\ 
5. El ecualizador E de dos morfismos de monoids \textit{integral}es $\theta_{i}:P \to Q$ es un submonoide exacto de P. Si P es fino (saturado), entonces lo mismo es cierto para E. \\ 
6. El producto fibrado de unn par de monoides finitamente generados (fino, saturado) sobre un monoide \textit{integral} es finitamente generado (fino, saturado). \\ 
7. Si $P \to Q$ es un homomorfismo de monoides finos, entonces la relación de congruencia $E.=P\times_{Q}P$ es finitamente generada como un monoide y en particular como una congruencia en P. \\ 
8. Sea P y Q monoides. Si Q es fino y P es finitamente generado, entonces Hom(P,Q) es tamién fino. Si Q es saturado, Hom(P,Q) también es saturado.  
\end{Theorem}

\begin{proof}
Los detalles de la demostración se pueden encontrar en ~\cite{Ogus:2018aa}. 
\end{proof}

\begin{cor}\label{res:3.21}
Si P es un monoide \textit{integral} y E es una relación de congruencia en P, entonces P/E es \textit{integral} si y solo si $E \to P\times P$ es exacto. En particular, si P y P/E son finos, entonces E también es fino. 
\end{cor}

\begin{proof}
De hecho la relación de congruencia E determinada por un homomorfismo sobre $\pi:P \to Q$ de monoides \textit{integral}es es solo el ecualizador de los dos mapeos $\pi \circ p_{1},\pi \circ p_{2}$ y entonces por el teorema~\ref{res:3.20} tenemos que es un submonoide exacto de $P\times P$. Por otro lado, sea E una relación de congruencia en P la cual es un submonoide exacto de $P\times P$ y $\pi:P \to Q$ el mapeo cociente. Entonces $\pi$ es el coecualizador de los dos mapeos $E \to P$. Para probar que Q es \textit{integral}, sean $q_{1},q_{2},q\in Q$ tales que $q_{1}+q=q_{2}+q$, y elegimos $p_{i}$ y p en P tales que $\pi(p_{i})=q_{i}$ y $\pi(p)=q$. Entonces $e:=(p_{1},p_{2})+(p,p)\in E$, por lo que se sigue que $(p_{1},p_{2})\in E^{gp}\cap (P\times P)$. Ya que E es un submonoide exacto, se sigue que $(p_{1},p_{2})\in E$ y entonces $\pi(p_{1})=\pi(p_{2})$. Si P y $P/E$ son finos, entonces E también es fino.  
\end{proof}

\begin{cor}\label{res:3.22}
Sea Q un monoide \textit{integral} finitamente generado, L un grupo abeliano finitamente generado, y $L \to Q^{gp}$ un homomorfismo. Entonces $L \times_{Q^{gp}} Q$ es un monoide finitamente generado. 
\end{cor}

\begin{proof}
Este corolario es un caso especial de (6) del teorema~\ref{res:3.20}. 
\end{proof}

\begin{obs}\label{res:3.23}
Sea Q un monoide \textit{integral}. Un subconjunto K de $Q^{gp}$ que es invariante bajo la acción de Q es llamado un ideal fraccional. El mapeo natural $\pi: Q \to \overline{Q}$ induce una biyección entre el conjunto de ideales fraccionales de Q y de $\overline{Q}$, y esta biyección toma ideales fraccionales finitamente generados y los envía a ideales fraccionales finitamente generados. 
\end{obs}

\begin{prop}\label{res:3.24}
Supongamos que $\theta: P  \to Q$ es un homomorfismo exacto de monoides finos, y que K es un ideal fraccional finitamente generado de Q. Entonces $J:= \theta^{-1}(K)$ es un ideal fraccional de P finitamente generado.  
\end{prop}

\begin{proof}
Este corolario es una consecuencia inmediata de (1) del Teorema~\ref{res:3.20}.
\end{proof}

\subsection{Dualidad}

\begin{Def}\label{res:3.25}
Sea Q un monoide. Definimos $H(Q):=Hom(Q,\bb{N})$.
\end{Def}

\begin{prop}\label{res:3.26}
Sea Q un monoide fino y F una cara de Q. Entonces existe un homomorfismo $h:Q \to \bb{N}$ tal que $h^{-1}(0)=F$.  
\end{prop}

\begin{proof}
Los detalles de la demostración se pueden encontrar en ~\cite{Ogus:2018aa}. 
\end{proof}

\begin{cor}\label{res:3.27}
Sea Q un monoide fino y x un elemento de $Q^{gp}$. Entonces $x\in Q^{sat}$ si y solo si $h^{gp}(x)\geq 0$ para todo $h\in H(Q)$. 
\end{cor}

\begin{proof}
Si $x\in Q^{sat}$ entonces $nx\in Q$ para algún $n\in \bb{Z}^{+}$, y entonces $h^{gp}(x)\geq 0$ para todo $h\in H(Q)$. Supongamos ahora que $h^{gp}(x)\geq 0$ para todo $h\in H(Q)$. Sea $Q^{\prime}$ el submonoide de $Q^{gp}$ generado por Q y $-x$ y escogemos un homomorfismo local $h:Q^{\prime} \to \bb{N}$. Entonces $h^{gp}(x)\geq 0$ y $h(-x)\geq 0$, por lo que de hecho, $h(-x)=0$. Ya que h es local, $-x \in Q^{\prime \ast}$. Entonces $x\in Q^{\prime}$, y escribimos $x=-xm+q$ con $m\in \bb{N}$ y $q\in Q$, por lo que vemos que $(m+1)x=q$ así que $x\in Q^{sat}$.  
\end{proof}

\begin{Theorem}\label{res:3.28}
Sea Q un monoide fino. \\
1. El monoide H(Q) es fino, saturado y \textit{sharp}. \\
2. El mapeo natural $H(Q)^{gp} \to Hom(Q^{gp},\bb{Z})$ se factoriza a través de un isomorfismo
$$\epsilon:H(Q)^{gp} \to Hom(\overline{Q}^{gp},\bb{Z}).$$
3. El mapeo de evaluación $ev: Q \to H(H(Q))$ se factoriza a través de un isomorfismo 
$$\overline{ev}:\overline{Q}^{sat} \to H(H(Q)).$$
Entonces el funtor H induce una involución contravariante de la categoría de monoides tóricos \textit{sharp}. En particular, esta categoría es dual en si misma.  
\end{Theorem}

\begin{proof}
El primer enunciado se sigue inmediatamente de (8) del teorema~\ref{res:3.20}. Ya que $H(Q) \to Hom(Q,\bb{Z})$ es inyectivo, también lo es el mapeo $H(Q)^{gp} \to Hom(Q,\bb{Z})$. Cualquier elemento h de H(Q) necesariamente aniquila a $Q^{\ast}$, así que este mapeo se factoriza a través de un mapeo $\epsilon$ siendo además inyectivo. Para probar que es sobre, sea $h:\overline{Q} \to \bb{N}$ un homomorfismo local y sea S un conjunto finito de generadores para $\overline{Q}$. Para $g\in Hom(\overline{Q},\bb{Z})$, existe $n\in \bb{Z}^{+}$ tal que $nh(s)\geq g(s)$ para cada $s\in S$. Entonces $nh(\overline{q})\geq g(\overline{q})$ para todo $\overline{q}\in \overline{Q}$, así que $h^{\prime}:=nh-g\in H(\overline{Q})$. Entonces $g=nh-h^{\prime}\in H(\overline{Q})^{gp} \cong H(Q)^{gp}$ como queríamos. \\
Ya que $H(H(Q))$ es fino,saturado y \textit{sharp}, ev se factoriza a través de un mapeo $\overline{ev}$. Sea $x_{1},x_{2}\in Q^{sat}$ con $ev(x_{1})=ev(x_{2})$, y sea $x_{1}-x_{2}\in Q^{gp}$ Entonces $h(x)=0$ para todo $h\in H(Q)$. Se sigue del corolario anterior que x y $-x$ pertenecen a $Q^{sat}$, así que $x\in (Q^{sat})^{\ast}$. Entonces $\overline{x_{1}}=\overline{x_{2}}$ en $\overline{Q^{sat}}$ lo cual prueba la inyectividad de $\overline{ev}$. Para la suprayectividad supongamos que $g\in H(H(Q))$. Ya que $Q^{gp}$ es un grupo finitamente generado, el mapeo que va de $Q^{gp}$ a su doble dual es sobre. Entonces existe un elemento $q\in Q^{gp}$ tal que $ev(q)=g$, es decir, tal que $h(q)=g(h)$ para todo $h\in H(Q)$. Entonces $H(q)\geq 0$ para todo $h\in H(Q)$, por lo que $q\in Q^{sat} $ como se quería.   
\end{proof}

\begin{cor}\label{res:3.29}
Sea Q un monoide fino. Un subconjunto S de Q es una cara de Q si y solo si existe un $h\in H(Q)$ tal que $S=h^{-1}(0)$. Para cada $S \subset Q$, sea $S^{\perp}$ el conjunto de $h\in H(Q)$ tal que $h(s)=0$ para todo $s\in S$, y para $T \subset H(Q)$, sea $T^{\perp}$ el conjunto de $q\in Q$ tal que $t(q)=0$ para todo $t\in T$. Entonces $F \mapsto F^{\perp}$ induce una biyección entre el conjunto de caras de Q y el conjunto de caras de H(Q), y $F=(F^{\perp})^{\perp}$ para cualquier cara de ambos.  
\end{cor}

\begin{proof}
El primer enunciado se sigue de la Proposición~\ref{res:3.26}. Es claro que si S es cualquier subconjunto de Q, entonces $S^{\perp}$ es una cara de H(Q) y $T^{\perp}$ es una cara de Q si T es cualquier subconjunto de H(Q). Más aún, $S_{2}^{\perp}\subset S_{1}^{\perp}$ si $S_{1} \subset S_{2}$ y $S \subset (S^{\perp})^{\perp}$. Sea F una cara de Q. Por la proposición~\ref{res:3.26} existe $h\in H(Q)$ tal que $h^{-1}(0)=F$. Entonces $h\in F^{\perp}$ y si $q\in (F^{\perp})^{\perp}$, $h(q)=0$ por lo que $q\in F$. Entonces $F=(F^{\perp})^{\perp}$ y entonces el mapeo $\perp$ que va del conjunto de caras de Q al conjunto de caras de H(Q) es inyectivo y el mapeo $\perp$ del conjunto de caras de H(Q) al conjunto de caras de Q es sobre. Entonces el mapeo que va del conjunto de caras de H(Q) al conjunto de caras de H(H(Q)) también es inyectivo. Por un corolario el mapeo $Q \to \overline{Q^{sat}}$ induce una biyección entre los correspondientes conjuntos de caras, y entonces por (3) del teorema~\ref{res:3.28}, $\overline{ev}$ identifica el conjunto de caras de Q con el conjunto de caras de H(H(Q)). Se sigue que $\perp$ es biyectivo. 
\end{proof}

\begin{cor}\label{res:3.30}
Si Q es un monoide fino entonces $Q^{sat}$ es fino. De hecho, la acción de Q en $Q^{sat}$ definida por el morfismo $Q \to Q^{sat}$ hace a $Q^{sat}$ un Q-conjunto finitamente generado.  
\end{cor}

\begin{proof}
Ya que $(Q^{sat})^{\ast}$ está contenido en $Q^{gp}$, es un grupo abeliano finitamente generado. El teorema~\ref{res:3.28} implica que $\overline{Q^{sat}}$ es fino, y ya que $Q^{sat}$ es \textit{integral}, se sigue de la proposición~\ref{res:3.1} que $Q^{sat}$ es finitamente generado, y entonces fino. Elegimos un conjunto finito de generadores de T de $Q^{sat}$ como monoide, y para todo $t\in T$, elegimos $n_{t}\in \bb{N}^{+}$ tal que $n_{t}t\in Q$. Entonces $\set{\sum j_{t}t:j_{t}\leq n_{t}}$ genera a $Q^{sat}$ como Q-conjunto.    
\end{proof}

\begin{cor}\label{res:3.31}
Si $\pi:Q^{\prime} \to Q$ es un homomorfismo sobre de monoides, entonces $H(\pi): H(Q) \to H(Q^{\prime})$ es inyectivo y exacto.  
\end{cor}

\begin{proof}
Es claro que $H(\pi)$ es inyectiva si $\pi$ es sobre. Más aún, por (2) del teorema~\ref{res:3.28}, podemos ver un elemento h de H(Q) como un homomorfismo $Q \to \bb{Z}$, y vemos que $h\in H(Q)$ si y solo si $h\circ \pi \in H(Q^{\prime})$.  
\end{proof}

\begin{cor}\label{res:3.32}
Sea Q un monoide fino \textit{sharp}. Entonces Q es isomorfo a un submonoide de $\bb{N}^{r} \otimes T$ para algún $r\in \bb{N}$ y algún grupo finito T. Si $Q^{gp}$ es libre de torsión, podemos tomar $T=0$. Si Q es saturado, entonces es isomorfo a un submonoide exacto de algún $\bb{N}^{r}$.  
\end{cor}

\begin{proof}
Supongamos que P es fino y \textit{sharp}. Por (8) del teorema~\ref{res:3.20}, el monoide $P:=H(Q)$ es fino y \textit{sharp}, y por (2) del teorema~\ref{res:3.28}, $P^{gp}:=H(Q)^{gp} \cong Hom(P^{gp},\bb{Z})$. Entonces $Hom(P^{gp},\bb{Z})\cong H(P)^{gp}$ es un grupo libre finitamente generado. Se sigue que el kernel del mapeo natural
$$Q^{gp} \to Hom(P^{gp},\bb{Z}) \cong H(P)^{gp}$$
es el subgrupo de torsión T de $Q^{gp}$. Eligiendo una descomposición,
$Q^{gp} \cong H(P)^{gp} \oplus T$
y un morfismo sobre $\bb{N}^{r} \to P$. Por el corolario~\ref{res:3.31}, H(P) es entonces un submonoide exacto de $H(\bb{N}^{r}) \cong \bb{N}^{r}.$ Entonces la inyección natural,
$$Q \to Q^{gp} \to H(P)^{gp}\oplus T$$
se factoriza a través de la inclusión $H(P)\oplus T \subset H(P)^{gp}\oplus T$, y Q es un submonoide de $\bb{N}^{r}\oplus T$ ya que $H(P) \subset \bb{N}^{r}$. Más aún, por (3) del teorema~\ref{res:3.28}, el mapeo natural $Q \to H(P)$ se factoriza a través de un isomorfismo $\overline{Q^{sat}} \to H(P)$. Entonces si Q también  es saturado, $Q\cong H(P)$, un submonoide exacto de $\bb{N}^{r}$. 
\end{proof}

\begin{obs}\label{res:3.33}
Recordemos que el interior de un monoide es el complemento de todas sus subcaras propias. Notemos que si Q es fino, entonces un elemento $h\in H(Q)$ está en el interior de H(Q) si y solo si $h:Q \to \bb{N}$ es un morfismo local. De hecho, por el corolario~\ref{res:3.29}, h está en el interior de H(Q) si y solo si $h^{\perp}$ no contiene ninguna cara no trivial de Q, es decir, si y solo si $h^{\perp}=Q^{\ast}$ 
\end{obs}

\begin{prop}\label{res:3.34}
Sea Q un monoide \textit{sharp} fino, d el rango de $Q^{gp}$ y $h:Q \to \bb{N}$ un homomorfismo local. Para cada número real r definimos
$$B_{h}(r):=\set{q\in Q:h(q)<r}.$$
Entonces existen constantes reales positivas c y C tales que, para todo $r\gg 0$,
$$cr^{d}< \sharp B_{h}(r)<Cr^{d}.$$
\end{prop}

\begin{proof}
El subgrupo de torsión T de $Q^{gp}$ es finito, y el cociente $Q^{gp}/T$ es un grupo abeliano libre de rango d. Por el teorema~\ref{res:3.28}, H(Q) es finitamente generado y \textit{sharp}, y entonces por la proposición~\ref{res:3.3} este tiene un único conjunto de generadores minimal $\set{h_{1},...,h_{m}}$. Ya que h es local, pertenece al interior de H(Q), por una observación anterior. Entonces la cara de H(Q) generada por h es todo H(Q), y en particular contiene a cada $h_{i}$. Por la proposición~\ref{res:2.51}, esto quiere decir que para cada i existe un $n_{i}$ tal que $n_{i}h\geq h_{i}$ en H(Q). Elijamos $n\geq n_{i}$ para todo i. Entonces $B_{h}(r) \subset \cap_{i} B_{h_{i}}(nr)$ para todo $r\in \bb{R}^{+}$. Ya que Q es \textit{sharp}, (2) del teorema~\ref{res:3.28} implica que $H(Q)^{gp} \cong Hom(Q^{gp},\bb{Z})$ y consecuentemente $\set{h_{1},...,h_{m}}$ genera al $\bb{Q}$-espacio vectorial $Hom(Q^{gp},\bb{Q})$. Ya que este espacio tiene dimensión d, el conjunto $\set{h_{1},...,h_{m}}$ contiene una base, la cual podemos suponer es $\set{h_{1},...,h_{d}}$. Entonces el mapeo $Q^{gp}\otimes \bb{Q} \to \bb{Q}^{d}:x\otimes 1 \mapsto (h_{1}(x),...,h_{d}(x))$  es un isomorfismo e induce una inyección entre la imagen $Q^{\prime}$ de Q en $Q^{gp}/T$ a $\bb{N}^{d}$. Si $q\in B_{h}(r)$, su imagen en $\bb{N}^{d}$ vive en $\set{(I_{1},...,I_{d}):0\leq I_{i} \leq nr}$, un conjunto de cardinalidad $(1+[nr])^{d}$. Ya que la cardinalidad de las fibras del mapeo $Q \to Q^{\prime}$ está acotada por el orden de t en T, se sigue que $\sharp B_{h}(r) \leq t(1+nr)^{d}$. Entonces si $C:= t(1+n)^{d}$ y $r\geq 1$, la cardinalidad de $B_{h}(r)$ está acotada por $Cr^{d}$. \\ 
Por otro lado, cualquier conjunto de generadores para el monoide Q también genera al d-dimensional $\bb{Q}$-espacio vectorial $\bb{Q}\otimes Q^{gp}$, y entonces contiene un subconjunto $\set{q_{1},...,q_{d}}$ cuya imagen forma una base. Entonces el homomorfismo $\theta:\bb{N}^{d} \to Q$ enviando $(I_{1},...,I_{d})$ a $\sum I_{i}q_{i}$ es inyectivo. Sea $n:=max \set{h_{1}(q_{1}),...,h_{d}(q_{d})}$ y sea $c:=(2nd)^{-d}$. Entonces $\sharp B_{h}(r) \geq cr^{d}$ para $r \geq nd$. De hecho, si $r \geq nd$, sea $m:=[r/nd]$, y notemos que 
$$r/nd \geq m \geq r/nd-1 \geq r/2nd.$$
El conjunto $\set{I:I_{i}\leq m,i=1,...,d}$ tiene cardinalidad $(1+m)^{d}$ y su imagen en Q también tiene cardinalidad $(1+m)^{d} \geq m^{d} \geq cr^{d}$ y está contenido en $B_{h}(r)$. 
\end{proof}

\subsection{Monoides valuativos y valuaciones}
Recordemos que un monoide Q es valuativo si para todo $x\in Q^{gp}$ se tiene que x o $-x$ están en Q. \\

\begin{Def}\label{res:3.35}
Si L es un grupo abeliano y $Q_{1},Q_{2}$ son submonoides de L, diremos que $Q_{2}$ domina a $Q_{1}$ si $Q_{1} \subset Q_{2}$ y el mapeo inclusión es local. 
\end{Def}

Esto define un orden el conjunto de submonoides de L.

\begin{prop}\label{res:3.36}
Sea L un grupo abeliano. Entonces un submonoide Q de L es un elemento maximal para la relación de dominación si y solo si es valuativo y $Q^{gp}=L$. Todo submonoide de L está dominado por tal monoide. Si l es finitamente generado, entonces todo submonoide finitamente generado de L es dominado por un monoide valuativo finitamente generado Q con $Q^{gp}=L$. 
\end{prop}

\begin{proof}
Supongamos primero que Q es un submonoide de L y es maximal bajo dominación. Entonces $Q \subset Q^{sat} \subset L$ y por (4) de una proposición anterior, el homomorfismo $Q \to Q^{sat}$ es local. Entonces $Q^{sat}$ domina a Q, y entonces $Q=Q^{sat}$. Ahora para probar que Q es valuativo, sea $x\in L$, y consideremos el submonoide P de L generado por Q y x. Si P domina a Q, entonces $P=Q$ y claramente $x\in Q$. Si no, existe $q\in Q$ y $p\in P$ con $p+q=0$. Escribimos $p=q^{\prime}+nx$ con $n>0$. Entonces $q+q^{\prime}+nx=0$. Esto muestra que $-nx\in Q$, y entonces $-x\in Q^{sat}=Q$. Ya que x fue un elemento arbitrario de L, se sigue que $Q^{gp}=L$. Si P es un submonoide de L dominando a Q y $p\in P$, entonces p o $-p$ pertenecen a Q. Pero si $-p$ pertenece a Q, entonces al ser $-p$ una unidad en P también lo será en Q, así que p pertenece a Q en cualquier caso. Esto muestra que $Q=P$. \\ 
Supongamos que P es un submonoide de L. Entonces el conjunto $S_{P}$ de submonoides de L que dominan a P es no vacío, y la unión de cualquier cadena en $S_{P}$ contiene un elemento maximal. Si L y P son finitamente generados, entonces por un corolario
$$P^{v}:= \set{w\in Hom(L,\bb{Z}):<w,p> \geq 0,p\in P}$$
es un submonoide finitamente generado de $Hom(L,\bb{Z})$. Sea $w:L \to \bb{Z}$ un elemento en el interior de $P^{v}$. Entonces w define un homomorfismo local $P \to \bb{N}$. Ya que $\bb{N}$ es valuativo, $Q:=w^{-1}(\bb{N})$ es valuativo y claramente $Q^{gp}=L$. Más aún, $Q=\set{w}^{v}$ y entonces es finitamente generado. Ya que w se factoriza a través del homomorfismo $P \to Q \to \bb{N}$ y es local, se sigue que $P \to Q$ también es local, así que Q domina a P.
\end{proof}

\begin{prop}\label{res:3.37}
Sea Q un monoide no cero, \textit{integral} y \textit{sharp}. Entonces las siguientes condiciones son equivalentes: \\
1. Q es valuativo y finitamente generado. \\ 
2. Q es isomorfo a $\bb{N}$. \\ 
3. $dim(Q)=1$ y Q es saturado y finitamente generado. \\ 
4. Q es saturado y $Q^{gp}\cong \bb{Z}$. 
\end{prop}

\begin{proof}
Dado que Q es no cero y \textit{sharp}, su ideal maximal $Q^{+}$ es no vacío. Si Q es finitamente generado, entonces por la proposición~\ref{res:3.3}, Q es generado por los elementos minimales de $Q^{+}$. Si Q es valuativo la relación de orden en Q es un orden total y $Q^{+}$ tiene un solo elemento minimal q. Entonces existe un homomorfismo sobre $\bb{N} \to Q$, y ya que Q es \textit{integral} y \textit{sharp}, este homomorfismo es un isomorfismo. Entonces (1) implica (2). Es obvio que (2) implica (1) y  (3). Si (3) se cumple entonces por la proposición~\ref{res:2.56}, $Q^{gp}$ tiene rango uno. Ya que Q es \textit{sharp} y saturado, $Q^{gp}$ es libre de torsión, y entonces (4) se cumple. Si (4) es cierto, tomamos un isomorfismo  $\phi:Q^{gp} \to \bb{Z}$. Ya que $\phi$ es inyectivo, $\phi(q)\neq 0$ para todo $q\in Q^{+}$; elegimos $q\in Q^{+}$ con mínimo $|\phi(q)|$ y ajustamos el signo de $\phi$ tal que $n=\phi(q)>0$. Entonces se sigue por ser \textit{sharp} Q que $\phi$ mapea Q a $\bb{N}$. Si $x\in Q^{gp}$ es el elemento tal que $\phi(x)=1$, entonces $\phi(nx)=\phi(q)$, entonces $nx=q\in Q$. Ya que Q es saturado tenemos que $x\in Q$, y por la minimalidad de $|\phi(x)|$, necesariamente $n=1$. Entonces $\phi_{Q}:Q \to \bb{N}$ es inyectivo y suprayectivo, por tanto un isomorfismo.   
\end{proof}

\begin{cor}\label{res:3.38}
Todo monoide fino saturado de dimensión uno es isomorfo a $\Gamma \oplus \bb{N}$ para algún grupo abeliano finitamente generado $\Gamma$. Si Q es tórico entonces $Q\cong \bb{Z}^{r}\oplus \bb{N}$ para algún número natural r.  
\end{cor}

\begin{proof}
Si Q es tal monoide, el resultado anterior implica que existe un isomorfismo $\phi:\overline{Q}\to \bb{N}$. Sea $q\in Q$ con $\phi(q)=1$. Entonces el homomorfismo $Q^{\ast}\oplus \bb{N} \to Q$ que envía 1 a q es un isomorfismo. Ya que Q es fino, $Q^{\ast}$ es un grupo abeliano finitamente generado y si Q es tórico, entonces $Q^{\ast}$ es libre de torsión, por tanto libre.  
\end{proof}

\begin{cor}\label{res:3.39}
Sea $\fk{p}$ un ideal primo de altura uno de un monoide fino Q. Entonces $Q^{sat}_{\fk{p}}$ es valuativo y existe un único epimorfismo
$$v_{\fk{p}}:Q^{gp} \to \bb{Z}$$
tal que $v_{\fk{p}}^{-1}(\bb{N}^{+})\cap Q =\fk{p}$. Más aún, $Q^{sat}_{\fk{p}}=\set{x\in Q^{gp}: v_{\fk{p}}(x)\geq 0}.$ 
\end{cor}

\begin{proof}
Sabemos que $Q^{gp}\cong (Q^{sat})^{gp}$, que $Q^{sat}$ es fino por el corolario~\ref{res:3.30} y que $\Spec(Q^{sat}) \to \Spec(Q)$ es un homeomorfismo. Entonces podemos reemplazar Q por $Q^{sat}$, y entonces asumir que Q es saturado. Entonces $Q_{\fk{p}}$ es saturado, así que $\overline{Q_{\fk{p}}}^{gp}$ es libre de torsión, y ya que $\fk{p}$ tiene altura uno, $dim(\overline{Q_{\fk{p}}})=1$. Por la proposición anterior, $\overline{Q_{\fk{p}}}$ es valuativo y entonces también lo es $Q_{\fk{p}}$. Más aún, existe un único isomorfismo $\overline{Q_{\fk{p}}}\cong \bb{N}$, y sea $v_{\fk{p}}$ la composición de $Q^{gp}\cong Q_{\fk{p}}^{gp} \to \overline{Q_{\fk{p}}}^{gp}$ con el isomorfismo inducido $\overline{Q_{\fk{p}}}^{gp} \to \bb{Z}$. Es claro que $v_{\fk{p}}^{-1}(\bb{N})=Q_{\fk{p}}$ y que $v_{\fk{p}}^{-1}(\bb{N}^{+})\cap Q=\fk{p}$. La unicidad de $v_{\fk{p}}$ es verificada fácilmente.   
\end{proof}

\begin{cor}\label{res:3.40}
Sea Q un monoide fino saturado. Entonces,
\begin{equation*}
\begin{split}
Q &= \set{x\in Q^{gp}:v_{\fk{p}}(x)\geq 0,ht(\fk{p})=1} \\
Q^{\ast} &= \set{x\in Q^{gp}:v_{\fk{p}}(x) = 0,ht(\fk{p})=1} 
\end{split}
\end{equation*}
En particular, Q es la intersección en $Q^{gp}$ del conjunto de todas sus localizaciones en sus ideales primos de altura uno.
\end{cor}

\begin{proof}
Sea $\fk{p}$ un ideal primo de altura uno y G su cara complementaria, una cara de Q. Entonces $\set{x\in Q^{gp}:v_{\fk{p}}(x)\geq 0}=Q_{G}$ y $\set{x\in Q^{gp}:v_{\fk{p}}(x)=0}=G^{gp}$. Entonces los enunciados (3) y (2), respectivamente implican los dos enunciados del corolario.  
\end{proof}

\begin{prop}\label{res:3.41}
Sea Q un monoide fino y $W_{Q}^{+}$ el monoide libre en el conjunto de ideales primos de altura uno de Q. (Este conjunto es finito.) Para $q\in Q$ sea
$$v(q)=\sum \set{v_{\fk{p}}(q) \fk{p}: ht(\fk{p})=1}\in W_{Q}^{+}.$$
Entonces $v:Q \to W_{Q}^{+}$ es un homomorfismo local. Más aún, $v(q_{1})=v(q_{2})$ si y solo si existe algún $n\in \bb{Z}^{+}$ tal que $n\overline{q_{1}}=n\overline{q_{2}}$ en $\overline{Q}$, y v es exacta si y solo si Q es saturado.  
\end{prop}

\begin{proof}
Es claro que $v:Q \to W_{Q}^{+}$ es un homomorfismo de monoides. Para ver que v es local, notemos que su objetivo es \textit{sharp} y por el corolario anterior, $v(q)=0$ implica que q es una unidad de $Q^{sat}$ y entonces también de Q. El mismo corolario también implica que v induce un homomorfismo exacto $Q^{sat} \to Q$. Entonces si $q_{1},q_{2}\in Q$ con $v(q_{1})=v(q_{2})$, tenemos que $q_{1}-q_{2}$ es una unidad en $Q^{sat}$ y así que existe $n\in \bb{Z}^{+}$ tal que $nq_{1}-nq_{2}\in Q^{\ast}$. Esto implica que $n\overline{q_{1}}=n\overline{q_{2}}$. Hemos visto ya que v es exacto si Q es saturado, y el converso se sigue del hecho de que un submonoide exacto de un monoide saturado es saturado.     
\end{proof}

\section{Variedades tóricas afines}

\subsection{Álgebras monoidales y esquemas monoidales} 

Sea R un anillo conmutativo fijo, usualmente el anillo de enteros $\bb{Z}$ o un campo, y sea $Alg_{R}$ la categoría de R-álgebras. 
\begin{Def}\label{res:4.1}
Sea Q un monoide. El álgebra R-monoidal de Q es la R-álgebra $R[Q]$, la cual al ser tratada como R-módulo es libre con base Q, equipada con la estructura de anillo única que hace al mapeo inclusión $e:Q \to R[Q]$ un homomorfismo de monoides entre Q y el monoide multiplicativo de $R[Q]$. 
\end{Def}

Si $p,q\in Q$ usaremos la notación aditiva para Q por lo que tenemos que $e(p+q)=e(p)e(q)$. Y escribiremos $e^{p}$ en lugar de $e(p)$. Por ejemplo $R[\bb{N}]$ es el álgebra polinomial $R[T]$ donde $T=e^{1}$. \\

Sea $\underline{A}_{m}$ el funtor que toma una R-álgebra y la envía a su monoide multiplicativo. Notemos que dar un homomorfismo de monoides $Q \to \underline{A}_{m}(A)$ es equivalente a dar un morfismo de R-álgebras $R[Q] \to A$. Por lo que $Hom(R[Q],A) \cong Hom(Q,\underline{A}_{m}(A))$. De este modo, el funtor $Q \to R[Q]$ es adjunto izquierdo del funtor $\underline{A}_{m}$.  

\begin{Def}\label{res:4.2}
Sea Q un monoide y A una R-álgebra. Un homomorfismo $Q \to \underline{A}_{m}(A)$ lo llamaremos un carácter A-valuado de Q. 
\end{Def}

\begin{Def}\label{res:4.3}
Una R-álgebra equipada con un carácter de Q es llamada una Q-álgebra. 
\end{Def}

Dado que el funtor $Q \mapsto R[Q]$ es una adjunta izquierda, este automaticamente conmuta con límites y colímites. Por ejemplo si Q es la suma amalgamada de $\theta_{i}:P \to Q_{i}$, entonces $R[Q]\cong R[Q_{1}]\otimes_{R[P]}R[Q_{2}]$. \\

El conjunto $\underline{A}_{Q}(A)$ de caracteres A-valuados de Q tiene una estructura natural de monoide:
\begin{equation*}
\begin{split}
\cdot: \underline{A}_{Q}(A)\times \underline{A}_{Q}(A) &\to \underline{A}_{Q}(A): (\theta_{1},\theta_{2}) \mapsto \theta_{1}\theta_{2} \\
1:Q &\to \underline{A}_{m}(A):a \mapsto 1 
\end{split}
\end{equation*}

Entonces podemos ver a $\underline{A}_{Q}$ como un funtor
\begin{equation*}
\begin{split}
\underline{A}_{Q}:Alg_{R} &\to Mon \\
A &\mapsto \underline{A}_{Q}(A).
\end{split}    
\end{equation*}
Y dado un morfismo de R-álgebras $f:A \to B$, le asociamos un morfismo

$$g:\underline{A}_{Q}(A) \to \underline{A}_{Q}(B):\phi:Q \to \underline{A}_{m} \mapsto f \circ \phi$$  
donde f lo consideramos como un morfismo entre sus respectivos monoides multiplicativos. 

El funtor $\underline{A}_{Q}$ es representado por el par $(R[Q],e)$ ya que la transformación natural
\begin{equation*}
\begin{split}
e:Hom(R[Q],-) &\to \underline{A}_{Q} \\
e_{A}:Hom(R[Q],A) &\to \underline{A}_{Q}(A) \\
\phi:R[Q] \to A &\mapsto e_{A}(\phi):Q \to \underline{A}_{m}(A)
\end{split}
\end{equation*}

es un isomorfismo natural por la adjunción de $Q \to R[Q]$. 

Utilizaremos la misma notación $\underline{A}_{Q}$ tanto para el funtor como para el esquema $\Spec{R[Q]}$ que representa.

\begin{Def}\label{res:4.4}
Un S-esquema es un esquema X con un morfismo de esquemas $X \to S$. 
\end{Def}

\begin{Def}\label{res:4.5}
Un objeto de monoide en una categoría C con productos finitos y objeto final S es una tripleta $(M,m,1)$, con M un objeto en C, $m:M\times M \to M$  y $1:S \to M$, que satisface que los siguientes diagramas conmuten: \\
i) Asociatividad
$$\xymatrix{
M\times M\times M \ar[d]_{m\times id}  \ar[r]^{id\times m} & M\times M \ar[d]^{m} \\
M\times M \ar[r]_{m} & M.
}
$$

ii) Identidad
$$\xymatrix{
S\times M= M\times S \ar[d]_{1\times id} \ar[dr] \ar[r]^{id\times 1} & M\times M \ar[d]^{m} \\
M\times M \ar[r]_{m} & M
}
$$

iii) Conmutatividad
$$\xymatrix{
M\times M \ar[dr]_{m} \ar[r]^{(p_{1},p_{2})} & M\times M \ar[d]_{m}\\
& M 
}
$$
\end{Def}

\begin{obs}\label{res:4.6}
Dar un objeto de monoide en C es equivalente a dar un funtor contravariante que va de C a la categoría de monoides tal que el funtor restringido a la categoría de conjuntos es representable. Para más detalles consulte ~\cite{Ogus:2018aa}. 
\end{obs}

\begin{Def}\label{res:4.7}
Un S-esquema monoidal es un objeto de monoide en la categoría de S-esquemas. 
\end{Def}

Cuando consideremos $S=\Spec{R}$, escribiremos R-esquema monoidal en lugar de $\Spec{R}$-esquema monoidal.

Consideremos el R-esquema afín $\underline{A}_{Q}=\Spec{R[Q]}$ y veamos que es un R-esquema monoidal. Consideremos los morfismos de R-álgebras:

\begin{equation*}
\begin{split}
m_{Q}^{\sharp}:R[Q] &\to R[Q] \otimes R[Q]:e^{p} \mapsto e^{p}\otimes e^{p} \\
1_{Q}^{ \sharp}:R[Q] &\to R: \sum_{q} a_{q}e^{q} \mapsto \sum_{q}a_{q} 
\end{split}
\end{equation*}

Dado que el funtor Spec es contravariente este induce morfismos de esquemas

$$m_{Q}:\underline{A}_{Q}\times \underline{A}_{Q} \to \underline{A}_{Q}, ~~~~ 1_{Q}:\Spec{R} \to \underline{A}_{Q}.$$

Tenemos que $\Spec{R}$ es un objeto final en la categoría de R-esquemas. Dado que $(id\otimes m_{Q}^{\sharp}) \circ m_{Q}^{\sharp}=(m_{Q}^{\sharp}\otimes id) \circ m_{Q}^{\sharp}$, esto implica que $m_{Q} \circ (id\times m_{Q}) =  m_{Q} \circ (m_{Q}\times id)$, por lo que se satisface i). De igual forma se puede ver que la condición ii) y iii) se satisface. Por tanto, $\underline{A}_{Q}=\Spec{R[Q]}$ es un R-esquema monoidal.

Entonces la sección identidad $1_{Q}$ del esquema monoidal $\underline{A}_{Q}$ está dada por el homomorfismo de monoides $Q \to \underline{A}_{m}(R):q \mapsto 1$.

\begin{Def}\label{res:4.8}
Se define el vértice de $\underline{A}_{Q}$ como la sección $v_{Q}$ de $\underline{A}_{Q}$ inducida por el morfismo $Q \to \underline{A}_{m}(R)$ que envía $q\in Q$ a 0 si $q\in Q^{+}$ y a 1 si $q\in Q^{\ast}$. 
\end{Def}

\begin{obs}\label{res:4.9}
El funtor $\underline{A}_{m}$ es isomorfo al funtor $\underline{A}_{\bb{N}}$. Sea A una R-álgebra, entonces
$$\underline{A}_{m}(A)=(A,\cdot,1) \cong Hom(\bb{N},A)$$
donde un elemento a corresponde con el homomorfismo de monoides $n \mapsto a^{n}$, y un morfismo $\bb{N} \to A:1 \to b$ lo envíamos al elemento b en $\bb{N}$. 
\end{obs}

La siguiente proposición muestra que Q puede ser recuperado usando el funtor $\underline{A}_{Q}$ con su estructura de esquema monoidal. 

\begin{prop}\label{res:4.10}
Supón que $\Spec{R}$ es conexo. Entonces el funtor 
$$Q \mapsto \underline{A}_{Q}$$
que va de la categoría de monoides a la de esquemas monoidales sobre R es completamente fiel:
$$Hom(P,Q) \cong Hom(\underline{A}_{Q},\underline{A}_{P})$$
\end{prop}

\begin{proof}
Sea $\theta^{\sharp}:\underline{A}_{Q} \to \underline{A}_{P}$ un homomorfismo de esquemas. Este corresponde con un morfismo de anillos $\theta:R[P] \to R[Q]$. Dado que $Q \to R[Q]$ es inyectivo, un morfismo $P\to Q$ es determinado por el homomorfismo correspondiente $R[P] \to R[Q]$. Esto prueba que el funtor es fiel. \\
Sea $\theta:R[P] \to R[Q]$ un homomorfismo de R-álgebras. Para cada $p\in P$ escribimos,
$$\theta(e^{p})=\sum_{q\in Q}a_{q}(p)e^{q}$$
con $a_{q}(p)\in R$. Decir que $\theta$ corresponde con un homomorfismo de monoides es equivalente a decir que los siguientes diagramas conmutan:
$$\sqrd{R[P]}{R[Q]}{R}{R}{\theta}{1_{Q}^{\sharp}}{id}{1_{P}^{\sharp}}$$
$$\sqrd{R[P]}{R[Q]}{R[Q]\otimes R[Q]}{R[P]\otimes R[P]}{\theta}{m_{Q}^{\sharp}}{\theta \otimes \theta}{m_{P}^{\sharp}}$$

El segundo diagrama dice que, para cada $p\in P$,

$$\sum_{q,q^{\prime}}a_{q}(p)a_{q^{\prime}}(p)e^{q}\otimes e^{q^{\prime}}= \sum_{q}a_{q}(p)e^{q}\otimes e^{q},$$
es decir, que $a_{q}(p)a_{q^{\prime}}(p)$ es igual a cero si $q\neq q^{\prime}$ y igual a $a_{q}(p)$ si $q=q^{\prime}$. El primer diagrama dice que $\sum_{q\in Q}a_{q}(p)=1$ para cada $p\in P$. Dado que $\Spec{R}$ es conexo, todo idempotente es 0 o 1, ya que $a_{q}(p)$ y $a_{q^{ \prime}}(p)$ son ortogonales si $q\neq q^{\prime}$ existe un único elemento $\beta(p)\in Q$ tal que $a_{q}(p)=0$ si $q\neq \beta(p)$ y $a_{q}(p)=1$ si $q=\beta(p)$. Entonces $\beta$ es una función $P \to Q$ tal que $\theta \circ e = e \circ \beta$. Ya que $\theta$ es un homomorfismo de monoides $\beta$ también lo es como se quería.   
\end{proof}

\begin{cor}\label{res:4.11}
Supon que el $\Spec{R}$ es conexo y Q es un monoide. \\ 
1. El monoide de caracteres de $\underline{A}_{Q}$, es decir, de morfismos $\underline{A}_{Q} \to \underline{A}_{m}$, es canonicamente isomorfo a Q. \\
2. El monoide de cocaracteres de $\underline{A}_{Q}$, es decir, de morfismos $\underline{A}_{m} \to \underline{A}_{Q}$, es canonicamente isomorfo a $H(Q):=Hom(Q,\bb{N})$.
\end{cor}

\begin{proof}
Notemos que,
\begin{equation*}
\begin{split}
Hom(\bb{N},Q) &\to Q \\
f:\bb{N} \to Q:1 \mapsto a &\mapsto a \\
g:\bb{N} \to Q:1 \mapsto a &\mapsfrom a
\end{split}
\end{equation*}
son homomorfismos de monoides e inversas una de la otra. Por lo que $Hom(\bb{N},Q) \cong Q$ como monoides. Tomando $P=\bb{N}$ en la proposición anterior tenemos que,
\begin{equation*}
\begin{split}
Q &\cong Hom(\bb{N},Q) \\
&\cong Hom(\underline{A}_{Q},\underline{A}_{\bb{N}}) \\
&\cong Hom(\underline{A}_{Q},\underline{A}_{m}).
\end{split}
\end{equation*}
El punto es simplemente tomar los monoides Q y $\bb{N}$ en ese orden en la proposición anterior. 
\end{proof}

Si Q es un monoide y A una R-álgebra, $\underline{A}_{Q^{gp}}(A)$ es precisamente el conjunto de elementos invertibles de $\underline{A}_{Q}(A)$, es decir, $\underline{A}_{Q^{gp}}(A)=\underline{A}_{Q}^{\ast}(A)$.

\subsection{Conjuntos de monoides y módulos de monoides}
\begin{Def}\label{res:4.12}
Sea E un R-módulo. Definimos VE como el funtor que toma una R-álgebra A y lo envía al conjunto de funciones R-lineales $E \to A$. 
\end{Def}
Este funtor es representado por el álgebra simétrica SE y el elemento universal de VE(SE) es la inclusión $E \to SE$. Para más detalles sobre esto consultar ~\cite{dieudonne1971elements}. Esto quiere decir,
$$VE(A)=Hom_{R}(E,A)=Hom_{Alg_{R}}(SE,A)$$
para todo R-álgebra A. El conjunto VE(A) tiene una estructura natural de A-módulo definiendo las operaciones puntualmente. \\
Esto define transformaciones naturales,
$$\mu:VE\times VE \to VE,~~~~ \psi:Id\times VE \to VE,$$
tal que para toda R-álgebra A,

\begin{equation*}
\begin{split}
\mu_{A}:VE(A)\times VE(A) &\to VE(A) \\
(f,g) &\mapsto fg
\end{split}
\end{equation*}
y
\begin{equation*}
\begin{split}
\psi_{A}:A\times VE(A) &\to VE(A) \\
(a,g) &\mapsto ag
\end{split}
\end{equation*}

donde Id es el funtor identidad en la categoría de R-álgebras. Notemos que este funtor es representado por la línea afín sobre R, es decir, por el $\Spec{R[x]}$. Esto es cierto al definir la siguiente transformación natural,
$$\phi:Id \to Hom_{Alg_{R}}(R[x],-)$$
tal que para todo $A\in Alg_{R}$,

\begin{equation*}
\begin{split}
\phi_{A}:A &\to Hom_{Alg_{R}}(R[x],A) \\
a &\mapsto \phi_{A}(a):R[x] \to A: \sum_{i=1}^{n}b_{i}x^{i} \mapsto \sum_{i=1}^{n}b_{i}a^{i}
\end{split}
\end{equation*}
con inversa,
\begin{equation*}
\begin{split}
\phi_{A}^{-1}: Hom_{Alg_{R}}(R[x],A) &\to A \\
\psi:R[x] \to A:x \mapsto b &\mapsto b.
\end{split}
\end{equation*}
Además para todo $f:A \to B$, si $g:Hom_{Alg_{R}}(R[x],A) \to Hom_{Alg_{R}}(R[x],B)$ tenemos que claramente $\phi_{B} \circ f= g \circ \phi_{A}.$ \\
Llamaremos también VE al R-esquema que representa, es decir, $VE=Spec(SE)$. Notemos que un elemento $e\in E$ define una transformación natural $\overline{e}:VE \to \underline{A}_{m}$ donde $\overline{e}_{A}(\alpha)=\alpha(e)$ para toda R-álgebra A y todo $\alpha:E \to A\in VE(A)$. Notemos también que $\overline{e}$ es compatible con las acciones de $\underline{A}_{m}:a\overline{e}_{A}(\alpha)=\overline{e}_{A}(a\alpha)$ para todo A y $a\in A$.   

\begin{prop}\label{res:4.13}
Para cualquier R-módulo E, el mapeo natural
$$E \to Mor_{\underline{A}_{m}}(VE,\underline{A}_{m}):e \mapsto \overline{e}$$
es un homomorfismo de R-conjuntos. Por consiguiente, el funtor V de la categoría de R-módulos a la categoría de $\underline{A}_{m}$-conjuntos es completamente fiel. 
\end{prop}

\begin{proof}
Es fácil ver que $e \mapsto \overline{e}$ es un morfismo de R-conjuntos ya que claramente $\overline{re}=r\overline{e}$ para todo $e\in E$ y $r\in R$. Para probar que es biyectivo, comencemos notando que la acción de $\underline{A}_{m}$ en VE está dada por el morfismo único de R-álgebras $\lambda_{E}:SE \to R[T]\otimes SE$ enviando cada $e\in E$ a $T \otimes e$. Se sigue que $\lambda_{E}(f)=\sum_{d} T^{d} \otimes f_{d}$ para todo $f=\sum_{d} f_{d} \in SE$, donde $f_{d}\in S^{d}E$. Por Yoneda, cualquier homomorfismo de funtores $\eta: VE \to \underline{A}_{m}$ está dado por un homomorfismo de anillos $\theta:R[T] \to SE$ el cual es determinado por el elemento $f=\theta(T)$ en SE. Un elemento $\alpha\in VE(A)$ induce un morfismo de R-álgebras $\overline{\alpha}:SE \to A$, y $\eta(\alpha)=\overline{\alpha}(f)$. Entonces $\eta$ es un morfismos de $\underline{A}_{m}$-conjuntos si y solo si el diagrama

$$\xymatrix{
R[T] \ar[d]_{m^{\sharp}} \ar[r]^{\theta} & SE \ar[d]_{\lambda_{E}} \\
R[T] \otimes R[T] \ar[r]_{id\otimes \theta} & R[T] \otimes SE
}
$$

conmuta. Pero $\lambda_{E} \circ \theta(T)= \sum_{d} T^{d}\otimes f_{d}$ y $id\otimes \theta \circ m^{\sharp}(T) = id\otimes \theta(T\otimes T)=T\otimes f$. Entonces el diagrama conmuta si y solo si $f=f_{1}$, es decir, si y solo si $f\in E \subset SE$.  
\end{proof}

\begin{Def}\label{res:4.14}
Sea Q un monoide. Una acción de $\underline{A}_{m}$ en VE es una transformación natural $m_{E}:\underline{A}_{Q}\times VE \to VE$ que satisface las reglas usuales de monoides y que es compatible con la estructura de $\underline{A}_{m}$-conjunto en VE.
\end{Def}

En otras palabras la definición anterior quiere decir que $VE$ debe ser un $(\underline{A}_{Q},\underline{A}_{m})$-biconjunto. El mapeo $m_{E}$ está dado por un homomorfismo $SE \to R[Q]\otimes SE$, y la compatibilidad con la $\underline{A}_{m}$-acción implica que $m_{E}$ es inducido por un morfismo $\theta_{E}:E \to R[Q]\otimes _{R} E$. El morfismo $m_{E}$ va a ser una acción en de $\underline{A}_{m}$ en VE si y solo si $\theta_{E}$ es una coacción de R[Q] en E, es decir, que los siguientes diagramas conmuten

$$\xymatrix{
E \ar[dr]_{id_{E}} \ar[r]^{\theta_{E}} & R[Q] \otimes E \ar[d]^{1_{Q}\otimes id_{E}} \\
& E
}
$$

$$\xymatrix{
E \ar[d]_{\theta_{E}} \ar[r]^{\theta_{E}} & R[Q] \otimes E \ar[d]^{m_{Q}^{\sharp} \otimes id_{E}} \\
R[Q] \otimes E \ar[r]_{id_{R[Q]} \otimes \theta_{E}} & R[Q]\otimes R[Q] \otimes E.
}$$
Para todo $e\in E$ escribimos,
$$\theta_{E}(e)=\sum q\otimes \pi_{q} (e) \in R[Q] \otimes_{R} E.$$
 
Entonces cada $\pi_{q}$ es un endomorfismo R-lineal de E. Los dos diagramas aseguran que $\pi_{p} \pi_{q} = 0$ si $p\neq q$ y que $\sum \pi_{q}= id_{E}$. Entonces la familia $\set{\pi_{q}:q\in Q}$ define una descomposición en suma directa $E = \oplus \set{E_{q}:q\in Q}$, donde $E_{q}$ es la imagen del idempotente $\pi_{q}$. Esta descomposición no es nada más que una Q-graduación del R-módulo E, y la correspondiente acción de $\underline{A}_{Q}(A)$ en $VE(A)$ está dada por 
$$(\alpha,v)(e)= \sum_{q} \alpha(q)v(e_{q}) $$
para todo $\alpha \in \underline{A}_{Q}(A),v \in VE(A)$ y $e= \sum e_{q}$ en E. Todo esto se encuentra en la siguiente proposición.

\begin{prop}\label{res:4.15}
 Sea Q un monoide. Entonces la construcción $E \mapsto VE$ induce un funtor completamente fiel de la categoría de R-módulos Q-graduados a la categoría de biconjuntos $(\underline{A}_{Q},\underline{A}_{m})$. 
\end{prop}

Generalicemos la idea anterior considerando R-módulos graduados por un Q-conjunto S en lugar de Q. 

\begin{Def}\label{res:4.16}
Si S es un Q-conjunto, R[S] es el R-módulo libre con base S y equipado con la estructura única de R[Q]-módulo compatible con la acción de Q en S. 
\end{Def}

El funtor $S \mapsto R[S]$ es adjunto izquierdo del funtor que envía un R[Q]-módulo a si mismo tratándolo como Q-conjunto. Se sigue que si $\set{s_{i}:i\in I}$ es una base para S como Q-conjunto, entonces $\set{e^{s_{i}}:i\in I}$ es una base para R[S] como un Q-módulo. Si S y $S^{\prime}$ son Q-conjuntos, existe un isomorfismo natural $R[S\otimes_{Q} S^{\prime}] \cong R[S] \otimes_{R[Q]} R[S^{\prime}]$. 
\begin{Def}\label{res:4.17}
Sea Q un monoide y S un Q-conjunto. Entonces un R[Q]-módulo S-graduado es un funtor que va de la categoría $\fk{T}S$ a la categoría de R-módulos.  
\end{Def}

La información de un R[Q]-módulo S-graduado E es equivalente a tener una colección de R-módulos $\set{E_{s}:s\in S}$ tal que para todo $q\in Q$ y $s\in S$ existe una función R-lineal $h_{q}:E_{s} \to E_{s+q}$ tal que $h_{q} \circ h_{q^{\prime}}=h_{q+q^{\prime}}$ y $h_{0}=id$. Entonces $\oplus_{s} E_{s}$ es un R[Q]-módulo en el sentido usual. \\ 
Para ver esto supongamos que $E:\fk{T}S \to Mod_{R}$ es un funtor. Entonces para todo $s\in S$ definimos $E_{s}:=E(s)$ un R-módulo. Y definimos 
$$h_{q}:E_{s} \to E_{s+q}:e \mapsto qe.$$
Claramente esta definición respeta la composición y es tal que $h_{0}=Id.$ Por lo que $\oplus E_{s}$ es un R[Q]-módulo.

Consideremos el funtor
$$VS: Alg_{R[Q]} \to Ens$$
el cual envía una R[Q]-álgebra A al conjunto de todos los Q-morfismos $S \to A$, donde A es visto como un Q-conjunto via el carácter $Q \to A$ que viene de su estructura de R[Q]-álgebra. Notemos que VS(A) tiene una estructura natural de A-módulo al definirse puntualmente. El funtor VS es representable por un esquema afín sobre $\underline{A}_{Q}$, el cual vamos a denotar también por VS. De hecho, si A es cualquier R[Q]-álgebra, entonces dar un morfismo de R[Q]-módulos $R[S] \to A$ es equivalente a dar un morfismo de Q-conjuntos $S \to A$. Entonces tenemos un isomorfismo de funtores $VS \cong V(R[S])$ en la categoría $Alg_{R[Q]}$. Por las observaciones iniciales sabemos que 
$$VS\cong V(R[S]) \cong Hom(S(R[S]),-)$$

por lo  que el esquema asociado es $VS := \Spec{S(R[S])}$. Se sigue por la proposición~\ref{res:4.13} que R[S] puede ser identificado con el conjunto de $\underline{A}_{m}$-morfismos $VS \to \underline{A}_{m}$. \\ 
Tenemos un morfismo natural $VS \to \underline{A}_{Q} \to \Spec{R}$. Sea $A_{S}$ el esquema VS visto sobre $\Spec{R}$. Explicitamente, si A es una R-álgebra,
$$A_{S}(A):= \set{(\alpha,\sigma):\alpha \in \underline{A}_{Q}(A), \sigma \in VS(A,\alpha)}.$$

La multiplicación puntual define un mapeo,
$$m_{S}:A_{S} \times A_{S} \to A_{S}$$ 

el cual es bilineal con respecto a la acción de $\underline{A}_{m}$. El correspondiente homomorfismo  de R[Q]-módulos es $R[S] \to R[S] \otimes R[S]:s \mapsto  s\otimes s$. La función constante $S \to R:s \mapsto 1$ induce el morfismo $1_{S}: \Spec{R} \to A_{S}$ y entonces $(m_{S},1_{S})$ da a $A_{S}$ la estructura de R-esquema monoidal. El morfismo natural $A_{S} \to \underline{A}_{Q}$ es un morfismo de R-esquemas monoidales.

\begin{Def}\label{res:4.18}
 Sea E un R[Q]-módulo y S un Q-conjunto. Entonces una estructura equivariante de $(A_{S},\underline{A}_{m})$-biconjuntos en VE es una acción monoidal,
 $$m_{E}: A_{S} \times VE \to VE$$
que es bilineal con respecto a la $\underline{A}_{m}$ estructura en $A_{S}$ y VE, y tal que el diagrama
\end{Def}

$$\xymatrix{
A_{S} \times VE \ar[d] \ar[r]^{m_{E}} & VE \ar[d] \\
\underline{A}_{Q} \times \underline{A}_{Q} \ar[r]_{m_{Q}} & \underline{A}_{Q}
}
$$
conmuta. \\
Si E es S-graduado y A es una Q-álgebra, entonces un elemento de VE(A) puede ser visto como una colección de mapeos R-lineales $n_{s}:E_{s} \to A$ tal que para cada $s\in S$ y $q\in Q$ el diagrama
$$\xymatrix{
E_{s} \ar[d]_{h_{q}} \ar[r]^{n_{s}} & A \ar[d]^{\alpha(q)} \\ 
E_{q+s} \ar[r]_{n_{q+s}} & A
}
$$

conmuta, donde $\alpha:Q \to A$ es el morfismo de la estructura Q-álgebra de A. Sea $(\beta, \sigma)$ un elemento de $A_{S}(A)$. Entonces $\alpha^{\prime}:=\beta \alpha$ es un elemento de $\underline{A}_{Q}(A)$ y podemos definir 
$$n_{s}^{\prime}:E_{s} \to A:= \sigma(s)n_{s}$$
Uno puede ver fácilmente que $\eta^{\prime}$ satisface el diagrama anterior con $\alpha^{\prime}$ en lugar de $\alpha$. Entonces $(\beta,\sigma)\eta:= \eta^{\prime}$ define la acción deseada; la bilinealidad y la compatibilidad con $m_{Q}$ son inmediatas.

\begin{prop}\label{res:4.19}
 Sea S un Q-conjunto. La construcción anterior define un funtor completamente fiel que va de la categoría de R[Q]-módulos S-graduados a la categoría de $(A_{s},\underline{A}_{m})$-biconjuntos equivariantes. 
\end{prop}

\begin{proof}
Sea E un R[Q]-módulo y $\mu:A_{S}\times VE \to VE$ una estructura equivariante de $(A_{S},\underline{A}_{m})$-biconjuntos en VE. Entonces $A_{S}\times VE$ es representado por $S(R[S])\otimes_{R[Q]} SE$. Dado que $\mu$ es compatible con la acción de $\underline{A}_{m}$, se sigue por la proposición~\ref{res:4.13} que $\mu$ es inducido por un homomorfismo de R[Q]-módulos:
$$\mu^{\sharp}:E \to R[S]\otimes_{R} E.$$
La compatibilidad de $\mu$ con $m_{Q}$ implica que este homomorfismo es lineal sobre el homomorfismo $\mu^{\sharp}_{Q}$. El hecho de que $\mu$ es una acción monoidal implica que los siguientes diagramas conmuten:

$$\xymatrix{
E \ar[dr]_{id_{E}} \ar[r]^{\mu^{\sharp}} & R[S] \otimes_{R} E \ar[d]^{1_{S}^{\sharp}\otimes id_{E}} \\
& E
}
$$

$$\xymatrix{
E \ar[d]_{\mu^{\sharp}} \ar[r]^{\mu^{\sharp}} & R[S]\otimes_{R} E \ar[d]_{m_{s}^{\sharp}} \\
R[S]\otimes_{R} E \ar[r]^{id_{R[S] \otimes \mu^{\sharp}}} & R[S]\otimes_{R} R[S] \otimes_{R} E  
}
$$

Si $e\in E$, escribimos $\mu^{\sharp}(e)=\sum s\otimes \pi_{s}(e)$. Entonces cada $\pi_{s}:E \to E$ es un mapeo R-lineal, y estos diagramas dicen que $\sum_{s} \pi_{s}=id_{E}$ y $\pi_{s}\circ \pi_{t}=\delta_{s,t}\pi_{t}$. En otras palabras, $\set{\pi_{s}:s\in S}$ es la familia de proyecciones correspondientes a la descomposición en suma directa $E=\oplus E_{s}$. Si $e\in E_{s}$, entonces $\mu^{\sharp}(e)=s\otimes e$ y, dado que $\mu^{\sharp}$ es lineal sobre $m_{Q}$, 

$$ \mu^{\sharp}(qe)=(q\otimes q)\mu^{\sharp}(e)=(q\otimes q)(s\otimes e) = (q+s)\otimes (qe),$$
asi que $qe\in E_{q+s}$. Entonces el mapeo $h_{q}:E_{s} \to E_{q+s}:e \to eq$ es R-lineal, y $\set{E_{s},h_{q}}$ define una S-graduación en el R[Q]-módulo E. Esto muestra como la S graduación es determinada por la estructura de biconjunto, y se sigue que el funtor es completamente fiel.  
\end{proof}

\begin{obs}\label{res:4.20}
Sea F una cara de un monoide \textit{integral} Q, S un Q-conjunto, y E un R[Q]-módulo S-graduado. Entonces la localización $E_{F}$ de E en F, es naturalmente un $R[Q_{F}]$-módulo $S_{F}$-graduado. Para cualquier $t\in S_{F}$ y $f\in F$, la multiplicación por $e^{f}$ define un isomorfismo de la componente de $E_{F}$ de grado t a la componente de grado f+t. Dado que todos estos isomorfismos conmutan, podemos identificar todas las componentes $E_{F,f+t}$ para todo $f\in F$. Estas componentes puede ser calculadas como sigue. Para cada $s\in S$, hay un sistema directo $\set{e^{f}: E_{s} \to E_{f+s}: f\in F}$, donde F es equipada con el orden monoidal estandar. Entonces uno ve fácilmente que, para cualquier $s\in S$ mapeado a cualquier $f+t\in S_{F}$,
$$\varinjlim \set{e^{f}:E_{s} \to E_{s+f}:f\in F} \cong E_{F,t}.$$
\end{obs}

\subsection{Caras, orbitas y trayectorias}
El funtor $Alg_{R} \to Mon$ que envía una R-álgebra a su monoide multiplicativo se puede levantar naturalmente por el funtor $Alg_{R} \to Moni$ enviando una R-álgebra A al monoide idealizado $((A,\cdot,1),\set{0})$. Si K es un ideal del monoide Q, entonces el R[Q]-módulo R[K] es un ideal de R[Q], y el cociente 
$$R[Q,K]:=R[Q]/R[K]$$
es un R-módulo libre con base $Q\setminus K$. \\

\begin{Def}\label{res:4.21}
Sea Q un monoide y K un ideal de Q. Llamaremos a $R[Q,K]$ el álgebra monoidal del monoide idealizado (Q,K). 
\end{Def}

Sea $\underline{A}_{Q,K}$ el funtor que toma una R-álgebra A y la envía al conjunto de morfismos idealizados $(Q,K) \to (A,0)$. Veamos que este funtor esta representando por $R[Q,K]$, es decir, que 
$$\underline{A}_{Q,K}(A)=Hom_{Moni}((Q,K),(A,0)) \cong Hom_{Alg_{R}}(R[Q,K],A).$$
Esto se cumple al considerar la transformación natural,
$$\phi:Hom_{Moni}((Q,K),-) \to Hom_{Alg_{R}}(R[Q,K],-)$$
tal que para todo $A\in Alg_{R}$,
\begin{equation*}
\begin{split}
\phi_{A}:Hom_{Moni}((Q,K),(A,0)) &\to Hom_{Alg_{R}}(R[Q,K],A) \\
\phi:(Q,K) \to (A,0) &\mapsto  \varphi:R[Q,K] \to A:\overline{\sum a_{i}e^{q_{i}}} \mapsto \sum a_{i}\phi(q_{i})
\end{split}
\end{equation*}
con inversa,
\begin{equation*}
\begin{split}
\phi_{A}^{-1}:Hom_{Alg_{R}}(R[Q,K],A) &\to Hom_{Moni}((Q,K),(A,0)) \\
\varphi:R[Q,K] \to A:\overline{\sum r_{i}e^{q_{i}}} \mapsto \sum r_{i}a(q_{i}) &\mapsto \phi:(Q,K) \to (A,0): q \mapsto a(q).   
\end{split}
\end{equation*}

Además si $f:A \to B$ es un morfismo de R-álgebras, y $$g:Hom((Q,K),(A,0) \to Hom((Q,K),(B,0)):\varphi:(Q,K) \to (A,0) \mapsto f \circ \varphi$$ 

$$h:Hom(R[Q,K],A) \to Hom(R[Q,K],B): \lambda:R[Q,K] \to A \mapsto f\circ \lambda$$. 

Entonces, claramente $\phi_{B} \circ g=h\circ \phi_{A}$, por lo que $\phi$ si define un isomorfismo natural.

De igual forma denotaremos por $\underline{A}_{Q,K}$ al R-esquema que representa, es decir, $\underline{A}_{Q,K}=\Spec{R[Q,K]}$. Entonces tenemos que $\underline{A}_{Q,K}$ es un R-subesquema cerrado del esquema $\underline{A}_{Q}$, y la proposición~\ref{res:4.15} muestra que es invariante bajo la acción de $\underline{A}_{Q}$ en si mismo. Dicho de otro modo, $\underline{A}_{Q,K}$ es un ideal del monoide $\underline{A}_{Q}$. \\
En particular, sea $\fk{p}$ un ideal primo de Q y $F=Q\setminus \fk{p}$ su correspondiente cara. La inclusión $F \to Q$ define un homomorfismo de R-álgebras $R[F] \to R[Q]$, y esto a su vez define un homomorfismo de esquemas,
$$r_{F}:\underline{A}_{Q} \to \underline{A}_{F}.$$
La composición del mapeo $R[F] \to R[Q]$ con el homomorfismo 
$$i_{\fk{p}}^{\sharp}:R[Q] \to R[Q,\fk{p}]:f \mapsto \overline{f}$$
induce un isomorfismo de R-álgebras $R[F]\to R[Q,\fk{p}]$, ya que induce una biyección en los elementos de la base. Este isomorfismo de anillos define un isomorfismo de esquemas $\underline{A}_{Q,\fk{p}} \to \underline{A}_{F}$, y entonces tenemos 
$$i_{F}:\underline{A}_{F} \to \underline{A}_{Q}$$
el cual es la composición del inverso del isomorfismo anterior con la inmersión cerrada $i_{\fk{p}}$. Entonces,

\[ 
i_{F}^{\sharp}(e^{q})=
\begin{cases} 
      e^{q} & q\in F \\
      0 & otro~caso 
   \end{cases}
\]

Veamos por ejemplo que si Q es \textit{sharp} entonces $\underline{A}_{Q,Q^{+}}\cong \Spec{R}$. Dado que Q es \textit{sharp}, tenemos que $Q^{+}=Q\setminus \set{0}$. Consideremos el homomorfismo,
\begin{equation*}
\begin{split}
\phi:R[Q]/R[Q^{+}] &\to R \\
\overline{a_{0}e^{0}} &\mapsto a_{0}.
\end{split}
\end{equation*}
el cual claramente es biyectivo. Por tanto, al ser $R[Q]/R[Q^{+}] \cong R$, tendremos que $\Spec{R}\cong \Spec{R[Q]/R[Q^{+}]}$.

\begin{prop}\label{res:4.22}
Sea F una cara de un monoide \textit{integral} Q, $i_{F},r_{F}$ los morfismos definidos anteriormente, y sea $i_{Q/F}$ la inmersión cerrada inducida por la proyección $Q \to Q/F$. Sea S el espectro del anillo base R. \\
1. Estos homomorfismos se ajustan a un diagrama conmutativo con cuadrados cartesianos: 

$$\xymatrix{
S \ar[d]_{1_{F}} \ar[r]^{v_{Q/F}} & \underline{A}_{Q/F} \ar[d]_{i_{Q/F}} \ar[r]^{\pi} & S \ar[d]_{1_{F}} \\
\underline{A}_{F} \ar[r]_{i_{F}} & \underline{A}_{Q} \ar[r]_{r_{F}} & \underline{A}_{F}
}
$$
En este diagrama $1_{F}$ es el punto S-valuado correspondiente a la sección identidad del esquema monoidal $\underline{A}_{F}$ y $v_{Q/F}$ es el vértice del esquema monoidal $\underline{A}_{Q/F}.$
La composición $r_{F} \circ i_{F}$ es $id_{\underline{A}_{F}}$. \\
2. El mapeo $r_{F}$ es un morfismo de esquemas monoidales, y el morfismo $i_{F}$ es compatible con las acciones del esquema monoidal $\underline{A}_{Q}$ en sí mismo y en el esquema ideal $\underline{A}_{F}$. \\ 
3. Si Q es fino, entonces $i_{F}$ es una strong deformation retrac. Esto es, existe un morfismo 
$$f:\underline{A}_{Q} \times \underline{A}_{m} \to \underline{A}_{Q}$$
tal que $f \circ j_{0}=i_{F} \circ r_{F}$, $f\circ j_{1}=id$, y $f \circ (i_{F}\times id)= i_{F}\circ pr_{1}$, donde $j_{0},j_{1}:\underline{A}_{Q} \to \underline{A}_{Q}\times \underline{A}_{m}$ están dados respectivamente por las secciones 0 y 1 de $\underline{A}_{m}$  y $pr_{1}: \underline{A}_{F}\times \underline{A}_{m} \to \underline{A}_{F}$ es la proyección. 
\end{prop}

\begin{proof}
La inmersión cerrada $i_{F}$ preserva la ley de conmposición para el esquema monoidal $\underline{A}_{F}$ y $\underline{A}_{Q}$ pero no la sección identidad de la estructura de esquema monoidal, así que $\underline{A}_{F}$ no puede ser considerado como un submonoide de $\underline{A}_{Q}$. Por otro lado, la inclusión $F \to Q$ define un morfismo $R[F] \to R[Q]$ y entonces un mapeo $r_{F}:\underline{A}_{Q} \to \underline{A}_{F}$. Dado que $r_{F}$ es inducido por un homomorfismo de monoides, es un homomorfismo de esquemas monoidales. Se sigue de la definición que $r_{F} \circ i_{f}=id_{\underline{A}_{F}}$. Entonces $r_{F},i_{F}$ son morfismos de $\underline{A}_{Q}$-conjuntos, y $r_{F}(a)=r_{F}(a\cdot 1)=ar_{F}(1_{A})$ para todo $a\in \underline{A}_{Q}(A)$. \\
Uno puede ver fácilmente los dos cuadrados conmutan. Todos los morfismos en el rectángulo exterior son mapeos identidad, así que es cartesiano, así que automáticamente el cuadrado izquierdo es cartesiano si el de la derecha lo es. El ultimo enunciado asegura que el ideal de la inmersión cerrada $i_{Q/F}$ es el ideal I generado por el conjunto de todos los $e^{f}-1$ tal que $f\in F$. De hecho es evidente que $i_{Q/F}^{\sharp}$ aniquila todos estos elementos y se factoriza a través del morfismo $R[Q]/I \to R[Q/F]$. Por otro lado, el mapeo $Q \to R[Q]/I$ envía F a 1, y entonces se factoriza Q/F la propiedad de mapeo universal. Esto da el mapeo inverso $R[Q/F] \to R[Q]/I$. \\
Si Q es fino, entonces por la proposición~\ref{res:3.26}, existe un morfismo $h:Q \to \bb{N}$ tal que $h^{-1}(0)=F$. Este h define un morfismo $t:\underline{A}_{m} \to \underline{A}_{Q}$;  en puntos A-valuados $t(a)=a^{h}$, donde $a^{h}$ es el  morfismo $Q \to A:q \mapsto a^{h(q)}$. Sea
$$f:\underline{A}_{Q}\times \underline{A}_{m} \to \underline{A}_{Q}$$ la composición de $id_{\underline{A}_{Q}}\times t$ con el mapeo multiplicación $m_{Q}$ de la estructura monoidal de $\underline{A}_{Q}$. En puntos A-valuados, f envía $(\alpha,a)$ a $\alpha a^{h}$. Sea $i_{0},i_{1}$ las secciones correspondientes al punto 0 y 1, y $j_{0},j_{1}$ los mapeos correspondientes $\underline{A}_{Q} \to \underline{A}_{Q} \times \underline{A}_{m}$. Podemos ver que $f\circ j_{0}=i_{F}\circ r_{F}$ y que $f\circ j_{1}=id$ en los puntos A-valuados. La segunda ecuación es obvia, y para la primera solo tenemos que observar que $f(\alpha,a)=\alpha 0^{h}$. Finalmente, si $\alpha$ pertenece a la imagen de $i_{F}$ entonces para todo $a\in A$, $f(\alpha,a)(q)=\alpha(q)a^{h(q)}=\alpha(q)$, ya que $\alpha(q)=0$ cuando $h(q)\neq 0$. Esto prueba que $i_{F}$ es deformación retracta fuerte. 
\end{proof}

\begin{cor}\label{res:4.23}
Si Q es un monoide fino, entonces $\underline{A}_{Q}(\bb{C})$, equipada con la topología compleja, es conexo si y solo si $Q^{\ast}$ es libre de torsión. Si Q es \textit{sharp} entonces $\underline{A}_{Q}(\bb{C})$ es contractible. 
\end{cor}

\begin{proof}
La proposición anterior implica que $\underline{A}_{Q}(\bb{C})$ y $\underline{A}_{Q^{\ast}}(\bb{C})$ tienen el mismo tipo de homotopía. Ya que el subgrupo de torsión de $Q^{\ast}$ se identifica con el grupo de componentes conexas del grupo algebraico $\underline{A}_{Q^{\ast}}$, se sigue que $\underline{A}_{Q}(\bb{C})$ y $\underline{A}_{Q^{\ast}}(\bb{C})$ son conexas si y solo si $Q^{\ast}$ es libre de torsión. Si Q es \textit{sharp}, entonces el esquema $\underline{A}_{Q^{\ast}}$ se reduce a un solo punto, así que $\underline{A}_{Q}(\bb{C})$ es contractible. 
\end{proof}

Cuando k es un campo y Q es un monoide \textit{integral}, el monoide $\underline{A}_{Q}(k)$ admite una estratificación explícita indexada por las caras de Q. Si $x\in \underline{A}_{Q}(k)=Hom(Q,k)$, sea $F(x)=x^{-1}(k^{\ast})$, una cara de Q. Si x y z son puntos de $\underline{A}_{Q}(k)$, entonces
$$F(xz)=F(x) \cap F(z).$$ 
Notemos que x es cero fuera de F(x) e induce un mapeo $F(x)^{gp} \to k^{\ast}$ el cual de hecho determina a x. Entonces podemos ver un punto de $\underline{A}_{Q}(k)$ como un par $(F,x^{\prime})$ donde F es una cara de Q y $x^{\prime}$ es un homomorfismo $F^{gp} \to k^{\ast}$.

\begin{prop}\label{res:4.24}
Sea Q un monoide fino, k un campo y F una cara de Q. Entonces el conjunto de todos los $y\in \underline{A}_{Q}(k)$ tal que $F(y)=F$ es $$\underline{A}_{F}^{\ast}(k):= \underline{A}_{F^{gp}}(k) \subset \underline{A}_{F}(k),$$
un subconjunto de zariski abierto de $\underline{A}_{F}(k)$. Si x y y son dos puntos de $\underline{A}_{Q}(k)$, entonces los siguientes son equivalentes: \\
1. $F(y) \subset F(x)$ \\
2. $y\in \underline{A}_{F(x)}(k)$ \\
3. Existe un $z\in \underline{A}_{Q}(k)$ tal que $y=zx$. 
\end{prop}

\begin{proof}
Identificamos un punto $y\in \underline{A}_{Q}(k)$ con el correspondiente carácter $Q \to k$. Entonces $F(y) \subset F(x)$ si y solo si $y(Q\setminus F(x))=0$, es decir, si y solo si $y\in \underline{A}_{F(x)}(k)$; entonces se tiene la equivalencia de (1) y (2). Por (3) del teorema~\ref{res:3.20}, F es fino, así que $\underline{A}_{F}^{\ast}$ es Zariski denso y abierto en $\underline{A}_{F}$, y la inclusión $\underline{A}_{F}^{\ast}(k) \subset \underline{A}_{F}(k) \subset \underline{A}_{Q}(k)$ identifica $\underline{A}_{F}^{\ast}(k)$ con el conjunto de todos los y tal que $F(y)=F$.  Si $F(y) \subset F(x)$, definimos $z: Q \to k$ como $z(q)=0$ si $q\in Q\setminus F(x)$ y $z(q)=y(q)/x(q)$ si $q\in F(x)$. Entonces usando el hecho de que F es una cara de Q, uno puede ver que de hecho $z\in \underline{A}_{Q}(k)$ y $y=zx$. Entonces (1) implica (3) y el converso es claro. Si $F(x)=F(y)=F$, entonces $x/y$ define un morfismo $F^{gp} \to k$. Si k es algebraicamente cerrado, $k^{\ast}$ es divisible, y si $Q^{gp}/F^{gp}$ es libre de torsión, la sucesión $F^{gp} \to Q^{gp} \to Q^{gp}/F^{gp}$ se divide. En cualquier caso, existe una extensión z de $x/y$ a $Q^{gp}$, la cual define un punto de $\underline{A}_{Q}^{\ast}$ tal que $y=xz$.   
\end{proof}

Esta proposición describe la acción de $\underline{A}_{Q}^{\ast}(k)$ en $\underline{A}_{Q}(k)$ donde k es un campo.  

\begin{prop}\label{res:4.25}
Sea $\theta: P \to Q$ un homomorfismo de monoides \textit{integral}es, escribimos $Q^{gp}/P^{gp}$ por $Cok(\theta^{gp})$, y $\underline{A}_{Q/P}^{\ast}:= \underline{A}_{Q^{gp}/P^{gp}}$. Entonces $\underline{A}_{Q/P}^{\ast}$ actúa naturalmente en $\underline{A}_{Q}$, visto como un objeto sobre $\underline{A}_{P}$. Para cada cara F de Q, los subfuntores $\underline{A}_{Q_{F}} \subset \underline{A}_{Q}$ y $\underline{A}_{F}\cong \underline{A}_{Q,p_{F}}$ son estables bajo esta acción. Más aún, el subgrupo $\underline{A}_{Q/(P+F)}$ de $\underline{A}_{Q/P}^{\ast}$ actua trivialmente en el subfuntor cerrado $\underline{A}_{F} \cong \underline{A}_{Q,p_{F}}$ de $\underline{A}_{Q}$. 
\end{prop}

\begin{proof}
Sea A un anillo conmutativo y supongamos que $\gamma \in \underline{A}_{Q/P}^{\ast}(A)$ y $\alpha \in \underline{A}_{Q}(A)$. Entonces $\gamma(\theta(p))=1$ para todo $p\in P$ y entonces $(\gamma \alpha )(\theta(p))= \alpha(\theta(p))$. En otras palabras, $\gamma$ actúa en $\underline{A}_{Q}$ como un objeto sobre $\underline{A}_{P}$. Más aún, para todo $q\in Q$, $\gamma(q) \in A^{\ast}$. En particular, $\gamma(q) \alpha(q)=0$ si y solo si $\alpha(q)=0$, y aplicando esto para todo $q\in \fk{p}_{F}$, vemos que $\gamma \alpha \in \underline{A}_{Q,\fk{p}_{F}}(A)$ si y solo si $\alpha \in \underline{A}_{Q,\fk{p}_{F}}(A)$. Similarmente, $\gamma(q) \alpha(q) \in A^{\ast}$ si y solo si  $\alpha(q) \in A^{\ast}$, y aplicando esto para todo $q\in F$ vemos que $\alpha \gamma \in \underline{A}_{Q_{F}}(A)$ si y solo si $\alpha \in \underline{A}_{Q_{F}}(A)$. Ahora supongamos que $\gamma \in \underline{A}_{Q/P+F}(A)$ y $\alpha \in \underline{A}_{F}(A) \cong \underline{A}_{Q,\fk{p}_{F}}(A)$. Finalmente notemos que $\gamma(q) \alpha(q)=\alpha(q)$ para todo $q\in Q$ ya que si $q\in \fk{p}_{F}$ entonces $\alpha(q)=0$, y $q\in F$ entonces $\gamma(q)=1$.     
\end{proof}

\begin{obs}\label{res:4.26}
Observa que el homomorfismo $Q^{gp} \to Q^{gp}/P^{gp}$ induce un homomorfismo de esquemas $\underline{A}_{Q/P}^{\ast} \to \underline{A}_{Q}^{\ast}$ identificando $\underline{A}_{Q/P}^{\ast}$ como el kernel del morfismo $\underline{A}_{Q}^{\ast} \to \underline{A}_{P}^{\ast}$. La acción de $\underline{A}_{Q/P}^{\ast}$ en $\underline{A}_{Q}$ viene dada por,
$$m_{Q/P}: \underline{A}_{Q} \times \underline{A}_{Q/P}^{\ast} \to \underline{A}_{Q}.$$
Este homomorfismo corresponde al homomorfismo de monoides:
$$\theta: Q \to Q \otimes Q^{gp}/P^{gp}:q \mapsto (q\otimes [q]).$$
Sea $(\alpha,\gamma)\in \underline{A}_{Q}(A) \times \underline{A}_{Q/P}^{\ast}(A)$ y $q\in Q$,
\begin{equation*}
\begin{split}
(\alpha,\gamma)(m^{ \sharp}_{Q/P}(q)) &= m_{Q/P}(\alpha,\gamma)(q) \\
&= (\alpha \gamma)(q) \\
&= \alpha(q) \gamma([q]) \\
&= (\alpha,\gamma)(q,[q]) \\
&= (\alpha,\gamma)(\theta(q))
\end{split}
\end{equation*}
La acción de $\underline{A}_{Q/P}^{\ast}$ en $\underline{A}_{Q}$ induce un mapeo,
$$h_{Q/P}:\underline{A}_{Q} \times \underline{A}_{Q/P}^{\ast} \to \underline{A}_{Q} \times_{\underline{A}_{P}} \underline{A}_{Q}: (\alpha, \gamma) \mapsto (\alpha,\alpha \gamma)$$
Este mapeo corresponde al homomorfismo de monoides
$$\phi: Q \oplus_{P} Q \to Q\oplus Q^{gp}/P^{gp}:[q_{1},q_{2}] \mapsto (q_{1},+q_{2},[q_{2}]).$$
Notemos que los morfismos inducidos 
$$Q^{gp} \oplus_{P^{gp}} Q^{gp} \to Q^{gp} \oplus Q^{gp}/P^{gp},~~~~~ \underline{A}_{Q}^{\ast} \times \underline{A}_{Q/P}^{\ast} \to \underline{A}_{Q}^{\ast} \times_{\underline{A}_{P}^{\ast}} \underline{A}_{Q}^{\ast}$$
son isomorfismos. El inverso del morfismo $\phi^{gp}$ es
$$Q^{gp} \oplus Q^{gp} /P^{gp} \to Q^{gp} \oplus_{P^{gp}} Q^{gp}: (x,[y]) \mapsto [(x-y,y)].$$
\end{obs}

\subsection{Geometría local de las variedades tóricas afines}
\begin{prop}\label{res:4.27}
Sea P un monoide \textit{integral} y R un anillo. \\
1. Si $P^{gp}$ es libre de torsión y R es un dominio entero, entonces R[P] es un dominio entero. \\
2. Si en adición P es finitamente generado y R es normal, $R[P^{sat}]$ es la normalización de R[P]. Entonces en este caso R[P] es normal si y solo si P es saturado. \\
3. Si P es fino, el morfismo $\pi: \Spec{R[P]} \to \Spec{R}$ es  fielmente plano y de representación finita. Más aún, la dimensión de Krull de las fibras de $\pi$ es el rango de $P^{gp}$. Si en adición el rango del subgupo de torsión de $P^{gp}$ es invertible en R, entonces las fibras son geometricamente reducidas.  
\end{prop}

\begin{proof}
Los detalles de la demostración se pueden encontrar en ~\cite{Ogus:2018aa}. 
\end{proof}

Para ver que la hipótesis de que $P^{gp}$ es libre de torsión sí es necesaria, consideremos el submonoide P de $\bb{Z}\oplus \bb{Z}/ 2\bb{Z}$ generado por $p=(1,0)$ y $q=(1,1)$. Este es \textit{sharp} y fino, pero $R[P] \cong R[x,y]/(x^{2}-y^{2})$, el cual no es un dominio entero ya que $(x^{2}-y^{2})$ no es un ideal primo.

\begin{prop}\label{res:4.28}
Sea Q un monoide fino, k un campo y $\underline{A}_{Q}:=\Spec{k[Q]}$. \\
Si $\fk{p}$ es un ideal primo de Q y $F=Q \setminus \fk{p}$ es su cara correspondiente, entonces la altura de $\fk{p}$ en Q es la codimensión del subesquema cerrado $\underline{A}_{F}$ de $\underline{A}_{Q}$.
\end{prop}

\begin{proof}
Notemos que en estos enunciados estamos usando la inmersión cerrada $i_{F}$ para identificar $\underline{A}_{F}$ y $\underline{A}_{Q,\fk{p}}$. La dimensión de $\underline{A}_{F}$ es el rango de $F^{gp}$ y la dimensión de $\underline{A}_{Q}$ es el rango de $Q^{gp}$, así que la codimensión del subconjunto cerrado localmente $\underline{A}_{F}$ de $\underline{A}_{Q}$ es el rango de $Q^{gp}/F^{gp}$. Esta es la altura de $\fk{p}$ como ideal primo de Q por la proposición~\ref{res:2.56}. 
\end{proof}

\begin{Theorem}\label{res:4.29}
Si R es un anillo Cohen-Macaulay noetheriano, y P es monoide saturado fino, entonces el anillo monoidal $R[P]$ es también Cohen-Macaulay. 
\end{Theorem}

\begin{proof}
Los detalles de la demostración se pueden encontrar en ~\cite{hochster1972rings}.
\end{proof}

\subsection{Ideales en álgebras monoidales}
Sea P un monoide \textit{integral} y R un anillo no cero. Recordemos que si K es un ideal de P entonces $R[K]$ es un ideal de $R[P]$. Además, si K es la intersección de una familia de ideales $K_{\lambda}$, entonces 
$$R[K] = \cap R[K_{i}].$$

\begin{Def}\label{res:4.30}
Un anillo R es llamado reducido si no tiene elementos no cero idempotentes. 
\end{Def}

\begin{Def}\label{res:4.31}
Un anillo conmutativo con unidad es normal si para cada ideal primo $\fk{p}$ la localización $R_{\fk{p}}$ es un dominio entero y es cerrado en su campo de fracciones. 
\end{Def}

\begin{prop}\label{res:4.32}
Sea Q un monoide \textit{integral}, R un anillo reducido  y K un ideal de Q. Supongamos que el orden del subgrupo de torsión de $Q^{gp}$ es invertible en R. Entonces el anillo $R[Q,K]$ es reducido si y solo si K es un ideal radical de Q. 
\end{prop}

\begin{proof}
Recordemos que $Q\setminus K \to R[Q,K]:q \mapsto e^{q}$ es una base para $R[Q,K]$. Entonces es claro que si $R[Q,K]$ es reducida, entonces K debe ser un ideal radical, ya que de otro modo existiría $q\in Q\setminus K$ y $n\in \bb{N}$ tal que $nq\in K$, y entonces $e^{q}$ sería un nilpotente no cero de R[Q,K]. Para el converso, primero notemos que si K es primo, su complemento es una cara F de Q, y $R[F] \cong R[Q,K]$. Más aún, $R[F] \subset R[F^{gp}] \cong R[\bb{Z}^{r} \oplus T]$ donde T es un grupo finito cuyo orden es invertible en R. Se sigue que $R[F^{gp}]$ y R[F] son reducidos. Si K es un ideal radical de P, entonces es la intersección de todos los ideales primos $\fk{p}_{\lambda}$ que lo contienen, y se sigue que R[Q,K] es la intersección de todos los ideales $R[Q,\fk{p}_{\lambda}]$. Ya que cada uno de estos es reducido, también lo es R[Q,K].   
\end{proof}

\begin{prop}\label{res:4.33}
Sea R un dominio entero, P un monoide tórico y $p\in P$. Entonces las componentes irreducibles de $\Spec{R[P,(p)]}$ son precisamente los conjuntos cerrados definidos por los ideales primos de P de altura uno conteniendo a p. Si P es saturado y R es normal, cada una de estas componentes irreducibles es normal. 
\end{prop}

\begin{proof}
Sea $\fk{p}$ un ideal primo de P y G su cara correspondiente. Ya que $R[G] \cong R[P,\fk{p}]$, tenemos que es un dominio entero y es normal si R es normal y G es saturado por la proposición~\ref{res:4.27}. Por un corolario, el ideal $\sqrt{(p)}$ de P puede ser escrito como la intersección de primos de altura uno conteniendo a p, $\fk{p}_{1} \cap ... \cap \fk{p}_{n}$. Concluimos que $\Spec{R[P,(p)]}= \Spec{R[P,\fk{p}_{1}]} \cup ... \cup \Spec{R[P,\fk{p}_{n}]}$. Como hemos visto, cada una de estas partes es irreducible y entonces estas son de hecho las componentes irreducibles de $R[P,(p)]$. Si P es saturado, también lo es cada una de las caras $G_{i}=P \setminus \fk{p}_{i}$, y entonces si R es normal, también lo es cada $R[P,\fk{p}_{i}] \cong R[G_{i}]$.  
\end{proof}

Esta proposición muestra que los divisores de Weil con soporte en un subconjunto cerrado de $\Spec{R[P]}$ definido por un elemento de P provienen de ideales en P.

\begin{Def}\label{res:4.34}
Sea P un monoide \textit{integral} y R un anillo. Si $f=\sum_{p} a_{p}(f)e^{p} $ es un elemento de $R[P]$ y S es un subconjunto de $R[P]$, entonces: \\
1. $\sigma(f) := \set{p\in P: a_{p}(f) \neq 0}$, $\sigma(S):= \cup \set{\sigma(f):f\in S};$ \\
2. K(f) es el ideal de P generado por $\sigma(f)$, y $K(S):=\cup \set{K(f):f\in S}.$
\end{Def}

Notemos que $K(S)$ es un ideal de P, ya que la unión de ideales es un ideal. Más aún, para todo $p\in P$ y $f\in R[P]$,
$$\sigma(e^{p}f)=p+\sigma(f).$$
Sea $p+p^{\prime}\in p+\sigma(f)$. Entonces,
\begin{equation*}
\begin{split}
a_{p+p^{\prime}}(e^{p}f) &= a_{p+p^{\prime}}(\sum a_{r}(f)e^{p+r}) \\ &= a_{p^{\prime}}(f) \\
&\neq 0
\end{split}
\end{equation*}
por lo que $p+p^{\prime}\in \sigma(e^{p}f)$. Sea $q\in \sigma(e^{p}f)$. Notemos que,
\begin{equation*}
\begin{split}
0 &\neq a_{q}(e^{p}f) \\
&= a_{q}(\sum a_{r}(f)e^{p+r}) \\
&= a_{t}(f),~p+t=q.
\end{split}
\end{equation*}
Esto además implica que $t\in \sigma(f)$, por lo que $q\in p+\sigma(f)$. Se sigue que si I es un ideal de R[P], entonces $\sigma(I)=K(I)$, ya que si $k\in K(I)$, existe $f\in I$ y $p\in P$ tal que $k\in p+ \sigma(f)=\sigma(e^{p}f)\subset \sigma(I)$. Notemos que cualquier $f\in R[P]$ está contenido en $R[K(f)]$, que cualquier ideal I de $R[P]$ está contenido en $R[K(I)]$, y que de hecho K(I) es el ideal más pequeño de P tal que $I\subset R[K(I)]$.

\begin{prop}\label{res:4.35}
Supongamos que f y g son elementos de $R[P]$. \\
1. $\sigma(f+g) \subset \sigma(f) \cup \sigma(g)$, entonces $K(f+g)\subset K(f) \cup K(g)$. \\
2. $\sigma(fg)\subset \sigma(f)+\sigma(g)$, entonces $K(fg) \subset K(f)+K(g) \subset K(f) \cap K(g)$. \\
3. $K(f)=K((f))$, donde (f) es el ideal de $R[P]$ generado por f. \\
4. Si I y J son ideales de R[P], $K(IJ)\subset K(I)+K(J)$. 
\end{prop}
 
\begin{proof}
Claramente se cumple:
$$a_{p}(f+g)=a_{p}(f)+a_{p}(g),$$
$$a_{p}(fg)=\sum_{q+q^{\prime}=p} a_{q}(f)a_{q^{\prime}}(g).$$
Entonces,
\begin{equation*}
\begin{split}
\sigma(f+g) &= \set{p\in P:a_{p}(f+g)\neq 0} \\
&= \set{p\in P:a_{p}(f)+a_{p}(g) \neq 0} \\
&\subset \set{p\in P:a_{p}(f)\neq 0~o~a_{p}(g)\neq 0} \\
&= \sigma(f)\cup \sigma(g).
\end{split}
\end{equation*}
De igual forma.

\begin{equation*}
\begin{split}
\sigma(fg) &= \set{p\in P:a_{p}(fg)\neq 0} \\
&= \set{p\in P:\sum_{q+q^{\prime}=p}a_{q}(f)a_{q^{\prime}}(g)\neq 0} \\
&\subset \sigma(f)+\sigma(g).
\end{split}
\end{equation*}

Notemos que (4) es una consecuencia inmediata de (1) y (2) ya que,
\begin{equation*}
\begin{split}
K(IJ) &= \cup_{f\in I,g\in J} K(fg) \\
&\subset \cup_{f\in I,g\in J} K(f)+K(g) \\
&\subset \cup_{f\in I} K(f) + \cup_{g\in J} K(g) \\
&= K(I)+K(J).
\end{split}
\end{equation*}

Se sigue de la definición que $\sigma(f) \subset K((f))$, y entonces $K(f) \subset K((f))$. Por otro lado, para cualquier $h\in (f)$, se sigue de (2) que $\sigma(h) \subset K(f)$ y entonces que $K(h) \subset K(f)$. Al ser h arbitrario concluimos que $K((f)) \subset K(f)$.  
\end{proof}

\begin{prop}\label{res:4.36}
Sea P un monoide \textit{integral} y R un anillo. \\
1. Si $f\in R[P]$, K(f) es el ideal unidad de P si y solo si $f\notin R[P^{+}]$. \\
2. Si $f\in R[P]$, entonces K(f) es generado principalmente por un elemento p de P si y solo si $f=e^{p}\overline{f}$, donde $\overline{f}$ es algún elemento de $R[P]\setminus R[P^{+}]$. \\ 
3. Supongamos que R es un dominio entero y $P^{\ast}$ es libre de torsión. Entonces si f y g son elementos de R[P] tal que K(f) y K(g) son principales, lo mismo es cierto para fg, y $K(fg)=K(f)+K(g)$. 
\end{prop}

\begin{proof}
El enunciado (1) es una tautología. Si $K:=K(f)$ es generado por p, entonces $k-p\in P$ para todo $k\in K(f)$. Entonces, $f= \sum_{k\in K} a_{k}e^{k}=e^{p} \sum_{k} a_{k}e^{k-p}$, así que $f=e^{p}\overline{f}$ con $\overline{f}=\sum_{k} a_{k}e^{k-p}$. Entonces,
$$(p)=K(f) \subset K(e^{p}) + K(\overline{f}) =(p)+K(\overline{f})$$
y se sigue que $K(\overline{f})=P$. Por el contrario, si $f=e^{p}\overline{f}$ con $K(\overline{f})=P$ entonces ciertamente $K(f) \subset (p)$. Pero si $\overline{f}=\sum \overline{a_{q}}e^{q}$ existe un $q\in P^{\ast}$ tal que $\overline{a_{q}}\neq 0$, y entonces $p+q\in K(f)$, así $p\in K(f)$. Entonces $K(f)=(p)$, y esto prueba (2). Si K(f) es generado principalmente por p y K(g) es generado principalmente por q, entonces $f=e^{p} \overline{f}$ y $g=e^{q} \overline{g}$ donde $\overline{f},\overline{g} \in R[P]\setminus R[P^{+}]$. El cociente de R[P] por $R[P^{+}]$ es isomorfo a $R[P^{\ast}]$. Si $P^{\ast}$ es libre de torsión y R es un dominio entero, entonces $R[P^{\ast}]$ también es un dominio entero por la proposición~\ref{res:4.27}. Entonces $R[P^{+}]$ es un ideal primo y $\overline{f} \overline{g} \notin R[P^{\ast}]$. Ya que $fg=e^{p+q}\overline{f} \overline{g}$, se sigue que K(fg) es generado principalmente por $p+q$. 
\end{proof}

El enunciado (3) muestra que si R es un dominio y $P^{\ast}$ es libre de torsión, entonces el conjunto de todos los f tal que K(f) es principal es un submonoide de R[P]. Esto no es cierto en general. \\

Recordemos del corolario~\ref{res:3.39} que asociado a cada ideal primo $\fk{p}$ de altura uno de un monoide fino P existe un homomorfismo $v_{\fk{p}}:P \to \bb{N}$; este homomorfismo es suprayectivo si P es saturado, como vamos a suponer. La imagen bajo $v_{\fk{p}}$ de un ideal no vacío K de P es entonces un ideal $K_{\fk{p}}$ de $\bb{N}$, generado principalmente por 
$$v_{\fk{p}}(K):=inf \set{v_{\fk{p}}(k):k\in K}.$$
Para cada $f\in R[P]$, sea 
$$v_{\fk{p}}(f):=v_{\fk{p}}(K(f)).$$
Es decir, $v_{\fk{p}}(f)$ es el mínimo del conjunto de todos los $v_{\fk{p}}$ tal que $p\in \sigma(f)$.

\begin{prop}\label{res:4.37}
Sea P un monoide tórico y R un dominio entero. \\
1. Si K es un ideal no vacío de P y $p\in K$ es un elemento tal que $v_{\fk{p}}(p)=v_{\fk{p}}(K)$ para todo ideal primo $\fk{p}$ de altura uno, entonces $K=(p)$. \\ 
2. Si f y g son elementos de $R[P]$, entonces $v_{\fk{p}}(fg)=v_{\fk{p}}(f)+v_{\fk{p}}(g)$, para todo ideal primo de altura uno. Más aún, K(fg) es principal si y solo si K(f) y K(g) lo son. 
\end{prop}

\begin{proof}
Supongamos que se cumple (1)y sea $k\in K$. Entonces $v_{\fk{p}}(k-p) \geq 0$ para todo ideal primo de altura uno de P. Por el corolario~\ref{res:3.40}, $k-p\in P$, por lo que se sigue que K es generado principalmente por p. Esto prueba (1). Para (2), sea G la careta de P complementaria a $\fk{p}$. El homomorfismo $P \to P_{G}$ es inyectivo, y entonces si $\overline{f}\in R[P_{G}]$ es la imagen de f, $K(\overline{f})$ es principal, generado por cualquier $p\in \sigma(f)$ tal que $v_{\fk{p}}(f)=v_{\fk{p}}(p)$. Ya que $P_{G}^{\ast}=G^{gp}$ es libre de torsión, (3) de la proposición~\ref{res:4.36} implica que para todo $f,g\in R[P]$, $K(\overline{f}\overline{g})=K(\overline{f})+K(\overline{g})$ y entonces $v_{\fk{p}}(fg)=v_{\fk{p}}(f)+v_{\fk{p}}(g)$. Sabemos que K(fg) es principal si K(f) y K(g) lo son. Por el contrario, si K(fg) es princiapl generado por r, (2) de la proposición~\ref{res:4.35} muestra que r puede ser escrito como una suma $p+q$, con $p\in K(f)$ y $q\in K(g)$. Entonces para todo $\fk{p}$ de altura uno, $v_{\fk{p}}(p) \geq v_{\fk{p}}(f)$ y $v_{\fk{p}}(q) \geq v_{\fk{p}}(g)$. Por otro lado, $v_{\fk{p}}(p)+v_{\fk{p}}(q)=v_{\fk{p}}(r)=v_{\fk{p}}(fg)=v_{\fk{p}}(f)+v_{\fk{p}}(g)$. Entonces, $v_{\fk{p}}(p)=v_{\fk{p}}(f)$  $v_{\fk{p}}(q)=v_{\fk{p}}(g)$ para todo $\fk{p}$. Por (1), esto implica que K(f) y K(g) son principales. 
\end{proof}

\begin{cor}\label{res:4.38}
Sea R un dominio entero, P un monoide tórico y F una cara de P. Entonces el conjunto $\fk{F}$ de elementos de R[P] tal que K(f) es generado princiaplmente por un elemento de F es una cara del monoide multiplicativo de R[P]. 
\end{cor}

\begin{proof}
Si f y g pertenecen a $\fk{F}$, entonces $K(f)=(p)$ y $K(g)=(q)$ con $p,q \in F$, así por (3) de la proposición~\ref{res:4.36} $K(fg)=p+q$ y $p+q \in F$. Entonces $\fk{F}$ es un submonoide de $R[P]$. Por el contrario, si $fg\in \fk{F}$ entonces por (2) de la proposición~\ref{res:4.37} se tiene que K(f) y K(g) son principales, digamos generados por p y q. Entonces $p+q$ genera a $K(fg)$ y vive en F. Ya que F es una cara, p y q pertenecen a F, y entonces f y g pertenecen a $\fk{F}$. Entonces $\fk{F}$ es una cara de R[P]. 
\end{proof}

\subsection{Completaciones y series de potencias formales}

\begin{lem}\label{res:4.39}
Sea Q un monoide fino \textit{sharp} y sea $h:Q \to \bb{N}$ un homomorfismo local. Las siguiente familias de subconjuntos de Q son cofinales, es decir, dado cualquier miembro de una de estas familias, existe otro miembro de cada una de las otras familias el cual lo contiene: \\
1. $\set{J_{h,n}:n\in \bb{N}}$, donde $J_{h,n}=\set{q:h(q)>n}$; \\
2. $\set{(Q^{+})^{n}:n\in \bb{N}}$; \\
3. El conjunto de todos los subconjuntos de Q cuyo complemento es finito. 
\end{lem}

\begin{proof}
 Es claro que $(Q^{+})^{n} \subset J_{h,n-1}$ para todo n, ya que si $q=q_{1}+...+q_{n}\in (Q^{+})^{n}$ entonces $h(q)\geq n$. Por otro lado, si $(s_{1},...,s_{r})$ es un conjunto finito de generadores para $Q^{+}$ y $M:=max \set{h(s_{1}),..., h(s_{r})}$, entonces $J_{h,n} \subset (Q^{+})^{n/M}$. De hecho, cualquier $q\in Q^{+}$ puede ser escrito como $\sum m_{i}s_{i}$, con $m_{i} \in \bb{N}$, así que si $q\in J_{h,n}$,
 $$n < h(q) = \sum m_{i} h(s_{i}) \leq M \sum m_{i}$$
 y así $\sum m_{i} > n/M$ y $q\in (Q^{+})^{n/M}$. Hemos visto que el complemento de $J_{h,n}$ es finito. Por otro lado, si $\sigma$ es un conjunto finito y m es una cota para $\set{h(q):q\in \sigma}$, entonces $J_{h,m}$ está contenido en el complemento de $\sigma$. 
\end{proof}

\begin{prop}\label{res:4.40}
Sea Q un monoide \textit{sharp}, fino y R un anillo conmutativo. \\
1. Para cada $q\in Q$, $\set{(p,p^{\prime})\in Q\times Q:p+p^{\prime}=q}$ es finito. \\
2. El R-módulo topológico $R[[q]]$ definido anteriormente admite una única multiplicación continua con la propiedad de que el mapeo natural $\overline{e}:Q \to R[[Q]]$ es un homomorfismo de monoides. \\
3. El anillo topológico $R[[Q]]$ es naturalmente identificado con la completación formal de R[Q] por el ideal $R[Q^{+}]$. \\
4. Si $Q^{gp}$ es libre de torsión y R es un dominio entero, entonces R[[Q]] también es un dominio entero. \\
5. Si R es un anillo local con ideal maximal $\fk{m}$, entonces R[[Q]] es también un anillo local, cuyo ideal maximal es el conjunto de elementos de R[[Q]] cuyo termino constante pertenece a $\fk{m}$. Más aún, $Q \to R[[Q]]$ es un homomorfismo local. 
\end{prop}

\begin{proof}
Dado que Q es fino, existe un homomorfismo local $h: Q \to \bb{N}$. Entonces $J_{h,n}$ es un ideal de Q y por la proposición~\ref{res:3.34}, $Q \setminus J_{h,n}$ es un conjunto finito. Si $q=p+p^{\prime}$, entonces $p,p^{\prime}\in Q \setminus J_{h,h(q)}$, y entonces existen solo finitos $p,p^{\prime}$. Entonces si $f=\sum a_{p}e^{p}$ y $g=\sum b_{p^{\prime}}e^{p^{\prime}}$ son elementos de $R[[Q]]$, por lo que podemos definir $fg=\sum c_{q}e^{q}$ donde 
$$c_{q}:=\sum_{p+p^{\prime}=q}a_{p}b_{p^{\prime}}.$$
Para cada subconjunto finito $\sigma$ de Q,
$$U_{\sigma}:= \set{\sum a_{q}e^{q}\in R[[Q]]:a_{q}=0,q\in \sigma}$$
es abierto, y la familia de tales subconjuntos es una base de vecindades abiertas del cero. Se sigue que la operación de multiplicación definida es continua, y de hecho es la única operación continua compatible con la ley de Q. \\ 
Si K es un ideal de Q y n es un número natural, sea $K^{n}$ el conjunto de todos los elementos que se pueden escribir como sumas de n elementos de K. Entonces K es un ideal de Q y $R[K^{n}]=(R[K])^{n}$. Por el lema~\ref{res:4.39} tenemos que la topología $R[[Q^{+}]]$-ádica en R[Q] es compatible con la topología debil. \\
Es claro por la construcción que una inyección $Q \to Q^{\prime}$ induce una inyección $R[[Q]] \to R[[Q^{\prime}]]$. Si Q es fino y \textit{sharp} y $Q^{gp}$ es libre de torsión entonces puede ser incrustado en $\bb{N}^{r}$ para algún r natural. Entonces (4) es reducido al caso $R[[\bb{N}^{r}]]$, un anillo de serie de potencias formales, fácilmente visto como un dominio entero al hacer inducción sobre r. El enunciado (5) puede ser reducido al caso de un anillo de series de potencias formales, tomando una suprayección $\bb{N}^{r} \to Q$ 
\end{proof}

\begin{cor}\label{res:4.41}
Sea Q un monoide fino, sea $h:Q \to \bb{N}$ un homomorfismo, y para cada $n\in \bb{N}$ sea $J_{h,n}$ como se definió anteriormente. Entonces $J_{h,1}$ es un ideal primo $\fk{p}$ y la familia $\set{J_{h,n}:n\in \bb{N}}$ y $\set{\fk{p}^{n}:n\in \bb{N}}$ son cofinales.   
\end{cor}

\begin{proof}
 Si Q es \textit{sharp} y h es local, este corolario se sigue inmediatamente del lema anterior. En cualquier caso es claro que $\fk{p}$ es primo. Sea $F=Q \setminus \fk{p}$, así que h se factoriza a través de un homomorfismo local $\overline{h}:Q/F \to \bb{N}$. Entonces el resultado para Q y h se sigue del resultado para $Q/F$ y $\overline{h}$. 
\end{proof}

De manera más general, si S es un Q-conjunto finitamente generado, vamos a considerar el conjunto $R[[S]]$ de series de potencias formales indexada por S y con coeficientes en R, el cual va a formar un módulo sobre R[[Q]] siempre que, para todo $(q,s)\in Q^{+}\times S$, $qs\neq s$.

\begin{Def}\label{res:4.42}
Sea Q un monoide \textit{integral} y sea K un ideal de Q. El monoide de Rees de $(Q,K)$ es el monoide $B_{K}(Q)$ cuyos elementos son pares $(m,p)$ donde $m\in \bb{N}$ y $p\in K^{m}$, y cuya ley monoidal está dada por $$(m,p)+(n,q):=(m+n,p+q).$$
Si S es un Q-conjunto, el conjunto de Rees de (S,K) es el conjunto de pares (n,s) donde $n\in \bb{N}$ y $s\in K^{n}S$, con la acción de $B_{K}(Q)$ dada por 
$$(m,p)+(n,s):=(m+n,p+s).$$
\end{Def}

Notemos que tenemos 
$$h:B_{K}(Q) \to \bb{N}:(m,p) \mapsto m,~~~~ g:B_{K}(S) \to \bb{N}:(n,s) \mapsto n,$$
donde h es un homomorfismo de monoides y g es un morfismo de $B_{K}(Q)-conjuntos$.

\begin{prop}\label{res:4.43}
 Sea K un ideal de un monoide fino Q y sea S un Q-conjunto finitamente generado. Supongamos que $qs\neq s$ para cualquier $s\in S$ y $q\in Q^{+}$. Entonces los siguientes enunciados se cumplen: \\
 1. Si T es un Q-subconjunto de S, entonces existe un entero m tal que $T\cap K^{m+n}S\subset K^{n}T$ para todo n. \\
 2. Si K es un ideal propio de Q, entonces $\cap \set{K^{n}S:n\geq 0}=\emptyset$. \\
 3. Si Q es \textit{sharp}, R[[S]] puede ser identificado con la completación $R[[Q^{+}]]$-ádica de R[S]. En particular, si $(s_{1},...,s_{n})$ es una sucesión de generadores para S como un Q-conjunto, entonces $(e^{s_{1}},...,e^{s_{n}})$ es una sucesión de generadores para R[[S]] como un R[[Q]]-módulo. 
\end{prop}

\begin{proof}
 Si Q es un monoide finitamente generado, entonces es noetheriano, así que K es finitamente generado como ideal. Si $(p_{1},...,p_{r})$ es una sucesión de generadores para Q y $(k_{1},...,k_{s})$ es una sucesión de generadores para K, entonces
 $$((0,p_{1}),...,(0,p_{r}),(1,k_{1}),...,(1,k_{s}))$$
 es una sucesión de generadores para el monoide $B_{K}(Q)$, se sigue por el teorema~\ref{res:3.10} que $B_{K}(Q)$ es noetheriano. Más aún, dado que S es finitamente generado como un Q-conjunto, $B_{K}(S)$ es finitamente generado como $B_{K}(Q)$-conjunto y entonces es también noetheriano. Si $T \subset S$ es un Q-subconjunto, entonces 
 $$B_{K}(S,T)=\set{(n,t):n\in  \bb{N}, t\in T \cap K^{n}S}$$
 es naturalmente un $B_{K}(Q)$-subconjunto de $B_{K}(S)$, y por tanto es finitamente generado, digamos por $((m_{1},t_{1}),...,(m_{p},t_{p}))$. Entonces cualquier cota superior m para $(m_{1},...,m_{p})$ satisface (1). En particular, $T=\cap \set{K^{n}S:n\geq 0}$ es un Q-subconjunto de S, y (2) implica que $T=TK$. Ya que $K\subset Q^{+}$, la proposición~\ref{res:3.13} implica que $T=\emptyset$, probando (2). Ses sigue que para todo $s\in S$, existe un $n\in \bb{N}$ tal que $s\notin (Q^{+})^{n}S$, y para todo subconjunto finito $\sigma$ de S, existe un $n\in \bb{N}$ tal que $\sigma \cap (Q^{+})^{n}S= \emptyset$. Por otro lado, es fácil ver que el complemento de cada $(Q^{+})^{n}S$ es finito, Entonces la familia de subconjuntos de S cuyo complemento es finito es cofinal con la familia $\set{(Q^{+})^{n}S:n\in \bb{N}}$. Se sigue que, para cada $s\in S$, $\set{(q,t)\in Q\times S:q+t=s}$ es finito. Esto nos permite definir una acción de R[[Q]] en R[[S]] y ver que la topología producto es la topología $R[[Q^{+}]]$-ádica. 
\end{proof}

\section{Conclusiones}

El estudio de la teoría de monoides conmutativos es muy vasto y, como vimos, el considerar estructuras algebraicas más generales trae consigo grandes ventajas, como lo es el encontrar propiedades que caractericen un mayor número de estructuras. El estudio de los monoides no termina aquí, sino que es solo la base de una teoría mucho más general. Al igual que la geometría algebraica clásica —el estudio de los esquemas— comienza con el estudio de los anillos conmutativos con unidad, la teoría de monoides conmutativos da paso al estudio de las estructuras logarítmicas, como lo son los log esquemas. Estas analogías entre estructuras motivan además el desarrollo de una teoría de cohomología para dichas log estructuras. Uno de los resultados más importantes que provee la geometría algebraica logarítmica es proporcionar una generalización de la teoría clásica de variedades tóricas.

A lo largo de este trabajo sentamos las bases categóricas, algebraicas y geométricas de la teoría de monoides conmutativos: los límites y colímites en la categoría $\mathbf{Mon}$, la clase de los monoides finos y saturados, la dualidad para monoides tóricos \textit{sharp}, y la construcción de los esquemas monoidales afines junto con su relación con las variedades tóricas. Estos elementos constituyen el lenguaje sobre el que se edifica la Geometría Logarítmica de Fontaine--Kato--Illusie, cuyo estudio sistemático ---y el de la teoría de cohomología asociada a las log estructuras--- dejamos para un trabajo posterior.

\bibliographystyle{alpha}
\bibliography{referencias.bib}
\end{document}